\documentclass[12pt,a4paper,reqno]{amsart}
\usepackage{amsmath}
\usepackage{amsfonts}
\usepackage{amssymb}
\numberwithin{equation}{section}

\theoremstyle{plain}

\newtheorem{theorem}[subsection]{Theorem}
\newtheorem{proposition}[subsection]{Proposition}
\newtheorem{lemma}[subsection]{Lemma}
\newtheorem{corollary}[subsection]{Corollary}

\theoremstyle{definition}

\newtheorem{definition}[subsection]{Definition}

\newtheorem{example}[subsection]{Example}

\renewcommand{\leq}{\leqslant}
\renewcommand{\geq}{\geqslant}

\newsavebox{\proofbox}
\savebox{\proofbox}{\begin{picture}(7,7)%
  \put(0,0){\framebox(7,7){}}\end{picture}}

\newcommand{\md}[1]{\ensuremath{(\mbox{mod}\, #1)}}

\def\endproof{\hfill{\usebox{\proofbox}}}
\def\E{\mathbb{E}}
\def\Z{\mathbb{Z}}
\def\R{\mathbb{R}}

\def\C{\mathbb{C}}
\def\P{\mathbb{P}}
\def\D{\mathcal{D}}
\def\B{\mathcal{B}}
\def\Freq{\mbox{Freq}}
\def\eps{\varepsilon}
\def\Conj{\mathcal{C}}
\def\Bounded{{\mathbf b}}
\def\Phase{{\mathbf p}}

\def\Ffiven{\mathbb{F}_5^n}
\def\quarter{\textstyle \frac{1}{4} \displaystyle}
\def\half{\textstyle \frac{1}{2} \displaystyle}

\def\ni{\noindent}
\def\vs{\vspace{11pt}}

\begin{document}

\title{An inverse theorem for the Gowers $U^3(G)$ norm}

\author{Ben Green}
\address{Department of Mathematics, University of Bristol, University Walk, Bristol BS8 1TW, England
}
\email{b.j.green@bristol.ac.uk}

\author{Terence Tao}
\address{Department of Mathematics, UCLA, Los Angeles CA 90095-1555, USA.
}
\email{tao@math.ucla.edu}

\thanks{The second author is supported by a grant from the Packard Foundation.}

\begin{abstract}  There has been much recent progress in the study of arithmetic progressions in various sets, such as dense subsets of the integers or
of the primes.  One key tool in these developments has been the sequence of \emph{Gowers uniformity norms} $U^d(G)$, $d=1,2,3,\ldots$ on a finite additive group $G$; in particular, to detect arithmetic progressions of length $k$ in $G$
it is important to know under what circumstances the $U^{k-1}(G)$ norm can be large.\vs

\ni The $U^1(G)$ norm is trivial, and the $U^2(G)$ norm can be easily described in terms of the Fourier transform.  In this paper we systematically study the $U^3(G)$ norm, defined for any function $f: G \to \C$ on a finite additive group $G$ by the formula
\begin{eqnarray*} \Vert f \Vert_{U^3(G)} & := & |G|^{-4}\sum_{x,a,b,c \in G}(f(x)\overline{f(x+a)f(x+b)f(x+c)}f(x+a+b) \times \\
& & \qquad\qquad\qquad\qquad \times f(x+b+c)f(x+c+a) \overline{f(x+a+b+c)})^{1/8}.\end{eqnarray*}

\ni We give an inverse theorem for the $U^3(G)$ norm on a arbitrary group $G$.  In the finite field case $G = \Ffiven$ we show that a bounded function $f: G \to \C$ has large $U^3(G)$ norm if and only if it has a large inner product with a function $e(\phi)$, where $e(x) := e^{2\pi i x}$ and $\phi: \Ffiven \to \R/\Z$ is a quadratic phase function.  In a general $G$ the statement is more complicated -- the phase $\phi$ is quadratic only locally on a Bohr neighbourhood in $G$.\vs

\ni As an application we extend Gowers proof \cite{gowers-4-aps} of Szemer\'edi's theorem for progressions of length 4 to arbitrary abelian $G$. More precisely, writing $r_4(G)$ for the size of the largest $A \subseteq G$ which does not contain a progression of length four, we prove that 
\[ r_4(G) \ll |G| (\log \log |G|)^{-c},\]
where $c$ is an absolute constant.\vs

\ni We also discuss links between our ideas and recent results of Host-Kra and Ziegler in ergodic theory.\vs

\ni In future papers we will apply variants of our inverse theorems to obtain an asymptotic for the number of quadruples $p_1 < p_2 < p_3 < p_4 \leq N$ of primes in arithmetic progression, and to obtain significantly stronger bounds for $r_4(G)$.
\end{abstract}

\maketitle

\section{Background and Motivation}

\ni A famous and deep theorem of Szemer\'edi asserts that any set of integers of positive upper density contains arbitrarily long arithmetic progressions.
More precisely:

\begin{theorem}[Szemer\'edi's theorem, infinitary version]\label{szemeredi}\cite{szemeredi}  Let $A$ be a subset of the integers $\Z$ whose upper density $\limsup_{N \to \infty} (2N + 1)^{-1}|A \cap [-N,N]|$ is strictly positive. Then for any $k \geq 1$, the set $A$ contains infinitely many arithmetic progressions $\{ a, a+r, \ldots, a+(k-1)r\}$, $r \neq 0$, of length $k$. 
\end{theorem}

\ni The first non-trivial case of this theorem is when $k=3$, which was treated by Roth \cite{roth} using a Fourier-analytic argument.
The case of higher $k$ was more resistant to Fourier-analytic methods, and the first full proof of this theorem was achieved
by Szemer\'edi \cite{szemeredi} using combinatorial methods.  Later, Furstenberg \cite{furst,FKO}
introduced an ergodic theoretic proof of this theorem.  More recently, Gowers \cite{gowers-long-aps} gave a proof which was both combinatorial and Fourier-analytic in nature, and which is substantially closer in spirit to Roth's original argument than the other proofs.  Even more recently there have been a number of other proofs of this theorem by other methods, such as hypergraph regularity \cite{gowers-hyper,nagle-rodl-schacht1,nagle-rodl-schacht2,rodl,rodl2,rodl3} or ``discrete ergodic theory'' \cite{tao:ergodic}.
This theorem and its various proofs have in turn generated many other mathematical developments.  For instance, in \cite{green-tao-primes}, we were able to apply Theorem \ref{szemeredi} to demonstrate that the primes contain arbitrarily long arithmetic progressions.\vs

\ni In this paper we shall be interested primarily in the Fourier-analytic approach to this theorem, specifically in the $k=4$ case, which was treated separately by Gowers in \cite{gowers-4-aps} and then again in \cite{gowers-long-aps}. This latter paper will be our key reference.  However as we shall see later there are some strong connections between this approach and the ergodic one, especially after the work on characteristic factors by Host and Kra \cite{host-kra1,host-kra2} and Ziegler \cite{ziegler02,ziegler}, and on the connection to nilsequences by Bergelson, Host, and Kra \cite{bhk}.  Before we give our main new results, however, we first give some further historical background and motivation.\vs

\ni Gowers' proof of the full Szemer\'edi theorem in \cite{gowers-long-aps} is quite lengthy and involves many deep new ideas.  However, it is 
possible to split it up into a number of simpler steps, all but one of which are straightforward.  Firstly, it is easy to show that for any 
fixed $k$, Theorem \ref{szemeredi} is equivalent to the following finitary version.  

\begin{theorem}[Szemer\'edi's theorem, finitary version \cite{szemeredi}]\label{szemeredi-finite}  Let $\delta > 0$ and $k \geq 1$.  Then there exists an integer $N_0 = N_0(\delta,k)$ such that whenever $N \geq N_0$ and $A \subseteq  [1,N]$ is such that $|A| / |[1,N]| \geq \delta$, then $A$ contains at least one proper arithmetic progression of length $k$.
\end{theorem}

\ni The next observation, due to Roth, is that one can hope to prove this theorem by downwardly inducting on the density parameter $\delta$ 
(the case $\delta \geq 1$ being trivial or vacuous).  In particular, for any fixed $k$, Theorem \ref{szemeredi-finite} is equivalent to the 
following assertion.

\begin{theorem}[Lack of progressions implies density increment]\label{prog-dens}  Let $\delta > 0$ and $k \geq 1$.  Let $N \geq 1$, and let $A \subseteq  [1,N]$
be such that $|A| / |[1,N]| \geq \delta$, and such that $A$ contains no proper arithmetic progressions of length $k$.  Then, if $N$ is sufficiently large depending on $k$ and $\delta$,
there exists an arithmetic progression $P \subseteq [1,N]$ with $|P| \geq \omega(N,\delta)$ for some function $\omega(N,\delta)$ of $N$ which goes 
to infinity as $N \to \infty$ for each fixed $\delta$, such that we have the density increment
$|A \cap P| / |P| \geq \delta + c(\delta)$, where $c(\delta) > 0$ is a 
function of $\delta$ which is bounded away from zero whenever $\delta$ is bounded away from zero.
\end{theorem}

\ni The deduction of Theorem \ref{szemeredi-finite} from Theorem \ref{prog-dens} is a straightforward induction argument.
For details of arguments of this type any of \cite{gowers-4-aps,gowers-icm,gowers-long-aps,green-fin-field,roth} may be consulted, or indeed \S \ref{finite-sz} or \ref{sz-g} of this paper.  Of course, the final bound $N_0(\delta,k)$ obtained in
Theorem \ref{szemeredi-finite} will depend on the explicit bounds $\omega(N,\delta)$, $c(\delta)$ obtained in Theorem \ref{prog-dens}.
In Roth's $k=3$ argument in \cite{roth}, $c(\delta)$ was roughly $\delta^2$ and $\omega(N,\delta)$ was roughly $N^{1/2}$, which led
to a final bound of the form $N_0(\delta,3) \leq \exp(\exp(C/\delta))$.  In Gowers' extension of Roth's argument in \cite{gowers-long-aps},
$c(\delta)$ was roughly $\delta^{C_k}$ and $\omega(N,\delta)$ was roughly $N^{c_k}$ for some $C_k,c_k > 0$ depending only on $k$, which 
led to a final bound of the form $N_0(\delta,k) \leq \exp(\exp(C_k/\delta^{C_k}))$ (see \cite{gowers-long-aps} for a more precise
statement).  These are the best known bounds for $N_0(\delta,k)$ except in the $k=3$ case, where the current record is
$N_0(\delta,3) \leq (C/\delta)^{C/\delta^2}$, due to Bourgain \cite{Bou}.\vs

\ni The next step is to pass from the interval $[1,N]$ to a cyclic group $\Z/N\Z$ for some prime $N$.  Indeed, by using Bertrand's postulate\footnote{that is, there is always a prime between $X$ and $2X$} and a simple covering argument to split progressions in $\Z/N\Z$ into progressions in $[1,N]$, one can show that Theorem \ref{prog-dens} is in turn equivalent for each fixed $k$ (up to minor changes in the bounds $\omega(N,\delta)$ and $c(\delta)$) to the following statement.

\begin{theorem}[Lack of progressions implies density increment]\label{prog-dens-cyclic}  
Let $\delta > 0$ and $k \geq 1$.  Let $N \geq 1$ be a prime, and let $Q$ be a proper progression in $\Z/N\Z$ such that
$|Q| \geq c_0 N$ for some $0 < c_0 \leq 1$.  Let  $A \subseteq  Q$
be such that $|A|  \geq \delta N$, and such that $A$ contains no proper arithmetic progressions of length $k$.  Then, 
if $N$ is sufficiently large depending on $k$ and $\delta$, there exists a proper arithmetic progression $P \subseteq \Z/N\Z$ with 
$|P| \geq \omega(N,\delta,c_0,k)$ for some function $\omega(N,\delta,c_0,k)$ of $N$ which goes to infinity as $N \to \infty$ for each 
fixed $\delta$, $c_0$, $k$, such that we have the density increment
\[ \frac{|A \cap P|}{|P|} \geq \frac{|A \cap Q|}{|Q|} + c(\delta,c_0,k),\] where $c(\delta,c_0,k) > 0$ is 
bounded away from zero whenever $\delta$, $c_0$ are bounded away from zero and $k$ is fixed.
\end{theorem}

\ni The deduction of Theorem \ref{prog-dens} from Theorem \ref{prog-dens-cyclic} is not difficult, 
see \cite{gowers-4-aps,gowers-long-aps,roth}.  Of course, it remains to prove Theorem \ref{prog-dens-cyclic}.  This was achieved 
in the $k=3$ case by Roth using Fourier-analytic methods.  To extend these arguments to the case of higher $k$, Gowers introduced a collection of tools which form a part of a theory which might be termed ``higher-order Fourier analysis'' for reasons which will become clear later.  In particular, to handle the $k=4$ case required ``quadratic Fourier analysis''. \vs 

\ni While Gowers' original argument takes place in a cyclic group $\Z/N\Z$ of prime order, we will work in the more general setting of an arbitrary
finite additive group.  This might seem unnecessary, but is consistent with what we call the \textit{finite field philosophy}. This is the observation that many questions concerning the integers $\{1,\dots,N\}$ or the cyclic group $\Z/N\Z$ may be asked very naturally for an arbitrary finite abelian group $G$, and they may be answered there by modifying the proof for $\Z/N\Z$ in a straightforward way. Thus it is often the case that the passage
\[ \Z/N\Z \xrightarrow{\mbox{\scriptsize generalization}} G\]
is rather straightforward. \vs

\ni However, it may be that the question is significantly easier to answer when $G$ is some specific group, typically a vector space over a finite field such as $\mathbb{F}_2^n, \mathbb{F}_3^n$ or $\mathbb{F}_5^n$. This observation was made in such papers as \cite{mesh-ap, ruzsa3}. In this paper we will take a particular interest in $\mathbb{F}_5^n$ since this is the smallest characteristic field for which arithmetic progressions of length 4 are a sensible thing to discuss. Now the passage
\[ \mathbb{F}_5^n \xrightarrow{\mbox{\scriptsize generalization}} G\]
might not be at all easy. However, in attempting such a route one has split the problem into two presumably easier subproblems, and furthermore there is now a library of tools available for effecting the generalization. This started with the work of Bourgain \cite{Bou} (though he did not phrase it this way), and has continued with various works such as \cite{green1,green-konyagin}. The present paper, particularly \S \ref{sec8} and \S \ref{ggc}, is another example in this vein. For a longer discussion of the finite field philosophy, see \cite{green-fin-field}.

\begin{definition}[Additive groups]
Define an \emph{additive group} to be a group $G = (G,+)$ with a commutative group operation $+$; if $x \in G$ and $n \in \Z$ we can define
the product $nx \in G$ in the usual manner.  If $f: G \to H$ is a function from one additive group to another, and $h \in H$, we define the shift\footnote{This ``ergodic'' notation corresponds to the backwards shift $T^h x := x-h$ on the underlying group $G$.  We will discuss further connections with ergodic theory in \S \ref{ergodic-sec}.}
 operator $T^h$ applied to $f$ by the formula $T^h f(x) := f(x+h)$, and the difference operator $h \cdot \nabla := T^h - 1$ applied to $f$ by the formula
$(h \cdot \nabla) f(x) := f(x+h) - f(x)$.  We extend these definitions to functions of several variables by subscripting the variable to which the operator is applied, thus for instance if $f(x,y)$ is a function of two variables we define $T^h_x f(x,y) = f(x+h,y)$ and $h \cdot \nabla_x f(x,y) = f(x+h,y) - f(x,y)$ if $x, h$ range inside an additive group $G$, and similarly
for the $y$ variable.
\end{definition} 

\ni\textit{Remark.} Throughout the paper, we will write $N := |G|$ for the cardinality of $G$.\vs

\ni\textit{Remark.} 
The notation above is of course designed to mimic that of several variable calculus.  We caution however that we do not assign any independent meaning to the symbol $\nabla$, unless it is prepended with a shift $h$ to create a difference operator $h \cdot \nabla$, which is of course a discrete analogue of a directional derivative operator.\vs 

\ni We now introduce a multilinear form $\Lambda_k(f_1,\ldots,f_k)$ which is useful for 
counting arithmetic progressions.  Here, and throughout the paper, it is convenient to adopt the notation of conditional expectation, which allows one to hide some distracting normalizing factors such as $1/N$ in our arguments. Thus if $f: G \to \C$ is a complex-valued function on a finite set $G$, and $B \subseteq  G$ is a non-empty subset of $G$, we will use $\E_{x \in B} f(x) := \frac{1}{|B|} \sum_{x \in B} f(x)$
to denote the average of $f$ over $B$.  We will abbreviate $\E_{x \in G} f(x)$ as $\E(f)$ when the domain $G$ of $f$ is clear from context.\vs

\ni Now if $G$ is a finite additive group and $f_0,\ldots,f_{k-1}: G \to \C$ are complex-valued functions, 
we define the $k$-linear form $\Lambda_k(f_0,\ldots,f_{k-1}) \in \C$ by
$$ \Lambda_k(f_0,\ldots,f_{k-1}) := \E_{x,r \in G} f_0(x) T^r f_1(x) \ldots T^{(k-1)r} f_{k-1}(x).$$
Observe that if $A \subseteq  G$ and $f_0 = \ldots = f_{k-1} = 1_A$, where $1_A: G \to \{0,1\}$ denotes the indicator function of $A$, then $\Lambda_k(1_A,\ldots,1_A)$ is just the number of progressions of length $k$ (including those with common difference $0$), divided by the normalizing factor of $N^2$.  In particular, if $(N,(k-1)!) = 1$ and $A$ contains no proper progressions of length $k$ then we see that $\Lambda_k(1_A, \ldots, 1_A) = |A|/N^2$, which will be quite small when $N$ is large.\vs

\ni It is thus of interest to determine under what conditions $\Lambda_k(1_A, \ldots, 1_A)$ is small or large.  To this end, Gowers introduced (what are now known as) the \emph{Gowers uniformity norms} $\|f\|_{U^d(G)}$ for any complex function $f: G \to \C$, whose definition we now recall.

\begin{definition}[Gowers uniformity norm]\label{gundef}  Let $d \geq 0$, and let $f: G \to \C$ be a function.  We define the Gowers uniformity norm $\|f\|_{U^d(G)} \geq 0$ of $f$ to be the quantity
\[ \Vert f \Vert_{U^d(G)} := \big(\E_{x \in G, h \in G^d} 
\prod_{\omega \in \{0,1\}^d} {\Conj}^{|\omega|} T^{\omega \cdot h} f(x)\big)^{1/2^d},\]
where $\omega = (\omega_1,\ldots,\omega_d)$, $h \in (h_1,\ldots,h_d)$, $\omega \cdot h := \omega_1 h_1 + \ldots + \omega_d h_d$,
$|\omega| := \omega_1 + \ldots + \omega_d$, and ${\Conj}$ is the conjugation operator ${\Conj} f(x) := \overline{f(x)}$.
\end{definition}
\ni\textit{Remark.} An equivalent definition of the $U^d(G)$ norms is given by the recursive formulae
\begin{equation}\label{ud-recursive}
\|f\|_{U^0(G)} = \E(f); \quad \|f\|_{U^1(G)} = |\E(f)|; \quad \| f \|_{U^d} := (\E_{h \in G} \| T^h f \overline{f} \|_{U^{d-1}(G)}^{2^{d-1}})^{1/2^d}
\end{equation}
for all $d \geq 1$.\vs

\ni\textit{Remark.} A configuration of the form $(x + \omega \cdot h)_{\omega \in \{0,1\}^d}$ is called a \emph{cube of dimension $d$}.  Thus $\|f\|_{U^d(G)}^{2^d}$ is a weighted average of $f$ over cubes; for instance, $\|1_A\|_{U^d(G)}^{2^d}$ is equal to the number of cubes contained in $A$, divided by the normalizing
factor of $N^{d+1}$.  The cases $d=0,1$ are rather degenerate, and indeed $U^d(G)$ is not a norm in these cases.
However for $d > 1$, one can show that $\| \cdot \|_{U^d(G)}$ is indeed a norm, i.e. it is homogeneous, non-negative, non-degenerate, and obeys the triangle inequality, see \cite[Lemma 3.9]{gowers-long-aps}.  These norms have also appeared recently in ergodic theory, see for instance
\cite{host-kra2}, and (together with the dual norms $U^d(G)^*$) played a key role in \cite{green-tao-primes}.
It thus seems of interest to study these norms more systematically; the results here can be viewed as a step in that direction.\vs

\ni We will study these norms in detail later, but for now let us give an example to illustrate what they are trying to capture.  Suppose
$f$ has the form $f(x) := e(\phi(x))$ for some phase function $\phi: G \to \R/\Z$, where $e: \R/\Z \to \C$ is the exponential map $e(x) := e^{2\pi i x}$.  Then a simple calculation shows that
$$ \| f \|_{U^d}^{2^d} = \E_{x,h_1,\ldots,h_d \in G}  e((h_1 \cdot \nabla_x) \ldots (h_d \cdot \nabla_x) \phi(x)).$$
Thus the $U^d$ norm is in some sense measuring the oscillation present in the $d^{th}$ ``derivative'' of the phase.  In particular, we expect the $U^d$ norm to be large if the phase behaves like a ``polynomial'' of degree $d-1$ or less, but small if the phase is behaving like a polynomial of degree $d$ or higher. \vs

\ni We observe that as an immediate consequence of \eqref{ud-recursive} and induction we have the monotonicity property
\begin{equation}\label{ud-monotone}
\| f \|_{U^d(G)} \leq \|f\|_{U^{d+1}(G)} \hbox{ for } d = 0,1,2,\ldots.
\end{equation}

\ni The relevance of the Gowers uniformity norms to arithmetic progressions lies in the following result, which was stated explicitly in \cite[Theorem 3.2]{gowers-long-aps} (in the case of cyclic groups $G = \Z/N\Z$) 
but has been implicit in the ergodic theory literature for some time. Write
$\mathcal{D} := \{z \in \mathbb{C} : |z| \leq 1\}$ for the unit disk.

\begin{proposition}[Generalized von Neumann Theorem]\label{gvn}  Let $G$ be a finite abelian group with $(N,(k-1)!) = 1$.
Let $f_0,\dots,f_{k-1} : G \rightarrow \mathcal{D}$ be functions.  Then we have
\[ |\Lambda_k(f_0,\dots,f_{k-1})| \leq \min_{1 \leq j \leq k} \Vert f_j \Vert_{U^{k-1}(G)}. 
\]
\end{proposition}

\ni It is instructive to continue with the phase example given earlier.  If $f_j = e(\phi_j)$, then 
$$\Lambda_k(f_0,\ldots,f_{k-1}) = \E_{x,r \in G} e( \phi_0(x) + \phi_1(x+r) + \ldots + \phi_{k-1}(x+(k-1)r)).$$
Thus $\Lambda_k(f_0,\ldots,f_{k-1})$
is measuring the oscillation present in the expression $\phi_0(x) + \phi_1(x+r) + \ldots + \phi_{k-1}(x+(k-1)r)$.  Proposition \ref{gvn} can 
then be viewed as a statement that if this expression does not oscillate, then neither do the expressions $(h_1 \cdot \nabla_x)
\ldots (h_{k-1} \cdot \nabla_x) \phi_j(x)$
for any $1 \leq j \leq k$.  Note that such a fact morally follows by ``differentiating'' the expression $\phi_0(x) + \phi_1(x+r) + \ldots + \phi_k(x+(k-1)r)$ in $k-1$ different directions to eliminate all but one of the terms in this series.  For completeness we give a proof of
this Proposition in Section \ref{Avg-sec}. 

\begin{corollary}[Lack of progressions implies large uniformity norm \cite{gowers-long-aps}]\label{prog-unif}  Let $k \geq 3$, let $G$ be a finite additive group with $(N,(k-1)!) = 1$, and let $A \subseteq  G$, $|A| = \alpha N$, be a non-empty set such that $A$ has no proper arithmetic progressions of length $k$.
If $N \geq 2/ \alpha^{k-1}$, then we have $\| 1_A - \alpha \|_{U^{k-1}(G)} \geq 2^{-k-1} \alpha^{k-1}$.\vs

\ni More generally, let $P$ be a proper arithmetic progression in $G$ such that $|P| \geq c_0 N$.  Let $A \subseteq  P$, $|A| = \alpha |P|$, be a non-empty set which contains no proper arithmetic progressions.  If $N > N_0(c_0,k,\alpha)$ then we have
$\| 1_A - \alpha 1_{P} \|_{U^{k-1}(G)} \geq c(c_0, \alpha,k) > 0$, where the quantity $c(c_0, \alpha, k)$ stays bounded away from zero when $c_0, \alpha$ are bounded away from zero and $k$ is fixed.
\end{corollary}

\begin{proof}  We begin with the first claim.
Since $1_A = \alpha + (1_A-\alpha)$, we can split the expression $\Lambda_k(1_A,\ldots,1_A)$ as the sum of $2^k$ expressions,
one of which is $\Lambda_k(\alpha,\ldots,\alpha)$, and the other $2^{k-1}$ of which can be bounded in magnitude by $\|1_A - \alpha \|_{U^{k-1}(G)}$
thanks to Proposition \ref{gvn}.  In particular we conclude that
$$ |\Lambda_k(\alpha,\ldots,\alpha) - \Lambda_k(1_A,\ldots,1_A)| \leq 2^k \| 1_A - \alpha \|_{U^{k-1}(G)}.$$
But clearly $\Lambda_k(\alpha,\ldots,\alpha) = \alpha^k$, while since $A$ has no proper arithmetic progressions we have
$\Lambda_k(1_A,\ldots,1_A) = |A|/N^2 = \alpha/N \leq \alpha^k/2$.  The first claim follows.\vs

\ni The second claim proceeds similarly but is based upon the decomposition $1_A = \alpha 1_{P}
+ (1_A - \alpha 1_{P})$, and the observation that $\Lambda_k(1_{P},\ldots,1_{P}) \geq c(c_0,k) > 0$ for some
positive quantity $c(c_0,k)$ depending on $c_0$ and $k$; we leave the details to the reader.
\end{proof}\vs

\ni Comparing this Proposition with Theorem \ref{prog-dens-cyclic}, we thus see that in order to prove Szemer\'edi's theorem for a fixed $k$, it
suffices to prove the following:

\begin{theorem}[Large uniformity norm implies density increment \cite{gowers-long-aps}]\label{large-density}
Let $\eta > 0$ and $k \geq 3$.  Let $G = \Z/N\Z$ be a cyclic group of prime order, and let $f: G \to \D$ be a real-valued bounded function
such that $\E(f) = 0$ and $\|f\|_{U^{k-1}(G)} \geq \eta$.  Then, 
if $N > N_0(k,\eta)$, there exists a proper arithmetic progression $P \subseteq G$ with 
$|P| \geq \omega(N,\eta,k)$ such that $\E_{x \in P}  f(x) \geq c(\eta,k)$, where
\begin{itemize}
\item $\omega(N,\eta) \rightarrow \infty$ as $N \rightarrow \infty$ for fixed $\eta$;
\item $c(\eta,k) > 0$ is bounded away from zero when $\eta$ is bounded away from zero and $k$ is fixed.
\end{itemize}
\end{theorem}

\ni Indeed, Theorem \ref{prog-dens-cyclic} then follows by applying Corollary \ref{prog-unif} and then invoking Theorem \ref{large-density} with
$f := 1_A - (\E_{x \in Q} 1_A(x)) 1_Q$.  Theorem \ref{large-density} is in fact deduced in \cite{gowers-long-aps} from the following stronger theorem:

\begin{theorem}[Weak inverse theorem for $U^{k-1}(\Z/N\Z)$ \cite{gowers-long-aps}]\label{weak-inverse}
Let $\eta > 0$ and $k \geq 3$.  Let $G = \Z/N\Z$ be a cyclic group of prime order, and let $f: G \to \D$ be a bounded function
such that $\|f\|_{U^{k-1}(G)} \geq \eta$.  Then, if $N \geq \exp(C_k \eta^{-C_k})$ for some sufficiently large $C_k > 0$, one can partition $\Z/N\Z$ into arithmetic progressions
$(P_j)_{j \in J}$, each of size
$|P_j| \geq c_k \eta^{-C_k} N^{c_k \eta^{-C_k}}$ for some $c_k, C_k > 0$, such that
$$\sum_{j \in J} |\E( f  1_{P_j} )| \geq c_k \eta^{C_k}$$
for some $c_k, C_k > 0$.
\end{theorem}

\ni Theorem \ref{large-density} (and hence Theorem \ref{szemeredi}) follows quickly from this and the
mean zero hypothesis  
$\sum_{j \in J} \E( f  1_{P_j} )  = \E(f) = 0$.  Indeed one gets a fairly good quantitative result for Theorem \ref{szemeredi-finite}, with
$N_0 = \exp(\exp( C_k \delta^{-C_k} ))$ for some explicit $C_k > 0$; see \cite{gowers-long-aps}.\vs

\ni We refer to Theorem \ref{weak-inverse} as a \textit{weak} inverse theorem because it gives a necessary criterion in order for a bounded function $f$
to have large $U^{k-1}(G)$ norm, and hence a sufficient condition for the $U^{k-1}(G)$ norm to be
small.  As discussed above, this theorem is strong enough to imply Szemer\'edi's theorem.  Also, Theorem \ref{weak-inverse}
could potentially be useful, when combined with such tools as Theorem \ref{gvn}, for not only demonstrating the existence of progressions of length $k$
in a given set $A$, but in fact providing an accurate count as to \emph{how many} such progressions there are.  For instance, one might hope to count the number of progressions of length $k$ in the primes less than $N$ by using Theorem \ref{weak-inverse} to show that a certain counting function $f$ associated to the primes has small $U^{k-1}(\Z/N\Z)$ norm and hence its contribution to the count of progressions in the primes could be controlled using Theorem \ref{gvn}.  However, the sufficient condition for smallness of $U^{k-1}(\Z/N\Z)$ given by Theorem \ref{weak-inverse}
is very difficult to verify for sets such as the primes (being at least as difficult as the Elliott-Halberstam conjecture, which is not known to be implied even by the GRH). \vs

\ni It is thus of interest to obtain a better inverse theorem for the $U^{k-1}(\Z/N\Z)$ norm, which gives a more easily checkable condition for
when this norm is small.  Ideally we would like this condition to be both necessary and sufficient, at least up to constant losses.  In this paper we shall achieve these objectives for $k=4$. \vs

\ni In subsequent work we will give various applications of the results and methods of this paper. In \cite{green-tao-ffszem} we obtain a new bound on the size of the largest subset of the vector space $\Ffiven$ with no 4-term arithmetic progression, and we hope to generalize that result to arbitrary abelian groups $G$. In another series of papers, we will obtain an asymptotic formula for the number of quadruples $p_1 < p_2 < p_3 < p_4 \leq N$ of primes in arithmetic progression. 

\section{Inverse theorems for $U^{k-1}$ norms}

\ni We have now motivated why we are interested in an inverse theorem for the $U^{k-1}$ norms.  Before we state our main theorems, let us give some other examples and results which will illustrate what the inverse theorem should be.  Recall that the $U^d$ norm of a function $e(\phi)$ measures
the oscillation in the $d^{th}$ derivative of the phase $\phi$.  Also recall that a polynomial of degree at most $d-1$ is a function whose 
$d^{th}$ derivative vanishes.  We generalize this concept as follows.

\begin{definition}[Locally polynomial phase functions]
If $B$ is any non-empty subset of a finite additive group $G$ and $d \geq 1$, we say
that a function $\phi: B \to \R/\Z$ is a \emph{polynomial phase function of order at most $d-1$ locally on $B$} if we have 
$$ (h_1 \cdot \nabla_x) \ldots (h_{d+1} \cdot \nabla_x) \phi(x) = 0$$
whenever the cube $(x + \omega_1 h_1 + \ldots + \omega_{d} h_{d})_{\omega_1,\ldots,\omega_{d} \in \{0,1\}}$ is contained in $B$.
If $f: B \to \C$ is a function, we define the \emph{local polynomial bias of order $d$ on $B$} $\| f \|_{u^d(B)}$ to be the quantity
$$ \| f \|_{u^d(B)} := \sup |\E_{x \in B}( f(x) e(-\phi(x)) )|$$
where $\phi$ ranges over all local polynomial phase functions of order at most $d-1$ on $B$.
\end{definition}

\ni To begin with we will work in the global setting $B=G$, but as will become clear later we will need to also work in the local setting.
We will refer to polynomial phase functions of degree at most 1 as \emph{linear phase functions}, and of degree at most 2 as \emph{quadratic phase functions}, with the modifiers ``local'' or ``global'' as appropriate.\vs

\ni The quantity $\| \|_{u^d(B)}$ is clearly a seminorm.  It shares several features in common with the $U^d(G)$ norm.  First of all, like the $U^d(G)$
norm, we have the monotonicity $\|f\|_{u^d(B)} \leq \|f\|_{u^{d+1}(B)}$, and when $B = G$ we also have the shift invariance
$\|T^h f \|_{u^d(G)} = \|f\|_{u^d(G)}$.  We also have the conjugation symmetry $\|\overline{f}\|_{u^d(B)} = \|f\|_{u^d(B)}$,
and the phase invariance $\| f e(\phi)\|_{u^d(B)} = \|f\|_{u^d(B)}$ whenever $\phi$ is a locally polynomial phase of degree at most $d-1$ on $B$.
The latter invariance also extends to the $U^d(G)$ norm, thus
\begin{equation}\label{phase-invariance}
\| f e(\phi) \|_{U^d(G)} = \|f\|_{U^d(G)}
\end{equation}
whenever $\phi: G \to \R/\Z$ is a global polynomial phase function of degree at most $d-1$.  Indeed,
this invariance\footnote{This polynomial phase invariance also indicates why Fourier analysis - which is essentially invariant under modulation by linear phase functions but not by quadratic or higher phases - is only able to effectively deal with the $U^2(G)$ norm and not with higher norms.  To deal with the $U^3(G)$ norm thus requires some sort of ``quadratic Fourier analysis'' which is insensitive to phase modulations by quadratic phases.  The results here can be viewed as some preliminary steps towards establishing such a quadratic Fourier analysis theory.} 
can easily be seen from \eqref{ud-recursive} and induction, using the fact that the
derivative of a polynomial of degree at most $d-1$ is a polynomial of degree $d-2$.\vs

\ni From this invariance and \eqref{ud-recursive}, \eqref{ud-monotone} we conclude that
$$ \|f\|_{U^d(G)} = \| f e(-\phi) \|_{U^d(G)} \geq \| f e(-\phi) \|_{U^1(G)}
= |\E_{x \in G}( f(x) e(-\phi(x)) )|$$
whenever $\phi$ is a global polynomial phase of degree at most $d-1$.
Taking suprema over all $\phi$, we obtain the inequality
\begin{equation}\label{ud-inverse-easy}
\|f\|_{U^d(G)} \geq \|f\|_{u^d(G)}
\end{equation}
for all $d \geq 1$, all additive groups $G$, and all $f: G \to \C$.\vs

\ni It is now natural to ask whether the inequality \eqref{ud-inverse-easy} can be reversed.  
When $d=1$ it is easy to verify (using \eqref{ud-recursive} and the fact that polynomials of degree at most 0 are constant) that we in fact have equality:
$$
\|f\|_{U^1(G)} = \|f\|_{u^1(G)}.
$$

\ni Consider next the case $d = 2$.  For this we need the Fourier transform.
Let $\widehat G$ be the Pontryagin dual of $G$, in other words the space of 
homomorphisms $\xi: x \mapsto \xi \cdot x$ from $G$ to $\R/\Z$.  As is well known, $\widehat G$ is an additive group which is isomorphic to $G$.
If $\xi \in \widehat G$, we define the \emph{Fourier coefficient} $\widehat f(\xi)$ of $f$ at the frequency $\xi$ by the formula
$$ \widehat f(\xi) = \E_x f(x) e(-\xi \cdot x).$$
As is well known, we have the
\emph{Fourier inversion formula}
$$ f(x) = \sum_{\xi \in \widehat G} \widehat f(\xi) e(\xi \cdot x)$$
and the \emph{Plancherel identity}
\begin{equation}\label{plancherel}
\E( |f|^2 ) = \sum_{\xi \in \widehat G} |\widehat f(\xi)|^2.
\end{equation}
One can then easily verify the pleasant identity
\begin{equation}\label{u24}
\| f \|_{U^2(G)}^4 = \sum_{\xi \in \widehat G} |\widehat f(\xi)|^4.
\end{equation}
For instance, this can be achieved by first establishing the identity $\| f \|_{U^2}^4 = \E( |f*f|^2 )^2$, where $f*f(x) := \E_y f(y) f(x-y) $
is the convolution of $f$ with itself, and then using \eqref{plancherel}.  Next, we make the easy observation that if $\phi: G \to \R/\Z$ is a global polynomial phase function  
of degree at most 1, then $x \mapsto \phi(x) - \phi(0)$ is a homomorphism from $G$ to $\R/\Z$, and hence there exists $\xi \in \widehat G$
such that $\phi(x) = \xi \cdot x + \phi(0)$.  From this it is easy to see that
\begin{equation}\label{u2-infty}
\| f \|_{u^2(G)} = \sup_{\xi \in \widehat G} |\widehat f(\xi)|.
\end{equation}

\ni Combining \eqref{plancherel}, \eqref{u24}, \eqref{u2-infty} we readily conclude

\begin{proposition}[Inverse theorem for $U^2(G)$ norm]\label{u2-inverse}  Let $f: G \to \D$ be a bounded function.  Then
$$ \| f\|_{u^2(G)} \leq \|f\|_{U^2(G)} \leq \|f\|_{u^2(G)}^{1/2}.$$
\end{proposition}

\ni We remark that this Proposition easily implies the $k=3$ case of Theorem \ref{weak-inverse} (using Dirichlet's theorem on approximation by rationals to cover $G$ by progressions on which $\phi$ is close to constant), and hence also implies Szemer\'edi's theorem for $k=3$.  Indeed
this is essentially Roth's original argument \cite{roth}, albeit phrased in very modern language.\vs

\ni Based on evidence such as Proposition \ref{u2-inverse}, one is tempted to conjecture that the $U^d(G)$ and $u^d(G)$ norms are also related for higher $d$, in the sense that
if $f$ is bounded and one of the two norms $\|f\|_{U^d(G)}$, $\|f\|_{u^d(G)}$ is small, then the other is also.  From \eqref{ud-inverse-easy}
we already know that one direction is true: smallness of the $U^d(G)$ norm implies the smallness of the $u^d(G)$ norm.  
Our first main result establishes a converse to this in the $d=3$ case when $G = \Ffiven$, though with only partially satisfactory control on the constants.

\begin{theorem}[Inverse theorem for $U^3(\Ffiven)$]\label{mainthm1}  
Let $f: \Ffiven \to \D$ be a bounded function and let $0 < \eta \leq 1$.  
\begin{enumerate}

\item If $\|f\|_{U^3(\Ffiven)} \geq \eta$, then there exists a subspace $W \leq \Ffiven$ of codimension at most $(2/\eta)^C$ such that
\begin{equation}\label{efphi}
\E_{y \in \Ffiven} \| f \|_{u^3(y+W)} \geq (\eta/2)^C,
\end{equation}
where we can take $C = 2^{16}$.
In particular, there exists $y \in G$ such that $\| f \|_{u^3(y+W)} \geq (\eta/2)^C$.

\item Conversely, given any subspace $W \leq \Ffiven$ and any function $f: \Ffiven \to \C$ we have
$\|f\|_{U^3(\Ffiven)} \geq \|f\|_{u^3(\Ffiven)} \geq 5^{-n}|W| \| f \|_{u^3(y+W)}$
for any $y \in \Ffiven$.

\end{enumerate}
\end{theorem}

\ni Combining the two parts of the theorem together we see that
\begin{equation}\label{fF}
 \|f\|_{u^3(\Ffiven)} \leq \|f\|_{U^3(\Ffiven)} \leq \frac{C}{\log^c (1 + 1 / \| f \|_{u^3(\Ffiven)})}
\end{equation}
for some absolute constants $c, C > 0$.  Thus this does give a result which asserts that the smallness of the $U^3(\Ffiven)$ norm implies the smallness of the $u^3(\Ffiven)$ norm and vice versa, although the dependence of constants is poor\footnote{We conjecture that one can improve the upper bound in \eqref{fF} to a polynomial dependence (bringing this estimate in line with Proposition
\ref{u2-inverse}; see \S \ref{sec13} for further discussion.}.  Note however that the control is much better if one localizes
the quadratic bias norm $u^3(\Ffiven)$ to cosets $y+W$ of $W$.  We remark that there is nothing particularly special about the finite field $\Ffiven$, and one has similar results for any other finite field of odd characteristic, though the constants depend of course on the field.\vs

\ni We shall prove Theorem \ref{mainthm1} (i) in \S \ref{finite-sec} ((ii) is easier, and we will prove it in \S \ref{model-sec}); it contains many of the main ideas of this paper, which combine the Fourier and combinatorial analysis of Gowers in \cite{gowers-4-aps} with an additional ``symmetry argument'' which is necessary to obtain a strong inverse theorem instead of a weak inverse theorem.  As a consequence we obtain, in \S \ref{finite-sz}, a Szemer\'edi theorem for progressions of length 4 in $\Ffiven$.\vs

\ni Let us now discuss finite abelian groups $G$ in general, particular importance being attached to $\Z/N\Z$ on account of potential applications.  
It is tempting to conjecture, in light of the preceding results, that in such groups, any bounded function with small $u^d(G)$
norm must necessarily have small $U^d(G)$ norm.  Unfortunately, such a statement is false, 
even for $G = \Z/N\Z$.  This fact was essentially discovered by Furstenberg and Weiss \cite{furstenberg-weiss}, in the closely related context of determining characteristic factors for multiple recurrence in ergodic theory; a similar observation was also made in page 487 of
\cite{gowers-long-aps}.  See \S \ref{ergodic-sec} for some further discussion 
of this connection.  We give one instance of the Furstenberg-Weiss example as follows:

\begin{example}\label{fw-ex}  Let $N$ be a large prime number, and let $M$ be the largest integer less than $\sqrt{N}$.  Let $G := \Z/N\Z$, and let $f: G \to \C$ be the bounded function defined by setting
$f( yM + z ) := e(y z / M ) \psi(y/M) \psi(z/M)$ whenever $-M/10 \leq y, z \leq M/10$, and $f = 0$ otherwise; here $\psi: \R \to \R_{\geq 0}$ is a non-negative smooth cutoff function which equals one on the interval $[-1/20,1/20]$ and vanishes outside of $[-1/10,1/10]$.  Then a direct calculation shows that
$\|f\|_{U^3(G)} \geq c_0$ for some absolute constant $c_0 > 0$, basically because all the phases in the expression for $\|f\|_{U^3(G)}^8$
cancel out leaving only the non-negative cutoffs $\psi$, whereas a Weyl sum computation
 reveals that $\E(f e(-\phi)) = O( N^{-c} )$ for any quadratic phase function $\phi$ and some explicit constant $c > 0$.  (Note that when $N$ is prime, the only quadratic phase functions $\phi$ are those of the form $\phi(x) = ax^2 + bx + c$ where $a, b, c \in \R/\Z$ with $Na=Nb=0$; see Lemma \ref{quadratic-classify}).
 We omit the details.
\end{example}

\ni The heart of the difficulty here is that the function $yM + z \mapsto yz/M$ is locally quadratic on the set $B := \{ yM+z: -M/10 \leq y,z \leq M/10 \}$, which is a fairly large subset of $G$, but does not extend (even approximately) to a \emph{globally} quadratic phase function on all of $G$.  These locally quadratic phase functions are thus a genuinely new class of obstructions to having small $U^d(G)$ norm which must now also be accounted for
in order to produce a genuine inverse theorem for the $U^3(G)$ norm.  Similar considerations also apply, of course, to the $U^d(G)$ norms for $d > 3$.\vs

\ni We must therefore understand the proper generalization of sets such as $B = \{ yM+z: -M/10 \leq y,z \leq M/10 \}$.   It turns out that there are two
ways to obtain such a generalization, which are in a sense dual to one another, namely that of \emph{generalized arithmetic progressions}
and that of \emph{Bohr sets}.  
For technical reasons it is convenient to work in the first instance with the latter notion, but we will discuss generalized arithmetic progressions later in the paper.

\begin{definition}[Bohr sets]\label{bohr-set}  
Let $G$ be a finite additive group, and let $S \subseteq  \widehat G$, $|S| = d$ be a subset of the dual group.  We define a sub-additive quantity $\| \|_S$ on $G$
by setting
$$ \|x\|_S := \sup_{\xi \in S} \| \xi \cdot x \|_{\R/\Z},$$
where $\|x\|_{\R/\Z}$ denotes the distance to the nearest integer, and
define the \emph{Bohr set} $B(S,\rho) \subseteq  G$ for any $\rho > 0$ to be the set
\[ B(S,\rho) := \{ x \in G : \|x\|_S < \rho \;\; \mbox{for all } \xi \in S \}.\]
\end{definition}

\ni Note that the dependence of the Bohr set $B(S,\rho)$ on $\rho$ can be rather discontinuous , as can be seen rather dramatically in finite field geometries such as $\mathbb{F}_3^n$. This is inconvenient in applications, but fortunately it was noted by Bourgain\footnote{In fact rather earlier Gowers \cite[Lemma 10.10]{gowers-long-aps} employed an argument which establishes that \textit{all} Bohr sets in $\Z/N\Z$, $N$ prime, are regular in a very weak sense.} \cite{Bou} that one may restrict attention to ``regular'' Bohr sets which enjoy some limited
continuity properties in $\rho$.

\begin{definition}[Regular Bohr sets \cite{Bou}]\label{regular-def} Let $S \subseteq \widehat{G}$, $|S| = d$, be a set of characters, and suppose that $\rho \in (0,1)$. A Bohr set $B(S,\rho)$ is said to be \emph{regular} if one has
$$ (1-100d |\kappa|) |B(S,\rho)| \leq |B(S, (1+\kappa)\rho)| \leq (1+100d |\kappa|) |B(S,\rho)|$$
whenever $|\kappa| \leq 1/100d$.
\end{definition}

\ni Lemma \ref{make-regular} gives a plentiful supply of regular Bohr sets. The constant 100 can be lowered but this will not concern us here. With this definition in place, we can now give the generalization of Theorem \ref{mainthm1} to arbitrary groups.  

\begin{theorem}[Inverse theorem for $U^3(G)$]\label{mainthm2}  Let $G$ be an finite additive group of odd order,
let $f: G \to \D$ be a bounded function and let $0 < \eta \leq 1$.  
\begin{enumerate}

\item If $\|f\|_{U^3(G)} \geq \eta$, then there exists a regular Bohr set $B := B(S,\rho)$ in $G$ with $|S| \leq (2/\eta)^C$ and $\rho \geq (\eta/2)^C$ such that
\begin{equation}\label{efphi-general}
\E_{y \in G} \|  f\|_{u^3(y+B)} \geq (\eta/2)^C,
\end{equation}
where it is permissible to take $C = 2^{24}$.
In particular, there exists $y \in G$ such that $\|f\|_{u^3(y+B)} \geq (\eta/2)^C$.

\item Conversely, $B = B(S,\rho)$ is a regular Bohr set, if $f: G \to \D$ is a bounded function and if $\| f \|_{u^3(y + B)} \geq \eta$, then we have
$$  \|f\|_{U^3(G)} \geq (\eta^3 \rho^2/C'd^3)^d$$
for some absolute constant $C'$.

\end{enumerate}
\end{theorem}

\ni Note that the $u^3(G)$ norm is no longer involved in this inverse theorem; this is necessary as demonstrated by Example \ref{fw-ex}, and has to
do with the lack of extendibility of some local quadratic phases to global ones. In later sections we will prove other, related, inverse theorems for the $U^3(G)$ norm. In \S \ref{sec10} we will obtain a result in which the quadratic phases are given quite explicitly when $G = \Z/N\Z$. Then, in Theorem \ref{nilsequence-thm}, we will provide a link to recent ergodic-theoretic work of Host-Kra and Ziegler .\vs

\ni In a future series of papers we will prove an enhanced version of Theorem \ref{mainthm2}, and use it to establish an asymptotic for the number of quadruples $p_1 < p_2 < p_3 < p_4 \leq N$ of primes in arithmetic progression. The enhancement required is that we must be able to deal with functions $f : \Z/N\Z \rightarrow \C$ which are not necessarily bounded, in particular functions such as $f = \Lambda - 1$, where $\Lambda$ is the von Mangoldt function.
Once this is done, one may analyse  
the $u^3$ norm of $f$ using what are, in essence, rather classical methods of analytic number theory such as Vaughan's decomposition of $\Lambda$.\vs

\ni A word of reassurance is perhaps in order for the reader interested in this result concerning primes. Although the present paper is long, only a few sections of it, namely sections \S \ref{Avg-sec}, \ref{sec5}, \ref{sec8} and \ref{ggc}, are relevant to that work. In fact it is hoped that the subsequent papers on primes will be readable largely independently of the present work.\vs

\ni Let us briefly mention the connection between the results here and those in \cite{tao:ergodic}.  In that paper the second author introduced the concept
of a \emph{uniformly almost periodic function of order $k-2$}, which generalized the concept of a polynomial phase of order at most $k-2$
(and which incorporates the ``locally polynomial phases'' discussed above.  One also obtained (by very elementary means) an inverse theorem 
for the $U^{k-1}$ norms involving these uniformly almost periodic functions, see \cite[Lemma 5.11]{tao:ergodic}.  However, because the uniformly almost periodic functions are a larger class than the locally polynomial phases, those results are weaker than the inverse theorems presented here.  Nevertheless, with substantial additional effort (involving for instance the van der Waerden theorem) it is possible to use the inverse
theorem for uniformly almost periodic functions to obtain another proof of Szemer\'edi's theorem. See \cite{tao:ergodic} for more details.
We also remark that very similar objects (the \emph{anti-uniform functions}) were also utilized in \cite{green-tao-primes} in order to reduce the task of establishing arbitrarily long progressions in the primes to Szemer\'edi's theorem.\vs

\ni Finally, let us offer a word of explanation for our policy concerning constants. For many of the arguments of this paper we have supplied exact constants, eschewing excessive use of the $O$-notation. This perhaps allows one to better see how bounds from different lemmas combine with one another to influence later bounds. Some readers may, however,  prefer to replace such quantities as $2^{90} \eta^{-384}$ with $C \eta^{-C}$ when reading the paper. 

\section{A model problem: global quadratic phase functions}\label{model-sec}

\ni We now present a simple result, namely the classification of globally quadratic phase functions on an arbitrary additive group $G$ of odd order, which we will need later, and which will serve to illustrate our strategy for the more advanced results we give below. \vs 

\ni Let us call a homomorphism $M: G \mapsto \widehat G$ \emph{self-adjoint} if we have
$$ Mx \cdot y - My \cdot x = 0 \hbox{ for all } x, y \in G.$$

\begin{lemma}[Inverse theorem for globally quadratic phase functions]\label{quadratic-classify}  Let $G$ be a finite additive group of odd order,
and let 
$\phi: G \to \R/\Z$ be a globally quadratic phase function.  Then there exists $c \in \R/\Z$, $\xi \in \widehat G$, and a self-adjoint homomorphism $M: G \to \widehat G$ such that $\phi(x) = Mx \cdot x + \xi \cdot x + c$.  
Conversely, all such functions $x \mapsto Mx \cdot x + \xi \cdot x + c$ are globally quadratic phase functions.
\end{lemma}

\begin{proof}  The converse is easy, so we focus on the forward direction.
It is convenient to adopt the notation $\Phase(x_1,\ldots,x_n)$ to denote an arbitrary function 
of variables $x_1,\ldots,x_n$ which takes values in $\R/\Z$, where $\Phase$ can vary from line to line or even within the same line.
This is useful for handling expressions whose exact value is not important for the argument, but whose functional 
dependencies on other variables needs to be recorded.\vs

\ni We shall give an argument which may seem a bit cumbersome (and is certainly not the shortest proof of this lemma), but it will serve to motivate the proof of Theorems \ref{mainthm1} and \ref{mainthm2}.  Indeed, this lemma can be thought of in some sense as the $\eta = 1$ case of those theorems.
First observe that if $\phi$ is a quadratic phase function, then $(h \cdot \nabla_x) \phi$ is a linear phase function for each $h$.  
Thus for each $h \in G$ there exists $\xi_h \in \widehat G$ such that
\begin{equation}\label{phihx}
 (h \cdot \nabla_x) \phi(x) = \phi(x+h) - \phi(x) = \xi_h \cdot x + \Phase(h) \hbox{ for all } x, h \in G.
 \end{equation}
The next step is to obtain some linearity on the map $h \mapsto \xi_h$, by using difference operators to eliminate various terms.  
Let $k \in G$ be arbitrary.  If we apply the difference operator $(k \cdot \nabla_x)$ to \eqref{phihx}
we can eliminate the $\Phase(h)$ term to obtain
\begin{equation}\label{thk}
(k\cdot \nabla) \phi(x+h) - (k \cdot \nabla) \phi(x) = \xi_h \cdot k \hbox{ for all } x, h, k \in G.
\end{equation}
If we then apply the difference operator $(h_1 \cdot \nabla_h)$ some $h_1 \in G$ to eliminate the $(k \cdot \nabla) \phi(x)$ term, we obtain
\begin{equation}\label{th1}
(h_1 \cdot \nabla) (k \cdot \nabla) \phi(x+h) = (h_1 \cdot \nabla_h) \xi_h \cdot k \hbox{ for all } x, h, k, h_1 \in G.
\end{equation}
Making the substitution $y = x+h$ we obtain
\begin{equation}\label{th1-y}
(h_1 \cdot \nabla) (k \cdot \nabla) \phi(y) = (h_1 \cdot \nabla_h) \xi_h \cdot k \hbox{ for all } y, h, k, h_1 \in G.
\end{equation}
If we then apply the difference operator $(h_2 \cdot \nabla_h)$ for some $h_2 \in G$ to eliminate the remaining $\phi$ term, we conclude
\begin{equation}\label{xihh}
 0 = (h_2 \cdot \nabla_h) (h_1 \cdot \nabla_h) \xi_h \cdot k \hbox{ for all } y, h, k, h_1 \in G.
 \end{equation}
Since $k$ is arbitrary, we conclude that 
\begin{equation}\label{xihh-2}
(h_2 \cdot \nabla_h) (h_1 \cdot \nabla_h) \xi_h = 0 \hbox{ for all } h_1, h_2 \in G.
\end{equation}
Thus if we write
\begin{equation}\label{xhio}
\xi_h = 2Mh + \xi_0  
\end{equation}
then we see that $M: G \mapsto \widehat G$ is a group homomorphism.
Note that we can insert the $2$ in front of $M$ because $|G|$ is odd; this factor of 2 will be convenient later.
Inserting this back into \eqref{phihx} we obtain
\begin{equation}\label{phihx-mh}
 (h \cdot \nabla_x) \phi(x) = 2Mh \cdot x + \xi_0 \cdot x + \Phase(h) \hbox{ for all } x, h \in G.
 \end{equation}
This is almost what we want, but we must somehow ``integrate'' the partial derivative $h \cdot \nabla_x$.  To do this we must first establish that $M$ is self-adjoint.  Informally, this self-adjointness reflects the symmetry $(h' \cdot \nabla) (h \cdot \nabla) \phi = (h \cdot \nabla) (h' \cdot \nabla) \phi$ of the
second derivative.  To make this rigorous (and in a manner which will extend suitably to more general situations)
we shall use the following ``symmetry argument''.  In order to focus on the $Mh \cdot x$ term, we shall write \eqref{phihx-mh} as
\begin{equation}\label{phase-1}
 \Phase(x+h) + \Phase(x) + \Phase(h) - 2Mh \cdot x = 0 \hbox{ for all } x,h \in G.
 \end{equation}
Substituting $x$ by $y$ and subtracting to eliminate $\Phase(h)$, we obtain
\begin{equation}\label{phase-2}
\Phase(y+h) + \Phase(y) + \Phase(x+h) + \Phase(x) - 2Mh \cdot (y-x) = 0 \hbox{ for all } x,y,h \in G.
\end{equation}
Making the substitution $z = x+y+h$, we conclude
\begin{equation}\label{phase-3}
\Phase(x,z) + \Phase(y,z) - 2M(z-x-y) \cdot (y-x) = 0 \hbox{ for all } x,y,z \in G.
\end{equation}
Absorbing as many terms as possible into the unspecified functions $\Phase$, we conclude that
\begin{equation}\label{phase-4}
 \Phase(x,z) + \Phase(y,z) + 2 \{x,y\} = 0 \hbox{ for all } x,y,z \in G,
 \end{equation}
where $\{x,y\}$ is the anti-symmetric form
\begin{equation}\label{xy-def}
 \{x,y\} := Mx \cdot y - My \cdot x.
 \end{equation}
Freezing the value of $z$, we conclude that
\begin{equation}\label{phase-5}
\Phase(x) + \Phase(y) + 2 \{x,y\} = 0 \hbox{ for all } x,y \in G.
\end{equation}
Replacing $y$ by $y'$ and subtracting to eliminate the $\Phase(x)$ factor, we conclude
\begin{equation}\label{phase-6}
 \Phase(y) + \Phase(y') + 2\{x, y'-y\} = 0 \hbox{ for all } x,y,y' \in G.
 \end{equation}
Thus for each $y, y'$, the linear function $x \mapsto \{x,y'-y\}$ is independent of $x$ and is hence always zero\footnote{There appear to be some intriguing parallels with symplectic geometry here.  Roughly speaking, the vanishing \eqref{phase-7} is an assertion that the
 graph $\{(h,Mh): h \in G \}$ is a ``Lagrangian manifold'' on the ``phase space'' $G \times \widehat G$.  This graph can also be interpreted (essentially) as the ``wave front set'' $\{ (x, \nabla \phi(x)): x \in G \}$ of the original function $e(\phi)$.  A similar interpretation persists in the proofs of Theorem \ref{mainthm1} and \ref{mainthm2} below.  Thus we see hints of some kind of ``combinatorial symplectic geometry'' emerging, though we do not see how to develop these possible connections further.}:
\begin{equation}\label{phase-7}
\{ x, y' - y \} = 0 \hbox{ for all } x,y,y' \in G.
\end{equation}
Thus $M$ is self-adjoint.\vs

\ni We return now to \eqref{phihx-mh}, which we write as
\begin{equation}\label{bias-1}
\Phase(h) + \Phase(x+h) - \phi(x) - 2Mh \cdot x = 0 \hbox{ for all } x,h \in G.
\end{equation}
From the self-adjointness of $M$ we have
\begin{equation}\label{phase-decompose}
2Mh \cdot x = M(x+h) \cdot(x+h) - Mx \cdot x - Mh \cdot h
\end{equation}
and hence
\begin{equation}\label{bias-2}
 \Phase(h) + \Phase(x+h) - \phi(x) + Mx \cdot x = 0 \hbox{ for all } x,h \in G.
\end{equation}
In particular, if we apply difference operators in the $h$ and $x+h$ variables we see that
the phase function $-\phi(x) + Mx \cdot x$ is linear:
\begin{equation}\label{bias-3}
(h_1 \cdot \nabla_x) (h_2 \cdot \nabla_x) (-\phi(x) + Mx \cdot x) = 0 \hbox{ for all } x, h_1, h_2 \in G.
\end{equation}
Hence there exists $\xi \in \widehat G$ 
and a $c \in \R/\Z$ such that
\begin{equation}\label{bias-4}
\phi(x) -  Mx \cdot x = \xi \cdot x  + c.
\end{equation}
The claim follows.
\end{proof}\vs

\ni As one corollary of this classification, we obtain the following ``quadratic extension theorem''.

\begin{proposition}[Quadratic Extension Theorem]\label{QHB}
Let $G$ be an additive group, let $H \leq G$ be a subgroup, and suppose that $y \in G$. Then any quadratic phase function $\phi : y + H \rightarrow \R/\Z$ can be extended \textup{(}non-uniquely in general\textup{)} to a globally quadratic phase function on $G$.
\end{proposition}
\ni\textit{Remark.} An important theme of this paper is that this behaviour is specific to cosets $y+H$, and breaks down for other sets such as Bohr sets.\vs

\ni As a consequence of Proposition \ref{QHB}, let us now establish the (easy) second part of Theorem \ref{mainthm1}.
The first inequality $\|f\|_{U^3(G)} \geq \|f\|_{u^3(G)}$ follows from \eqref{ud-inverse-easy}, so we focus
on the second inequality $\|f\|_{u^3(G)} \geq \frac{|W|}{|G|} \| f \|_{u^3(y+W)}$.  It suffices to show that
$$\|f\|_{u^3(G)} \geq \frac{|W|}{|G|} \big|\E_{x \in y+W} f(x) e(-\phi(x))\big|$$
whenever $\phi$ is a locally quadratic function on $y+W$.  But by the preceding discussion, we can extend $\phi$ to all of $G$, and write
$$ \frac{|W|}{|G|} \big|\E_{x \in y+W} f(x) e(-\phi(x))\big| = \big|\E_{x \in G} f(x) 1_{W}(x-y) e(-\phi(x))\big|.$$
But by Fourier inversion we may write
$$ 1_{W}(x-y) = \E_{\xi \in W^\perp} e(\xi \cdot (x-y))$$
where $W^\perp := \{ \xi \in \widehat G: \xi \cdot x = 0 \hbox{ for all } x \in W \}$.  Thus we have
$$ \E_{x \in G} f(x) 1_{W}(x-y) e(-\phi(x)) = \E_{\xi \in W^\perp} \E_{x \in G} f(x) e(\xi \cdot(x-y) - \phi(x)).$$
But since $\xi \cdot(x-y) - \phi(x)$ is globally quadratic in $x$, we have
$$ |\E_{x \in G} f(x) e(\xi \cdot(x-y) - \phi(x))| \leq \|f\|_{u^3(G)}$$
and the claim follows from the triangle inequality.  Note that this argument in fact works for arbitrary groups $G$ (with $W$ now being 
a subgroup of $G$ rather than a subspace).

\section{Averaging lemmas}\label{Avg-sec}

\ni In this section we collect some very simple averaging estimates which we shall rely frequently on in the sequel.  

\begin{lemma}[Averaging on a subgroup]\label{avg-1}  Let $G$ be an additive group, let $H$ be a finite subgroup of $G$, and let $A \subseteq H$ be non-empty.  Let $f: H \to \C$ be a function.  Then
$$ \E_{x \in H} \E_{y \in x + A}f(y)  = \E_{y \in H}f(y).$$
In particular, by the pigeonhole principle there exists $x \in H$ such that
$$ |\E_{y \in x + A}f(y)| \geq |\E_{y \in H}f(y)|.$$
\end{lemma}

\begin{proof} Since $H+h = H$ for all $h \in A$ we have
$$ \E_{x \in H} f(x+h) = \E_{y \in H}f(y) \hbox{ for all } h \in A.$$
Averaging this over all $h \in A$ we obtain the first claim, and the second claim then follows from the pigeonhole principle.
\end{proof}

\ni This Lemma will be adequate for our purposes when we are in the finite field geometry case, because we will have plenty of subgroups available.  In the general group case, however, we will also need a more general type of averaging principle. The next lemma contains several rather similar formulations of such a result; all of them will be useful later on.

\begin{lemma}[Averaging on a Bohr set]\label{avg} Let $S \subseteq  \widehat{G}$ be a set of $d$ characters and let $0 < \rho < 1$ and $\eps \leq 1/200d$ be parameters. Suppose that the Bohr set $B := B(S,\rho)$ is regular, and let $A \subseteq  B(S,\eps\rho)$ be any set. Finally, let $f : G \rightarrow \D$ be any function. Then
\begin{enumerate}
\item $\big| \E_{x \in B} f(x) - \E_{x \in B} f(x) \big| \leq 200 d\eps$ if $y \in A$;
\item $\big| \E_{x \in B} f(x) - \E_{x \in B} \E_{y \in x + A} f(x) \big| \leq 200d \eps$;
\item There is some $x \in B$ for which $\E_{y \in x + A}f(y) \geq \E_{y \in B}f(y) - 200d\eps$;
\item There is some $x \in B(S,(1 - \eps)\rho)$ such that $\E_{y \in x + A}f(y) \geq \E_{y \in B}f(y) - 500 d\eps$.
\end{enumerate} 
\end{lemma}
\begin{proof} To prove (i) we must check that 
\[ \big| \sum_{x \in B + y} f(x) - \sum_{x \in B} f(x) \big| \leq 200 d\eps |B|.\]
This follows from the fact that $B(S,\rho)$ and $y + B(S,\rho)$ differ in at most $200d\eps$ elements. Indeed it is easy to see that
\[ B(S,\rho) \triangle (y + B(S,\rho)) \subseteq  B(S,(1 + \eps)\rho) \setminus B(S,(1 - \eps) \rho)\] ($\triangle$ denotes symmetric difference),
and this latter set has size no more than $200 d\eps$ by regularity.\\
(ii) follows from (i) and the triangle inequality. Indeed $|\E_{x \in B} f(x + y) - \E_{x \in B} f(x)| \leq 200d\eps$, whence
\[ \big| \E_{x \in B} \E_{y \in A} f(x + y) - \E_{x \in B} f(x) \big| \leq 200d\eps,\] which is (ii).\\
(iii) is immediate from (ii) and the pigeonhole principle.\\
(iv) follows from the pigeonhole principle and the fact that 
\[ \E_{x \in B} f(x) - \E_{x \in B(S,(1 - \eps)\rho)} \E_{y \in x + A} f(y) \big| \leq 500 d\eps,\]
which is in turn implied by (ii) and the bound
\[ \big| \E_{x \in B} F(x) - \E_{x \in B(S,(1 - \eps)\rho)} F(x) \big| \leq 300\epsilon d,\]
 valid for any $F : G \rightarrow \D$. To confirm this, note that 
 \[ \big| \sum_{x \in B} F(x) - \sum_{x \in B(S,(1 - \eps)\rho)} F(x)\big| \leq |B \setminus B(S, (1 - \eps)\rho)| \leq 100\eps d |B|.\]
 Therefore
 \begin{eqnarray*}
 \big| \E_{x \in B} F(x) - \E_{x \in B(S,(1 - \eps)\rho)} F(x)\big| & \leq & \frac{1}{|B|} \bigg| \sum_{x \in B} F(x) - \sum_{x \in B(S,(1-\eps)\rho)} F(x) \bigg| \\ & & \quad + \bigg| \frac{1}{|B|} - \frac{1}{|B(S,(1 - \eps)\rho)|}\bigg|\cdot  \bigg|\sum_{x \in B(S,(1 - \eps)\rho)} F(x)\bigg| \\ & \leq & 100\eps d + \bigg| \frac{|B|}{|B(S,(1 - \eps)\rho)|} - 1 \bigg| \\ & \leq & 100 \eps d + \frac{1}{1 - 100\eps d} - 1 \leq 300\eps d.
 \end{eqnarray*}
This completes the proof of Lemma \ref{avg}.\end{proof}\vs

\ni We adopt the following useful notation, analogous to the $\Phase$ notation used in proving Lemma \ref{quadratic-classify}.  
When $g(x_1,\ldots,x_n)$ is a complex-valued function of certain
variables $x_1,\ldots,x_n$ with $\Vert g \Vert_{\infty} \leq 1$, we shall refer to $g(x_1,\ldots,x_n)$
instead as $\Bounded(x_1,\ldots,x_n)$.  The notation $\Bounded$ thus denotes a function with $\Vert \Bounded \Vert_{\infty} \leq 1$, but the notation may refer to different functions from line to line, or even on the same line (similar to the $O$ notation,
or the use of the unspecified constants $C$).\vs

\ni We now record a basic application of the Cauchy-Schwarz inequality, whose proof is immediate.

\begin{lemma}[Cauchy-Schwarz]\label{cz}  
Let $X, Y$ be finite sets, and let $f: X \times Y \to \C$ be a function.  Then for any bounded function $\Bounded(x)$ of $X$, we have
\[
|\E_{x, y} f(x,y) \Bounded(x)|  \leq \E_x |\E_y f(x,y)| 
\leq (\E_x |\E_y f(x,y)|^2  )^{1/2} = (\E_{x,y,y'} f(x,y') \overline{f(x,y)})^{1/2}.
\]
\end{lemma}

\ni In the special case when $Y=G$ is a group, we conclude in particular the \emph{Van der Corput inequality}
\begin{equation}\label{vdc}
|\E_{x \in X,y \in G} f(x,y) \Bounded(x)| \leq |\E_{x \in X, y,h \in G} T^h_y f(x,y) \overline{f(x,y)}|^{1/2}
\end{equation}
whenever $f: X \times G \to \C$ is a function, which follows by using the substitution $y' = y+h$.  This inequality is very useful for eliminating
unknown bounded functions $\Bounded(x)$ in an expression to be estimated. Using the van der Corput inequality, we can now prove Proposition \ref{gvn}.\vs

\ni\textit{Proof of Proposition \ref{gvn}.} It suffices to prove the more general statement\footnote{Indeed, the $U^d(G)$ norms are capable of controlling even more general expressions, for instance when the linear shifts $a_j h$ are replaced by polynomial shifts with no constant coefficient; this is implicitly in \cite{bergelson-leibman1}.}
\begin{equation}\label{taj}
 |\E_{x,h} \prod_{j \in J} T^{a_j h} f_j(x)| \leq \| f_{j_0} \|_{U^{|J|-1}(G)}
\end{equation}
for all finite sets $J$ with $|J| \geq 1$, 
all $j_0 \in J$, all bounded functions $(f_j)_{j \in J}$, and all distinct integers $(a_j)_{j \in J}$ such that $a_j - a_{j'}$
is coprime to $|G|$ for all distinct $j, j' \in J$.\vs

\ni We induct on $|J|$.  When $|J|=1$ the claim is trivial from \eqref{ud-recursive}, 
so suppose $|J| \geq 2$ and the claim has already been proven for smaller
values of $J$.  Let $j_1$ be an element in $J \backslash \{j_0\}$.  By making the change of variables $x \to x + a_{j_1} h$ if necessary we may assume
that $a_{j_1} = 0$.  Since $f_{j_1}$ is bounded, we can then express the left-hand side of \eqref{taj} as
$$ \big|\E_{x,h \in G}\big( \Bounded(x) \prod_{j \in J \backslash \{j_1\}} T^{a_j h} f_j(x)\big)\big|,$$
which by \eqref{vdc} can be bounded by
$$ \E_{k \in G}\big( \E_{x,h \in G}\big( \prod_{j \in J \backslash \{j_1\}} T^{a_j h} (T^{a_j k} f_j \overline{f_j})(x)\big)\big)^{1/2}.$$
Applying the inductive hypothesis \eqref{taj} to the inner expectation, we can bound this in turn by
$$ \E_{k \in G}( \| T^{a_j k} f_j \overline{f_j} \|_{U^{|J|-2}(G)})^{1/2},$$
which by H\"older's inequality and the substitution $h := a_j k$ (noting that $(N,a_j - a_{j_1}) = (N,a_j) = 1$) is bounded by
$$ \E_{h \in G}\big( \| T^h f_j \overline{f_j} \|_{U^{|J|-2}(G)}^{2^{|J|-2}} |\big)^{1/2^{|J|-1}}.$$
The claim then follows from \eqref{ud-recursive}.
\endproof\vs

\ni Let us now consider averages of the form
$$ \E_{z \in B',x \in B}( f(z) \Bounded(x) \Bounded(z+x))$$
where $B, B'$ are non-empty subsets of an additive group $G$, and $f,g,h$ are bounded functions.  If $B=B'=G$, then a variant of Proposition \ref{gvn} shows that this quantity is bounded by $\|f\|_{U^2(G)}$, and hence (by Proposition \ref{u2-inverse}) if the above average is large, then $f$ must
correlate with a linear phase function.  It turns out that a similar statement is true for arbitrary $B, B'$, provided that $B+B'$ is only a little bit larger than $B$.

\begin{lemma}[Large trilinear form implies correlation with linear phase]\label{fgh-bohr}  Let $B$, $B'$ be two non-empty subsets of an additive group $G$.  Then we have
$$ \| f \|_{u^2(B')} \geq \frac{\E(1_{B})}{\E(1_{B+B'})} |\E_{z \in B',x \in B}( f(z) \Bounded_1(x) \Bounded_2(z+x))|$$
for any $f: G \to \C$ and any two bounded functions $\Bounded_1$, $\Bounded_2$.
\end{lemma}

\begin{proof}  Without loss of generality we may assume that $f,\Bounded_1,\Bounded_2$ vanish outside of $B'$, $B$, $B+B'$ respectively.
From Fourier expansion we have
\begin{align*}
\E_{z \in B', x \in B}( f(z) \Bounded_1(x) \Bounded_2(z+x) ) &= \frac{1}{\E(1_{B'}) \E(1_B)} \E_{y \in G}( f*\Bounded_1(y) \Bounded_2(y) ) \\
&= \frac{1}{\E(1_{B'}) \E(1_B)} \sum_{\xi \in \widehat G} \widehat f(\xi) \widehat \Bounded_1(\xi) \widehat \Bounded_2(-\xi).
\end{align*}
On the other hand, from Plancherel we have $\sum_{\xi \in \widehat G} |\widehat \Bounded_1(\xi)|^2 \leq \E(1_B)$
and $\sum_{\xi \in \widehat G} |\widehat \Bounded_2(-\xi)|^2 \leq \E(1_{B+B'})$.  From H\"older's inequality we thus conclude that
$$ |\E_{z \in B',x \in B}( f(z) \Bounded_1(x) \Bounded_2(z+x))| \leq \frac{\E(1_{B+B'})}{\E(1_{B'}) \E(1_B)} \sup_{\xi \in \widehat G} |\widehat f(\xi)|,$$
and hence there exists $\xi \in \widehat G$ such that
$$ |\E_{z \in B'}( f(z) e(-\xi \cdot z)  )| \geq |\E_{z \in B',x \in B}( f(z) \Bounded_1(x) \Bounded_2(z+x) )| \frac{\E(1_{B})}{\E(1_{B+B'})}.$$
The claim follows.
\end{proof}

\section{An argument of Gowers}\label{sec5}

\ni We now begin the proofs of the inverse theorems in Theorem \ref{mainthm1} and Theorem \ref{mainthm2}.  As we shall see, the arguments shall
be analogous to those used to prove Lemma \ref{quadratic-classify}, and closely follow the treatment in Gowers \cite{gowers-4-aps}.  
The first part of the argument is to establish a ``phase derivative'' $h \to \xi_h$ for the function $f$ and to establish some additivity properties on this phase derivative.  These arguments apply to arbitrary finite additive groups $G$; later on we shall treat the finite
field case separately.

\begin{proposition}[Large $U^3$ norm gives many additive quadruples \cite{gowers-4-aps}]\label{add-quad}  Let $G$ be an arbitrary finite additive group, and let $f: G \to \D$ be a bounded function such that $\|f\|_{U^3(G)} \geq \eta$ for some $\eta > 0$.  Then there exists a set $H \subseteq  G$, and 
a function $\xi: H \to \widehat G$ whose graph $\Gamma := \{ (h,\xi_h): h \in H \} \subseteq  G \times \widehat G$ obeys the estimate
\begin{equation}\label{Zzz}
 |\{ (z_1,z_2,z_3,z_4) \in \Gamma: z_1 + z_2 = z_3 + z_4 \}| \geq 2^{-8}\eta^{64}N^3.
 \end{equation}
Furthermore for each $(h,\xi_h) \in \Gamma$ we have
$$ |\E_x T^h f(x)  \overline{f}(x) e(- \xi_h \cdot x)| \geq \eta^4 / 2^{1/2}.$$
\end{proposition}

\begin{proof} 
As we shall see, this proposition corresponds fairly closely with the first part of the proof of Lemma \ref{quadratic-classify} (up to \eqref{xihh-2}). From \eqref{ud-recursive} we have
$$ \E_{h \in G} \| T^h f \overline{f} \|_{U^2(G)}^4 \geq \eta^8.$$
Applying Proposition \ref{u2-inverse} we conclude that
$$ \E_{h \in G} \| T^h f \overline{f} \|_{u^2(G)}^2 \geq \eta^8.$$
Thus if we let 
$$ H := \{ h \in G: \| T^h f \overline{f} \|_{u^2(G)}^2 \geq \eta^8/2 \}$$
then we have
$$ \E_{h \in G} \| T^h f \overline{f} \|_{u^2(G)}^2 1_{G \backslash H}(h) \geq \eta^8/2$$
and hence
$$ \E_{h \in G} \| T^h f \overline{f} \|_{u^2(G)}^2 1_H(h)  \geq \eta^8/2.$$
In particular we have 
\begin{equation}\label{eih}
\E(1_H) \geq \eta^8/2.
\end{equation}
By \eqref{u2-infty} and definition of $H$, we can find a map $h \mapsto \xi_h$ from $H$ to $\widehat G$
such that
\begin{equation}\label{eq5.2a} |\E_x T^h f(x) \overline{f}(x) e(-\xi_h \cdot x)| \geq \eta^4 / 2^{1/2} \hbox{ for all } h \in H.
\end{equation}  Let us fix this map $h \mapsto \xi_h$.
We square sum the above expression in $h$ and use \eqref{eih} to conclude
$$ \E_h |\E_x T^h f(x)  \overline{f}(x) e(-\xi_h \cdot x)|^2 1_H(h)  \geq \eta^{16} / 4.$$
But from the identity
$$
|\E_x T^h f(x) \overline{f}(x) e(-\xi_h \cdot x)|^2 
= \E_{x,k} T^{k}(T^h f)(x) \overline{T^h f}(x) \overline{T^k f(x)} f(x) e(\xi_h \cdot k)
$$
we conclude
\begin{equation}\label{eq5.2b} |\E_{x,h,k} T^{k}(T^h f)(x) \overline{T^h f}(x) e(\xi_h \cdot k) 1_H(h) \overline{T^k f(x)} f(x)| \geq \eta^{16}/4.\end{equation}
At this point we suppress the explicit mention of the functions $f$ and write this simply as
$$ |\E_{x,h,k \in G}( \Bounded(x+h,k) 1_H(h) e(\xi_h \cdot k) \Bounded(x,k))| \geq \eta^{16}/4.$$
Applying \eqref{vdc} to eliminate the $\Bounded(x,k)$ factor, we conclude
\begin{equation}\label{eq5.2c} \E_{h,h_1,x,k} \Bounded(x+h+h_1,k) \Bounded(x+h,k) e((h_1 \cdot \nabla_h) \xi_h \cdot k) 
1_H(h+h_1) 1_H(h) \geq \eta^{32}/16.\end{equation}
Making the substitution $y = x+h$, we obtain
\begin{equation}\label{eq5.2d} \E_{h,y,h_1,k} e((h_1 \cdot \nabla_h) \xi_h \cdot k) 1_H(h+h_1) 1_H(h) \Bounded(y,h_1,k)
 \geq \eta^{32}/16.\end{equation}
Applying \eqref{vdc} again we conclude
\begin{equation}\label{eq5.2e} \E_{h,h_1,h_2,y,k} e((h_2 \cdot \nabla_h) (h_1 \cdot \nabla_h) \xi_h \cdot k) 1_H(h+h_1+h_2) 1_H(h+h_2) 1_H(h+h_1) 1_H(h) \geq 2^{-28}\eta^{64}.\end{equation}
Summing this in $k$ using the Fourier inversion formula, and discarding
the irrelevant $y$ averaging, we infer
\begin{equation}\label{eq5.2f}
 \E_{h,h_1,h_2} 1_{(h_2 \cdot \nabla_h) (h_1 \cdot \nabla_h) \xi_h = 0} 1_H(h+h_1+h_2) 1_H(h+h_2) 1_H(h+h_1) 1_H(h) 
\geq \eta^{64}/256.
\end{equation}  The claim now follows by substituting $z_1 := (h,\xi_h)$, $z_2 := (h+h_1,\xi_{h+h_1})$, $z_3 := (h+h_2,\xi_{h+h_2})$, and 
$z_4 :=  (h+h_1+h_2, \xi_{h+h_1+h_2})$.
\end{proof}\vs

\ni As remarked, the analogy between this argument and the first part of the proof of Lemma \ref{quadratic-classify} is very close. In particular equations \eqref{eq5.2a}, \eqref{eq5.2b}, \eqref{eq5.2c}, \eqref{eq5.2d}, \eqref{eq5.2e} and \eqref{eq5.2f} are analogues of \eqref{phihx}, \eqref{thk}, \eqref{th1}, \eqref{th1-y}, \eqref{xihh} and \eqref{xihh-2} respectively.\vs

\ni In order to exploit the conclusion \eqref{Zzz}, we require two very useful results from additive combinatorics. The first result asserts that a set $\Gamma'$ with some partial additive structure (in the sense that it contains many additive quadruples) can be refined to have a more complete additive structure (in the sense that its sum set is small).  This type of result was first obtained by Balog and Szemer\'edi \cite{balog}, but with very weak constants.
A version of the theorem with polynomial dependencies between the constants was obtained by Gowers \cite{gowers-4-aps}. The version we quote below, with rather good powers in those polynomials, may be proved by a careful working of the argument in Chang \cite{chang-er}. Of course for the purposes of this paper the precise values of the constants are somewhat unimportant.\vs

\ni Write $\Gamma' - \Gamma' := \{ z - z': z, z' \in \Gamma'\}$ for the difference set of $\Gamma'$.

\begin{theorem}[Balog-Szemer\'edi-Gowers theorem \cite{chang-er}]\label{bsg}  Let $G$ be an additive group, and let $\Gamma$ be a finite non-empty subset of $G$ such that
$$ |\{ (z_1,z_2,z_3,z_4) \in \Gamma: z_1 + z_2 = z_3 + z_4 \}| \geq |\Gamma|^3 / K$$
for some $K \geq 1$.  Then there exists a subset $\Gamma' \subseteq \Gamma$ such that
$$ |\Gamma'| \geq 2^{-6} K^{-1}|\Gamma| \hbox{ and } |\Gamma' - \Gamma'| \leq 2^{42} K^6 |\Gamma'|.$$
\end{theorem}
  
\ni The other tool we need is the Pl\"unnecke inequality, a simple proof\footnote{Actually, since we only need this theorem for bounded $k, l$, and would not be concerned if the bound of $K^{k+l}$ were worsened to $K^{C(k+l)}$ for some absolute constant $C$, it is possible to modify the proof of Theorem \ref{bsg} in order to gain control on $|k\Gamma' - l\Gamma'|$ directly, without needing the Pl\"unnecke inequalities.  However this would not
simplify the remainder of the argument and so we do not give the details of this alternate approach here.} of which can be found in \cite{ruzsa-graph}.

\begin{theorem}[Pl\"unnecke inequalities \cite{plun,ruzsa-graph}]\label{plun}  Let $G$ be an arbitrary additive group, and let $\Gamma'$ be a finite non-empty subset of $G$ such that $|\Gamma' - \Gamma'| \leq K |\Gamma'|$ for some $K \geq 1$.  Then we have
$|k \Gamma' - l \Gamma'| \leq K^{k+l} |\Gamma'|$ for all $k, l \geq 1$, where $k \Gamma'$ denotes the $k$-fold sumset of $\Gamma'$.
\end{theorem}

\ni In our applications the integers $k, l$ will be quite small, in fact they will not exceed 9.  Combining Proposition \ref{add-quad} with
Theorem \ref{bsg} and Theorem \ref{plun}, we conclude

\begin{proposition}[$h \mapsto \xi_h$ is nearly affine-linear]\label{add-quad-2}  Let $G$ be an arbitrary finite additive group, and let $f: G \to \D$ be a bounded function such that $\|f\|_{U^3(G)} \geq \eta$ for some $\eta > 0$.  Then there exists a set $H' \subseteq  G$ 
a function $\xi: H' \to \widehat G$ whose graph $\Gamma' := \{ (h,\xi_h): h \in H' \} \subseteq G \times \widehat G$ obeys the estimates
\begin{equation}\label{gamma-small}
|\Gamma'| \geq 2^{-14}\eta^{64} N 
\end{equation}
and
$$
|k\Gamma' - l\Gamma'| \leq (2^{90}\eta^{-384})^{k+l} |\Gamma'|\hbox{ for all } k,l \geq 1.
$$
Furthermore for each $(h,\xi_h) \in \Gamma'$ we have
\begin{equation}\label{thfxi}
 |\E_{x \in G}( T^h f(x)  \overline{f(x)} e(-\xi_h \cdot x))| \geq \eta^4/2.
 \end{equation}
\end{proposition}

\ni Thus far we have not used anything about the underlying group $G$ other than it is finite and abelian.  We could continue doing this,
proving Theorem \ref{mainthm2} and extracting Theorem \ref{mainthm1} as a corollary.  However for expository reasons we will now restrict
to the finite field geometry case $G = \Ffiven$ and give a complete proof of Theorem \ref{mainthm1}.  The arguments there serve as a simplified model which
will help motivate the general case.

\section{The finite field case}\label{finite-sec}

\ni We now restrict attention to a finite field geometry setting $G = \Ffiven$, and prove Theorem \ref{mainthm1} (i). In this section, then, $N = 5^n$. Thanks to Proposition \ref{add-quad-2},
we have already isolated a phase derivative $h \mapsto \xi_h$ which exhibits some linear behavior.  The first step (following Gowers \cite{gowers-4-aps,gowers-long-aps}) is to show that $\xi$ in fact
matches up with a linear phase function on a large subspace; then we shall show that in fact $\xi$ matches up with a \emph{self-adjoint} linear phase function on a large subspace (this is the substantially new part of the argument).  
Finally, we shall again follow Gowers \cite{gowers-4-aps} and
conjugate $f$ by a quadratic phase to eliminate this linear phase derivative and conclude the argument.\vs

\ni\textit{Step 1: Linearization of phase derivative.} In this subsection we establish

\begin{proposition}[Graphs have large linear component \cite{gowers-4-aps,gowers-long-aps}]\label{quad-linear} 
Let $H' \subseteq  \Ffiven$, and let $\xi: H' \to (\Ffiven)^*$ be a function whose graph $\Gamma' := \{ (h, \xi_h): h \in H' \} \subseteq \Ffiven \times (\Ffiven)^*$ obeys the estimates
\begin{equation}\label{gamma-large}
K^{-1} N \leq |\Gamma'| \leq |9\Gamma' - 8\Gamma'| \leq K N
\end{equation}
for some $K \geq 1$.  Then there exists a linear subspace $V \leq \Ffiven$ with the codimension bound
$$ n - \dim(V) \leq (2K)^6$$
and a translate $x_0 + V$ of this subspace, together with a linear transformation $M: V \to (\Ffiven)^*$ and an element $\xi_0 \in (\Ffiven)^*$ such that
\begin{equation}\label{hxixi}
 \E_{h \in V}( 1_H(x_0 + h) 1_{\xi_{x_0+h} = 2Mh + \xi_0} ) \geq (2K)^{-4}.
\end{equation}
\end{proposition}

\begin{proof}
Note that while $\Gamma'$ is a graph, the slightly larger sets $k\Gamma' - l\Gamma'$ need not be a graph.  The next lemma shows that this can be rectified by passing to an appropriate subset $\Gamma'' \subseteq \Gamma'$.

\begin{lemma}\label{random-slice}  There exists a subset $\Gamma''$ of $\Gamma'$ with \begin{equation}\label{eq6.03}|\Gamma''| \geq N/5K^3\end{equation}
such that $4\Gamma''-4\Gamma''$ is a graph.
\end{lemma}

\ni \textit{Remarks.} When $G = \mathbb{Z}/N\mathbb{Z}$ is a cyclic group of prime order, this lemma is essentially \cite[Lemma 7.5]{gowers-long-aps}, and our arguments both here and in \S \ref{ggc} are in a similar spirit. Curiously, the argument of \cite{gowers-long-aps} uses the fact that $\mathbb{Z}/N\mathbb{Z}$ is a field, and so this is a rare instance of the cyclic group case being somewhat easier than the general group case. The conclusion can also be rephrased as an assertion that the map $h \mapsto \xi_h$ is
a \emph{Freiman $8$-homomorphism} on $H''$. \vs

\begin{proof}  Let $A \subseteq (\Ffiven)^*$ be the set of all $\xi$ such that $(0,\xi) \in 8\Gamma' - 8\Gamma'$.  Observe that since $\Gamma'$ is a graph,
we have $|\Gamma' + A| = |\Gamma'| |A|$.  On the other hand, $\Gamma' + A$ is contained in $9\Gamma' - 8\Gamma'$.  Applying
\eqref{gamma-large}, we conclude that $|A| \leq K^2$.  Now let $m = \lceil \log_5 |A| \rceil$, and let $\Psi: (\Ffiven)^* \to \mathbb{F}_5^m$ be a randomly chosen linear transformation from $(\Ffiven)^*$ to $\mathbb{F}^m$.  Observe that for each non-zero $\xi \in A$, we have $\Psi(\xi) \neq 0$ with probability at least $1 - 5^{-m}$.  Thus we see that with non-zero probability $\Psi$ is non-zero on all of $A \backslash \{0\}$.\vs

\ni Fix $\Psi$ with the above properties, and let $c$ be a randomly selected point in $\mathbb{F}^m$, and define $\Gamma'' := \{ (h,\xi_h) \in \Gamma': \Psi(\xi_h) = c \}$.  Then the expected size of $|\Gamma''|$ is at least $|\Gamma'| / 5^m \geq N/5K^2$.  Also, observe that if $(0,\xi)$ lies in
$8\Gamma'' - 8\Gamma''$ then $\xi$ lies in $A$ and (by linearity of $\Psi$) $\Psi(\xi) = 0$, hence $\xi = 0$.  Since $8\Gamma'' - 8\Gamma''$ is the difference set of $4\Gamma''-4\Gamma''$, this implies that $4\Gamma''-4\Gamma''$ is a graph.  The claim follows.
\end{proof}\vs

\ni Define $H''$ by $\Gamma'' := \{ (h,\xi_h): h \in H'' \}$.
Now we use a result of Bogolyubov \cite{bogolyubov} (which we shall also utilize later in this paper) to obtain some control on the set $2H''-2H''$.  

\begin{lemma}[Bogolyubov lemma \cite{bogolyubov}]\label{bog-lemma} Let $A$ be a subset of a finite additive group $G$ such that
$|A| \geq \delta |G|$.  Then there exists a set $S \subseteq  \widehat G$ with $|S| \leq 2\delta^{-2}$ such that
$B(S, \frac{1}{4}) \subseteq  2A-2A$.
\end{lemma}
\ni\textit{Remark.} Somewhat sharper versions of this lemma are known, see \cite{chang}, but we will not need these here.\vs

\begin{proof}
From Fourier inversion we have
$$ 1_A(x) = \sum_{\xi \in \widehat G} \widehat 1_A(\xi) e(\xi \cdot x).$$
Convolving this with itself four times, we obtain
$$ 1_A * 1_A * 1_{-A} * 1_{-A}(x) = \sum_{\xi \in G} |\widehat 1_A(\xi)|^4 e(\xi \cdot x)$$
where $f*g(x) := \E_y f(y) g(x-y)$ is the normalized convolution operation.  Since the left-hand side is only
non-zero on $2A-2A$, we conclude that
$$ \{ x \in G: \Re \sum_{\xi \in G} |\widehat 1_A(\xi)|^4 e(\xi \cdot x) > 0 \} \subseteq  2A-2A.$$
Let us now study the sum on the left-hand side.  Let $0 < \alpha < \delta$ be a parameter to be chosen later, and let
$S := \{ \xi \in G: |\widehat 1_A(\xi)| \geq \alpha \}$ denote the large Fourier coefficients of $A$.  Since $\alpha < \E(1_A)$, we have
$0 \in S$.  Also, from the Plancherel identity
$$ \sum_{\xi \in G} |\widehat 1_A(\xi)|^2 = \E( 1_A ) \geq \delta$$
and Chebyshev's inequality we have an upper bound on the cardinality of $S$:
$$ |S| \leq \delta / \alpha^2.$$
Now suppose that $x \in B(S, \frac{1}{4})$.  Then we have $\cos(2\pi \xi \cdot x) \geq 0$ for all $\xi \in S$.  Thus
\begin{align*}
\Re \sum_{\xi \in G} |\widehat 1_A(\xi)|^4 e(\xi \cdot x) &=
\sum_{\xi \in S} |\widehat 1_A(\xi)|^4 \cos(2\pi \xi \cdot x) + \sum_{\xi \not \in S} |\widehat 1_A(\xi)|^4 \cos(2\pi i \xi \cdot x) \\
&\geq |\widehat 1_A(0)|^4 - \sum_{\xi \not \in S} |\widehat 1_A(\xi)|^4 \\
&\geq \delta^4 - \sum_{\xi \in G} \alpha^2 |\widehat 1_A(\xi)|^2 \\
&\geq \delta^4 - \alpha^2\delta
\end{align*}
by another application of Plancherel's theorem.  Thus if we set $\alpha := \delta^{3/2}/\sqrt{2}$ (for instance),
the claim follows.  
\end{proof}\vs

\ni We apply Lemma \ref{bog-lemma} with $G = \Ffiven$ and $A = H''$, the set coming from Lemma \ref{random-slice}. The set $S$ we produce satisfies $|S| \leq 50 K^6$. Let $V \subseteq  \Ffiven$
be the subspace $V := \{ x \in \mathbb{F}_5^{n+m}: x \cdot \xi = 0 \hbox{ for all } \xi \in S \}$.
Then by linear algebra we have
\begin{equation}\label{diml}
\dim(V) \geq n - |S| \geq n - 50K^6
\end{equation}
and, since clearly $V \subseteq B(S,\frac{1}{4})$, we have
$$ V \subseteq  2H''-2H''.$$
Thus there exists a map $M: V \to \widehat G$ such that \[ Z := \{ (h, 2Mh): h \in V \}\] is a subset of  $2\Gamma''-2\Gamma''$; since $2\Gamma''-2\Gamma''$ is a graph
and contains the origin, we have $M0=0$.  Also, since $4\Gamma''-4\Gamma''$ is a graph and $V$ is closed under addition, we see that $M(h+h') = Mh + Mh'$ for all $h,h' \in V$, thus $M$ is linear.\vs

\ni Consider the set $Z + \Gamma''$.  On one hand, this set can be foliated into (say) $L$ disjoint cosets of the linear space $Z$.
On the other hand, it is contained in $3\Gamma'' - 2\Gamma''$, 
and \[ |3 \Gamma'' - 2\Gamma''| \leq |9 \Gamma' - 8 \Gamma'| \leq KN.\]
This gives the bound
$|Z| \leq KN/L$.  Since $Z + \Gamma''$ also contains $\Gamma''$, we thus see from the pigeonhole principle and \eqref{eq6.03} that there exists a coset $(x_0,\xi_0) + Z$ of $Z$ such that
\[ |\Gamma'' \cap ((x_0,\xi_0) + Z)| \geq |Z|/5K^4.\]
Since $\dim(Z) = \dim(V) \geq n - 50K^6$, the claim \eqref{hxixi} follows.
\end{proof}\vs

\ni Combining Proposition \ref{quad-linear} with Proposition \ref{add-quad-2}, we quickly conclude the following proposition.

\begin{proposition}[Large $U^3(\Ffiven)$ gives linear phase derivative]\label{add-quad-3}  Let $G = \Ffiven$, and let $f: G \to \D$ be a bounded function such that $\|f\|_{U^3(G)} \geq \eta$ for some $\eta > 0$.  Then there exists a linear subspace $V$ of $G$ with the codimension bound
\begin{equation}\label{ndimv}
 n - \dim(V) \leq 2^{C_1}\eta^{-C_1'}
 \end{equation}
and a translate $x_0 + V$ of this subspace, together with a linear transformation $M: V \to \widehat G$ and an element $\xi_0 \in \widehat G$ such that
\begin{equation}\label{tff}
 \E_{h \in V}  |\E_{x \in G} T^{x_0+h} f(x)  \overline{f(x)} e(- (2Mh + \xi_0) \cdot x)|  \geq 2^{-C_2} \eta^{C'_2}.
\end{equation}
It is permissible to take $C_i,C'_i = 2^{16}$ for $i = 1,2$.
\end{proposition}

\ni\textit{Step 2: The symmetry argument.} We now establish some symmetry properties on $M$, closely following the argument in Lemma \ref{quadratic-classify}.  As we shall be focusing more on $M$ than on $f$ in this step, we shall suppress the terms
involving $f$ (and $\xi_0$ and $x_0$) using the $\Bounded$ notation.  Indeed from \eqref{tff} we have
\begin{equation}\label{tff-1}
|\E_{x \in G, h \in V} \Bounded(x+h) \Bounded(x) \Bounded(h) e(- 2Mh \cdot x)| \geq 2^{-C_2}\eta^{C'_2}
\end{equation}
Once again we use Cauchy-Schwarz and similar tools to eliminate all the bounded functions.
By Lemma \ref{avg-1} there exists $x_1 \in G$ such that
\begin{equation}\label{tff-2}
|\E_{x,h \in V} \Bounded(x+x_1+h) \Bounded(x+x_1) \Bounded(h) e(-2Mh \cdot (x+x_1) )| \geq 2^{-C_2}\eta^{C'_2},
\end{equation}
which after redefining the bounded functions to absorb the $x_1$ terms implies that
\begin{equation}\label{tff-3}
|\E_{x,h \in V} \Bounded(x+h) \Bounded(x) \Bounded(h) e(-2Mh \cdot x ) | 
 \geq 2^{-C_2}\eta^{C'_2}.
\end{equation}
Applying Cauchy-Schwarz (Lemma \ref{cz}) to eliminate $\Bounded(h)$, we see that
\begin{equation}\label{tff-4}
|\E_{x,y,h \in V} \Bounded(y+h) \Bounded(y+h) \Bounded(x+h) \Bounded(x) e(-2Mh \cdot (y-x))| \geq 2^{-2C_2}\eta^{2C'_2};
\end{equation} Making the substitution $z=x+y+h$, this becomes
\begin{equation}\label{tff-5}
|\E_{x,y,z \in V} \Bounded(z,y) \Bounded(z,x) e(-2M(z-x-y) \cdot (y-x) )| \geq 2^{-2C_2}\eta^{2C'_2};
\end{equation} Absorbing as many phase terms into the functions $\Bounded$ as we can, we infer
\begin{equation}\label{tff-5a}
|\E_{x,y,z \in V} \Bounded(z,y) \Bounded(z,x) e(2 \{x,y\} )| \geq 2^{-2C_2}\eta^{2C'_2},
\end{equation}
where $\{x,y\}$ is the anti-symmetric form defined in \eqref{xy-def}, that is to say $\{x,y\} = Mx \cdot y - My \cdot x$.  
By the pigeonhole principle in $z$ we conclude that
\begin{equation}\label{tff-6}
|\E_{x,y \in V} \Bounded(y) \Bounded(x) e(2 \{x,y\} )| \geq 2^{-2C_2}\eta^{2C'_2}
\end{equation}
for some bounded functions $\Bounded(y)$, $\Bounded(x)$.  
Applying Cauchy-Schwarz again to eliminate the $\Bounded(x)$ factor, we deduce that 
\begin{equation}\label{tff-7}
|\E_{x,y,y' \in V} \Bounded(y) \Bounded(y') e(2 \{x,y'-y\} )| \geq 2^{-4C_2}\eta^{4C'_2}.
\end{equation}  By the triangle inequality we obtain
\begin{equation}\label{tff-8}
\E_{y,y' \in V} |\E_{x \in V} e(2\{x,y'-y\})| \geq 2^{-4C_2}\eta^{4C'_2};
\end{equation}
Making the substitution $h = y'-y$ we conclude
\begin{equation}\label{tff-9}
\E_{h \in V} |\E_{x \in V} e(2\{x,h\})| \geq 2^{-4C_2}\eta^{4C'_2}.
\end{equation}
The map $x \mapsto 2\{x,h\}$ is a homomorphism from $V$ to $\R/\Z$.  Thus if we write
$$ W := \{ h \in V: \{x,h\} = 0 \hbox{ for all } x \in V \}$$
then $W$ is a linear subspace of $V$ and 
$$ \E_{x \in V} e(2\{x,h\}) = 1_W(h).$$
Thus, in view of \eqref{tff-9}, we see that $W$ is extremely large relative to $V$:
$$ |W|/|V| = \E_{y,h \in V} 1_W(h) \geq 2^{-4C_2}\eta^{4C'_2}.$$
In particular, from \eqref{ndimv} we have
$$ n - \dim(W) \leq 2^{C_1} \eta^{-C'_1} + 4C'_2 \log_5 (2/\eta) + 4C_2 \log_5 2 \leq 2^{C_1 + 1}\eta^{-C'_1}.$$
By construction of $W$ we see that $M$ is self-adjoint on $W$, that is to say
\begin{equation}\label{self-adjoint}
M w \cdot w' = Mw' \cdot w \hbox{ for all } w,w' \in W.
\end{equation}
Let us remark that the analogy between this argument and that of \S \ref{model-sec} is exceptionally close. The seven equations \eqref{phase-1}, \eqref{phase-2}, \eqref{phase-3}, \eqref{phase-4}, \eqref{phase-5}, \eqref{phase-6} and \eqref{phase-7} are analogues of \eqref{tff-3}, \eqref{tff-4}, \eqref{tff-5}, \eqref{tff-5a}, \eqref{tff-6}, \eqref{tff-7} and \eqref{tff-8} respectively.\vs

\ni\textit{Step 3: Eliminating the quadratic phase component.} We now give the final part of the proof of Theorem \ref{mainthm1}, which also follows the proof of Lemma \ref{quadratic-classify}
closely. We return to \eqref{tff}, which we write as
$$ |\E_{x \in G; h \in V} \Bounded(h) \Bounded(x+h) \overline{f(x)} e(-2Mh \cdot x)| \geq 2^{-C_2}\eta^{C'_2},$$
where we have distributed the phase factor $e(- \xi_0 \cdot x) = e(-\xi_0 \cdot (x+h)) e(\xi_0 \cdot h)$ among
the functions $\Bounded$.  
Here the focus will be on the $\overline{f(x)}$ factor, the aim being to demonstrate that this function exhibits some quadratic bias. We first observe that the simple averaging argument of Lemma \ref{avg-1} allows us to find $h' \in V$ such that
$$ |\E_{x \in G; h \in W} \Bounded(h+h') \Bounded(x+h+h') \overline{f(x)} e(-2M(h+h') \cdot x)| \geq 2^{-C_2}\eta^{C'_2}.$$
Once again we can absorb $e(-2Mh' \cdot x)$ into the functions $\Bounded$, and conclude that
$$ |\E_{x \in G; h \in W} \Bounded(h) \Bounded(x+h) \overline{f(x)} e(-2Mh \cdot x)| \geq 2^{C_2}\eta^{C'_2},$$
which implies that 
$$ |\E_{y \in G; x,h \in W} \Bounded(h) \Bounded(x+h+y) \overline{f(x+y)} e(-2Mh \cdot (x+y))| \geq 2^{C_2}\eta^{-C'_2}.;$$
Using the triangle inequality we deduce
\begin{equation}\label{eq6.17a} \E_{y \in G} |\E_{x,h \in W} \Bounded(h,y) \Bounded(x+h,y) \overline{f(x+y)} e(-2Mh \cdot x)| \geq 2^{C_2}\eta^{-C'_2}.\end{equation}
Now we observe from \eqref{self-adjoint} that we have the identity 
\[ 2Mh \cdot x = M(x + h) \cdot (x + h) - Mx \cdot x - Mh \cdot h,\]
and hence
$$ e(-2Mh \cdot x) = \Bounded(x+h) \Bounded(h) e(Mx \cdot x).$$
Therefore \eqref{eq6.17a} implies that 
$$ \E_{y \in G} |\E_{x,h \in W} \Bounded(h,y) \Bounded(x+h,y) \overline{f(x+y)} e(-Mx \cdot x)| \geq 2^{-C_2}\eta^{C'_2}.$$
Applying Lemma \ref{fgh-bohr} for each $y \in G$ separately and with $B = B' = W$, we conclude that
$$ \E_{y \in G} \| \overline{f(x+y)} e(Mx \cdot x) \|_{u^2(W)} \geq 2^{-C_2}\eta^{C'_2} ,$$ and hence of course
$$ \E_{y \in G} \| \overline{f(x+y)} e(Mx \cdot x) \|_{u^3(W)} \geq 2^{-C_2}\eta^{C'_2}.$$
But the $u^3(W)$ norm is invariant under quadratic phase modulations, conjugation, and translation, and so we have
$$ \E_{y \in G} \| f \|_{u^3(y+W)} \geq 2^{-C_2}\eta^{C'_2}$$
which gives \eqref{efphi}.  
\endproof

\section{Application: Szemer\'edi's theorem in finite field geometries}\label{finite-sz}

\ni As a sample application of Theorem \ref{mainthm1}, we can now prove a quantitative Szemer\'edi theorem for progressions of length four in $\Ffiven$. If $G$ is any finite abelian group of order $N$, where $(6,N) = 1$, we define $r_4(G)$ to be the cardinality of the largest set $A \subseteq G$ which does not contain four distinct elements in arithmetic progression.

\begin{theorem}[Szemer\'edi theorem for $r_4(\Ffiven)$]\label{sz4-ff} Write $N = 5^n$. Then we have the bound
\[ r_4(\Ffiven) \ll N(\log \log N)^{-2^{-21}}.\]
\end{theorem}
\ni\textit{Remark.} In \S \ref{sz-g} below we will prove a similar result for an arbitrary $G$. Although that result will supersede the present one, the proof is quite a bit more complicated, and so we give the finite field argument separately now. In \cite{green-tao-ffszem} we improve the bound to $r_4(\Ffiven) \ll N(\log N)^{-c}$ using substantially lengthier arguments.\vs

\ni As with all analytic arguments for proving Szemer\'edi type theorems, the key step is the establishment of the following ``density increment'' result.

\begin{proposition}\label{sec7-dens-inc} Let $\delta > 0$, suppose that 
\begin{equation}\label{n-lower-bd}
n > 6(2/\delta)^{2^{20}},
\end{equation}and let $A \subseteq \Ffiven$ be a set with size at least $\delta N$. Suppose that $A$ contains no four-term arithmetic progressions. Then we can find an affine subspace $x_0 + V$ of $\Ffiven$ with dimension $\dim(V) \geq n/3$ such that we have the density increment
\[ \E_{x \in x_0 + V} 1_A(x) \geq \E_{x \in \Ffiven} 1_A(x) + (\delta/2)^{2^{20}}.\]
\end{proposition}
\begin{proof} Write 
\[ \alpha := \E_{x \in \Ffiven} 1_A(x) = |A|/N.\] By Corollary \ref{prog-unif} and the lower bound $N \geq 2/\delta^3$ (which is very much a consequence of \eqref{n-lower-bd}!) we conclude that
$$ \| 1_A - \alpha \|_{U^3(\Ffiven)} \geq \delta^3/8.$$
Applying Theorem \ref{mainthm1}, we can thus find a subspace $W \leq \Ffiven$ with codimension
at most $(2/\delta)^{4C}$ and quadratic phase functions $\phi_y$ for each $y \in \Ffiven$ such that
\begin{equation}\label{eaea}
\E_{y \in \Ffiven} \big| \E_{x \in y + W} (1_A(x) - \alpha)e(-\phi_y(x)) \big| \geq (\delta/2)^{4C}. 
\end{equation}
On the other hand, we may assume that
$$ \E_{x \in y + W} (1_A(x) - \alpha) < \quarter (\delta/2)^{4C}$$
for all $y \in \Ffiven$, since the proposition is immediate otherwise. 
Now since
$\E_{x \in y + W}( 1_A(x) - \alpha )$ has mean zero, we thus conclude that
$$ \E_{y \in \Ffiven} \big|\E_{x \in y + W}( 1_A(x) - \alpha) \big|  < \half (\delta/2)^{4C}.$$
Subtracting this from \eqref{eaea} and using the pigeonhole principle, we infer that there exists $y \in \Ffiven$ such that
$$ \big| \E_{x \in y + W} (1_A(x) - \alpha) e(-\phi_y(x)) \big| > \half (\delta/2)^{4C} + \big| \E_{x \in y + W} (1_A(x) - \alpha) \big|.
$$
Now observe (using Lemma \ref{quadratic-classify}) that $\phi_y$ takes values in the set $T = \{0,\frac{1}{5},\frac{2}{5},\frac{3}{5},\frac{4}{5}\}$.  For each $t \in T$ let
$S_t$ be the quadratic surface $S_t := \{ x \in y+W: \phi_y(x) = t \}$.  Then by the triangle inequality we have
$$ \big| \E_{x \in y + W} (1_A(x) - \alpha) e(-\phi_y(x)) \big| \leq \sum_{t \in T} \big| \E_{x \in y + W} (1_A(x) - \alpha) 1_{S_t}(x) \big|
$$
whilst
$$ 
\big| \E_{x \in y + W} (1_A(x) - \alpha) \big| \geq - \E_{x \in y + W} (1_A(x) - \alpha) = -\sum_{t \in T} \E_{x \in y + W} (1_A(x) - \alpha) 1_{S_t}(x).
$$
Combining these estimates, and using the pigeonhole principle, we deduce that there exists $t \in T$ such that
$$ 
\big| \E_{x \in y + W} (1_A(x) - \alpha) 1_{S_t}(x) \big| > \half(\delta/2)^{4C} \E_{x \in y + W} 1_{S_t}(x) - \E_{x \in y + W} (1_A(x) - \alpha) 1_{S_t}(x),
$$
and hence
\begin{equation}\label{eaeast}
 \E_{x \in y + W} (1_A(x) - \alpha)1_{S_t}(x) \geq \quarter(\delta/2)^{4C} \E_{x \in y + W} 1_{S_t}(x).
 \end{equation}
This gives a density increment on a quite large quadratic hypersurface $S_t \cap (y+W)$. Our job is now to convert this into a density increment on
a subspace.  Observe from Lemma \ref{quadratic-classify} that we can write $S_t \cap (y+W)$ in the form
$$ S_t \cap (y+W) = y + \{ x \in W : \half x \cdot Mx + v \cdot x = c \}$$
for some self-adjoint linear transformation $M: W \to \widehat W$, some $v \in W^*$, and some $c \in F$.  We now locate a large subspace on which $M$ is degenerate.
To this end we need a simple lemma\footnote{One could also proceed here using the theory of Witt groups, but that would be far more advanced technology than what is actually needed here.}.  For future reference we shall phrase this lemma for more general finite fields than $\mathbb{F}_5$.

\begin{lemma}[Gauss sum lemma]\label{gauss} Let $F$ be a finite field of odd characteristic, let $W$ be a vector space over $F$, and let $M: W \to \widehat W$ be a self-adjoint linear transformation.  If $\dim(W) \geq 3$, then
there exists a non-zero $x \in W$ such that $x \cdot Mx = 0$.
\end{lemma}

\begin{proof} If $M$ has a non-trivial kernel then the claim is easy, so assume $M$ has no kernel.  Then a standard Gauss sum computation (using  \eqref{vdc}, for instance) shows that $|\E_{x \in W} e( a(x \cdot Mx) )|^2 = |F|^{-\dim(W)}$ for all $a \in |F| \setminus \{0\}$.  But by Fourier inversion we have
$$ \E_{x \in W} 1_{x \cdot Mx = 0}  = \textstyle\frac{1}{|F|}\displaystyle \sum_{a \in F} \E_{x \in W} e( a(x \cdot Mx) ) \geq \textstyle\frac{1}{|F|}\displaystyle (1 - (|F|-1) \cdot  |F|^{-\dim(W)/2}).$$
Since $\dim(W) \geq 3$ we thus see that $\{ x \in W: x \cdot Mx = 0 \}$ must contain at least one non-zero element, and we are done.
\end{proof}\vs

\ni In our situation, the space $W$ has dimension substantially larger than 3 -- in fact $\dim(W) \geq n - (2/\delta)^{4C}$.  Let $U$ be a subspace of $W$
which is degenerate in the sense that $x \cdot My = 0$ for all $x,y \in U$, and which is maximal with respect to this property.
We claim that $\dim(U) \geq n/2 - (2/\delta)^{4C}$. Indeed if this were not the case the space $U^\perp := \{ x \in W: x \cdot My = 0 \hbox{ for all } y \in U \}$
would have at least three more dimensions than $U$ (in fact, vastly more than this), and one could apply the previous lemma to $U^\perp/U$ to contradict the maximality of $U$.\vs

\ni Let $U$ be as above.  By splitting \eqref{eaeast} into cosets of $U$ and applying the pigeonhole principle, we may find a coset $z+U$ of $U$ in
$y+W$ such that
$$ \E_{x \in z + U}( 1_A(x) - \alpha ) 1_{S_t}(x) \geq \quarter (\delta/2)^{4C}\E_{x \in z + U}1_{S_t}(x).$$
Note that on this space $z+U$, the form $x \cdot Mx$ becomes linear, which means that $S_t \cap (z+U)$ is an affine subspace of $z+U$ of codimension
at most $1$.  Letting $x_0+V$ denote this subspace, and recalling that the constant $C$ in Theorem \ref{mainthm1} could certainly be taken to be $2^{17}$, Proposition \ref{sec7-dens-inc} follows. Note that \eqref{n-lower-bd} suffices to guarantee the stated bound $\dim(V) \geq n/3$.\endproof\vs

\ni\textit{Proof of Theorem \ref{sz4-ff}.} Suppose that $A \subseteq \Ffiven$ has cardinality $\delta N$ yet contains no four-term arithmetic progression. Then by repeated application of Proposition \ref{sec7-dens-inc} we may construct a sequence
\[ \mathbb{F}_5^n \geq x_1 + V_1 \geq x_2 + V_2 \geq \dots\]
of affine subspaces such that $\dim(V_j) \geq n/3^j$ and
\begin{equation}\label{di} \E_{x \in x_j + V_j} 1_A(x) \geq \delta + j(\delta/2)^{2^{20}},\end{equation} provided only that at all stages the condition \eqref{n-lower-bd} is satisfied, that is to say
\begin{equation}\label{end-iteration}
\dim(V_j) > 6(2/\delta)^{2^{20}}.
\end{equation}
Equation \eqref{di} is impossible if $j \geq j_0 = (2/\delta)^{2^{20}}$, and so \eqref{end-iteration} must be violated by some $j \leq j_0$. This means that 
\[ n > 6 \times (2/\delta)^{2^{20}} \times 3^{(2/\delta)^{2^{20}}},\]
which certainly implies that $\delta \ll (\log \log N)^{-2^{-21}}$.
\end{proof}\vs

\ni\textit{Remark.} A more-or-less identical argument shows that 
\[ r_4(\mathbb{F}_p^n) \ll N(\log \log N/p)^{-c},\] uniformly for all primes $p \geq 5$, for some absolute constant $c > 0$. In combination with Gowers' result that $r_4(\Z/N\Z) \ll N(\log \log N)^{-c}$, this may be used to give a fairly cheap proof that $r_4(G) = o(N)$ for all finite abelian $G$, a result which was first obtained by Frankl and R\"odl \cite{fr} by rather different means. The key to the argument is that $G$ contains either a large subgroup of the form $\mathbb{F}_p^n$, or else a large cyclic subgroup. The bound obtained is of the form 
\[ r_4(G) \ll N(\log \log\log \log N)^{-c};\]
we suppress the details since a far superior bound will be obtained in \S \ref{sz-g}.

\section{Some results on Bohr sets}\label{sec8}

\ni Let $G$ be an arbitrary finite abelian group with $|G| = N$, let $S \subseteq \widehat{G}$ will be a set of $d$ characters, and suppose that $\rho \in (0,1)$ is a positive parameter. We will collect some basic facts about Bohr sets $B(S,\rho)$ which we will
need to prove Theorem \ref{mainthm2}.  These Bohr sets play the role that subspaces did in the finite geometry setting, with the quantity $d$ corresponding, roughly speaking, to the codimension of the subspace.\vs

\ni Note that a Bohr set always contains 0, and is symmetric around the origin.  In fact we have the following easy bounds on the size of Bohr sets:

\begin{lemma}[Bounds for size of Bohr sets]\label{bohr-bound} We have
$$ |B(S,\rho)| \geq \rho^{d} N$$
and
$$ |B(S,2\rho)| \leq 5^{d} |B(S,\rho)|.$$
\end{lemma}
\begin{proof}  By the triangle inequality, we see that for any $y = (y_\xi)_{\xi \in S}$ in the torus $(\R/\Z)^S$, we have
\begin{align*}
\sum_{x \in G} \prod_{\xi \in S} 1_{\| \xi \cdot x - y_\xi \|_{\R/\Z} < \rho/2} &=
|\{ x \in G: \sup_{\xi \in S} \| \xi \cdot x - y_\xi \|_{\R/\Z} < \rho/2 \}| \\
& \leq
 |\{ x \in G: \sup_{\xi \in S} \| \xi \cdot x \|_{\R/\Z} < \rho \}| \\
 &= |B(S,\rho)|.
\end{align*}
Integrating this over all $y \in (\R/\Z)^S$, we conclude that
$$
\sum_{x \in G} \rho^{d} \leq |B(S,\rho)|
$$
which gives the first bound.\vs

\ni To establish the second bound, we integrate the same expression but only over the cube $\{ y: \sup_{\xi \in S} \| y_\xi \| \leq \frac{5}{2} \rho \}$.
Note from the triangle inequality that this does not affect the components of the integral for which $x \in B(S, 2\rho)$.  Thus we have
$$ \sum_{x \in B(S,2\rho)} \rho^{d} \leq |B(S,\rho)| (5\rho)^{d}$$
which gives the second bound.
\end{proof}\vs

\ni Next, we establish, following Bourgain, that regular Bohr sets (as defined in Definition \ref{regular-def}) exist in abundance.

\begin{lemma}[Regular Bohr sets are ubiquitious \cite{Bou}]\label{make-regular}  Let $0 < \eps < 1$.  Then there exists $\rho \in [\eps,2\eps]$ such that $B(S,\rho)$ is regular.
\end{lemma}

\begin{proof}  We may assume $S$ is non-empty since the claim is trivial otherwise.
Let $f: [0,1] \to \R$ be the function $f(a) := \frac{1}{d} \log_2 |B(S, 2^a \eps)|$.  Observe that $f$ is non-decreasing in $a$,
and from Lemma \ref{bohr-bound} we have $f(1) - f(0) \leq \log_2 5$. \vs

\ni Suppose we could find $0.1 \leq a \leq 0.9$ is such that $|f(a') - f(a)| \leq 20 |a-a'|$ for all $|a| \leq 0.1$.  Then it
is easy to see that the Bohr set $B(S, 2^a \eps)$ is regular.  Thus, it suffices to obtain an $a$ with this property.  This can be done
directly from the Hardy-Littlewood maximal inequality (applied to the Lebesgue-Stieltjes measure $df$), or as follows.  If no such $a$ exists, then for every $a \in [0.1,0.9]$ there exists
an interval $I$ of length at most $0.1$ and with one endpoint equal to $a$, such that $\int_I df > \int_I 20\ dx$.
These intervals cover $[0.1,0.9]$, which has measure $0.8$.  By the Vitali covering lemma\footnote{One can also use the Besicovitch covering lemma at this point, which would in fact give slightly better bounds.  Indeed, one can improve the constant $1/5$ to $1/2$, see for instance \cite{Croft}.}, one can find thus find a finite subcollection of disjoint intervals $I_1,\ldots,I_n$ of total length $|I_1| + \ldots |I_n| \geq 0.8/5$ (say).
But then we have
$$ \log_2 5 \geq \int_0^1 df \geq \sum_{i=1}^n \int_{I_i} df \geq \sum_{i=1}^n \int_{I_i} 20 dx \geq \frac{0.8}{5} \times 20,$$
a contradiction.
\end{proof}\vs

\ni Next, we show that given any small set of points in $G$, one can find a large Bohr set which avoids all of them except possibly for zero.
We will need this to develop the analogue of Lemma \ref{random-slice} (see Lemma \ref{more-random-slice} below).

\begin{lemma}[Separation lemma]\label{separation}
Let $G$ be a finite additive group, and let $A \subseteq  G$ be a set of elements containing zero.  Then there
exists a set $S \subseteq \widehat G$ with $|S| \leq 1 + \log_2 |A|$ such that $A \cap B(S, \frac{1}{4}) = \{0\}$.
\end{lemma}

\begin{proof}  We induct on $|A|$.  When $|A| = 1$ the claim is trivial (set $S = \emptyset$).  Now suppose that $|A| \geq 2$ and the claim has already
been proven for smaller sets $A$.  Suppose that $2^n < |A| \leq 2^{n+1}$.  Let $\xi \in \widehat G$ be chosen randomly.  Observe that for each $x \in A \backslash \{0\}$, the map $\xi \mapsto \xi \cdot x$
is a non-trivial group homomorphism from $\widehat G$ to $\R/\Z$, thus the random variable $\xi \cdot x$ is uniformly distributed over a cyclic subgroup of $\R/\Z$.  In particular, we have 
$$\P( x \in B(\xi, \quarter )) = \P( \| \xi \cdot x \|_{\R/\Z} < \quarter ) \leq \half.$$
Summing this over all non-zero $A$, we conclude
$$ \E|(A \backslash \{0\}) \cap B(\xi,\quarter)| \leq (|A|-1)/2.$$
In particular, we can find $\xi \in \widehat G$ such that
$$ |A \cap B(\xi,\quarter)| \leq \lfloor \frac{|A|+1}{2} \rfloor \leq 2^n.$$
By induction hypothesis we conclude that there exists a set $S'$ of cardinality at most $n$ such that
$$ A \cap B(\xi,\quarter) \cap B(S',\quarter) = \{0\},$$
and the claim follows by setting $S := S' \cup \{\xi\}$.
\end{proof}\vs

\ni Our next task is to investigate the Fourier-analytic behavior of Bohr sets.  In the first instance we deal with a concept somewhat more general that that of a Fourier coefficient, replacing a linear phase function by a locally linear phase function.

\begin{lemma}[Generalized Fourier decay]\label{gfd}  Let $S \subseteq \widehat{G}$, $|S| = d$, be a set of characters. Let $B := B(S,\rho)$ be a regular Bohr set, and let $\phi: B(S,2\rho) \to \R/\Z$ be a function which is locally linear in the sense that that $\phi(x+y) = \phi(x)+\phi(y)$ whenever $x,y \in B(S,\rho)$.  
Suppose that 
\[ |\E_{x \in B} (e (\phi(x)))| \geq \eta\]
for some $0 < \eta \leq 1$ \textup{(}large generalized Fourier coefficient\textup{)}.  Then $\phi$ is close to constant in the sense that for every $y \in B$ we have
$$ \| \phi(y) \|_{\R/\Z} \leq \frac{2^{12}d\| y \|_S}{\rho \eta^2}.$$
\end{lemma}

\begin{proof}  Let $y \in B$, and let $M \geq 0$ be the largest integer such that $M\|y\|_S \leq \eta \rho / 400 d$.
If $M = 0$ then we have $\|y\|_S \geq \eta \rho / 400 d$, and the claim is trivial.  By Lemma \ref{avg} (ii) we have
$$ | \E_{x \in B}( e(\phi(x)) ) - \E_{x \in B} \E_{-M \leq n \leq M} e(\phi(x + ny))) | \leq \eta/2$$
and hence by the triangle inequality
$$ | \E_{x \in B} \E_{-M \leq n \leq M} e(\phi(x + ny)) | \geq \eta/2.$$
On the other hand, by the local linearity of $\phi$ we have $\phi(x+ny) = \phi(x) + n \phi(y)$ for all $|n| \leq M$
and hence $e(\phi(x+ny)) = \Bounded(x) e(n\phi(y))$. By the triangle inequality we conclude that
$$ |\E_{-M \leq n \leq M} e(n \phi(y))| \geq \eta/2.$$
But by the geometric series formula, the left-hand side is bounded by $2 / M \| \phi(y) \|_{\R/\Z}$.  This implies that
$$ \| \phi(y) \|_{\R/\Z} \leq \frac{4}{M\eta} \leq \frac{3200 d}{\rho \eta^2}\| y \|_S,$$
which implies the result.
\end{proof}\vs

\ni As a corollary we see that the normalized Fourier transform of a Bohr set decays away from the ``polar body'' of that Bohr set\footnote{One could obtain much better Fourier localization properties by replacing the Bohr sets by smoother weight functions; see for instance \cite{green1,tao1} for examples of this approach.  This also conveys the slight advantage that all weight functions can automatically be made regular.  However these functions have the disadvantage of being spread out in physical space, and we found it more convenient to use Bourgain's machinery of regular Bohr sets from \cite{Bou} instead.}. If $\xi \in \widehat{G}$ then define $$ \| \xi \|_{B(S,\rho)} := \sup_{y \in B(S,\rho)} \| \xi \cdot y \|_{\R/\Z}.$$ Note that if $\xi \in S$ then $\| \xi \|_{B(S,\rho)} \leq \rho$. 

\begin{corollary}[Fourier decay]\label{fd} Let $S \subseteq \widehat{G}$ be a set of $d$ characters, let $B := B(S,\rho)$ be a regular Bohr set, and let $0 < \theta \leq 1$.  Then for any $\xi \in \widehat G$, we have
$$ |\E_{x \in B} e(\xi \cdot x)| \leq 64\big(\frac{\theta d}{\| \xi \|_{B(S,\theta \rho)}}\big)^{1/2}$$
where 

\end{corollary}

\begin{proof} Apply Lemma \ref{gfd} with $\phi(x) := \xi \cdot x$, with
$\eta := |\E_{x \in B} e(\xi \cdot x)|$
and with $y$ being an arbitrary element of $B(S,\theta r)$.
\end{proof}\vs

\ni We can exploit this decay via a Tomas-Stein almost-orthogonality type argument (also used by Bombieri \cite{bombieri} in the context of the large sieve; see also \cite{Mont}) to conclude

\begin{corollary}[Local Bessel inequality]\label{lbi}   Let $S \subseteq \widehat{G}$ be a set of $d$ characters, 
let $B := B(S,\rho)$ be a regular Bohr set, let $0 < \theta \leq 1$, and let $\xi_1,\ldots,\xi_k\in \widehat G$
be frequencies such that $\| \xi_i - \xi_j \|_{B(S,\theta \rho)} \geq \delta$ for all $1 \leq i < j\leq k$ and some $\delta > 0$.  Then
$$ \E_{x \in B} |\sum_{j=1}^k\Bounded(j) e(\xi_j(x))|^2 \leq k + 2^7 k^2 (\theta d / \delta)^{1/2}$$
for any bounded complex numbers $\Bounded(j)$.
\end{corollary}

\begin{proof} We have
\begin{align*}
\E_{x \in B} |\sum_{j=1}^k \Bounded(j) e(\xi_j(x))|^2 
&= \sum_{1 \leq i,j \leq k} \Bounded(i,j) \E_{x \in B} e((\xi_i - \xi_j) \cdot x)  \\
&\leq k + 2 \sum_{1 \leq i < j \leq k} |\E_{x \in B} e(\xi_i - \xi_j) \cdot x)| \\
&\leq k + 2^7 k^2 (\theta d / \delta)^{1/2},
\end{align*}
thanks to Corollary \ref{fd}.  The claim follows.
\end{proof}\vs

\ni We can dualize the above corollary to give the following result. This allows us to generalize, to the relative setting, a frequently-used consequence of Parseval's identity: a large set $A \subseteq G$ cannot have too many large Fourier coefficients.

\begin{corollary}[Local Bessel inequality, dual version]\label{lbi-dual} Let $S \subseteq \widehat{G}$ be a set of $d$ characters, let $B := B(S,\rho)$ be a regular Bohr set, let $0 < \theta, \eta \leq 1$, and suppose that $A \subseteq B$. Let
\[ \Gamma := \{ \xi \in \widehat G : |\widehat{1_A}(\xi)| \geq \eta \E (1_B)\}.\]
Then there exist frequencies $\xi_1,\ldots, \xi_k \in\widehat G$ with $k \leq 2/\eta^2$ such that any $\xi \in \Gamma$ is close to some $\xi_i$ in the $\| \cdot \|_{B(S,\theta \rho)}$ norm:
\begin{equation}\label{xihat}\Gamma \subseteq  \{ \xi \in \widehat G: \| \xi - \xi_j \|_{B(S,\theta \rho)} \leq 2^{16}\theta d / \eta^4
\hbox{ for some } 1 \leq j \leq k \}.
\end{equation}
\end{corollary}

\begin{proof}  Let $\delta := 2^{16} \theta d / \eta^4$, and let $\xi_1,\ldots,\xi_k$ be frequencies in $\Gamma$
such that $\|\xi_i - \xi_j\|_{B(S,\theta \rho)} \geq \delta$, and which is maximal with respect to set inclusion.  Then it is clear that
\eqref{xihat} holds.  For each $1 \leq j \leq k$, we have $\xi_j \in \Gamma$.  Hence there exists a bounded complex number $\Bounded(j)$ such that
$$ \Re \E_{x \in B} \Bounded(j) e(\xi_j(x)) 1_A(x) \geq \eta.$$
Summing this in $j$ and applying Cauchy-Schwarz, we conclude that
$$ \E_{x \in B} |\sum_{j =1}^k \Bounded(j) e(\xi_j(x)) |^2 \geq \eta^2 k^2.$$
Applying Corollary \ref{lbi}, we conclude that
$\eta^2 k^2 \leq k + \frac{1}{2} \eta^2 k^2$, and hence $k \leq 2 / \eta^2$.
The claim follows.
\end{proof}\vs

\ni As a consequence, we can now generalize Bogolyubov's argument (Lemma \ref{bog-lemma}) to subsets of Bohr sets.

\begin{lemma}[Local Bogolyubov lemma]\label{bog-lemma-local} Let $S \subseteq \widehat G$ be a set of $d$ characters, and let $B := B(S,\rho)$. Let $A \subseteq  B$ be a set with $|A| = \delta |B|$. Then there exists a set $S' \subseteq  \widehat G$ with $|S'| \leq 2^7 \delta^{-3}$ such that
$B(S \cup S', 2^{-33}\delta^6 \rho/d ) \subseteq  2A-2A$.
\end{lemma}

\begin{proof}  It is convenient to replace $2A-2A$ by the slightly smaller set $A+A'-A-A'$.
Let $\eps = \delta/400d$. By Lemma \ref{make-regular}, there exists $\rho' \in [\eps \rho, 2\eps, \rho]$
such that the Bohr set $B' := B(S,\rho')$ is regular.  By Lemma \ref{avg} (iii) we can find $x \in B$
such that
$$ \E_{y \in x + B'}1_A(y) \geq \E_{y \in B}1_A(y) - 200 d \eps \geq \delta / 2.$$
Let $A' := A \cap (x+B')$.  From the Fourier inversion formulae
$$
1_A(x) = \sum_{\xi \in \widehat G} \widehat 1_A(\xi) e(\xi \cdot x); \quad
1_{A'}(x) = \sum_{\xi \in \widehat G} \widehat 1_{A'}(\xi) e(\xi \cdot x)
$$
we conclude that
\begin{equation}\label{aaaa}
1_A * 1_{A'} * 1_{-A} * 1_{-A'}(x) = \sum_{\xi \in \widehat G} |\widehat 1_A(\xi)|^2 |\widehat 1_{A'}(\xi)|^2 e(\xi \cdot x).
\end{equation}
In particular, applying \eqref{aaaa} with $x=0$ we conclude that
$$ \E_{x \in G} |1_A * 1_{A'}(x)|^2 = \sum_{\xi \in \widehat G} |\widehat 1_A(\xi)|^2 |\widehat 1_{A'}(\xi)|^2.$$
The function $1_A * 1_{A'}$ is supported on $B(S,\rho + \rho')$, which has cardinality at most $2|B|$ since $B$ is regular and $\eps < 1/200 d$.  Thus by Cauchy-Schwarz
\begin{eqnarray*}
\E_{x \in G} |1_A * 1_{A'}(x)|^2  &\geq  & (\E_{x \in G} 1_A * 1_{A'}(x))^2 / 2\E(1_B) \\
&= &   \E( 1_A )^2 \E(1_{A'})^2 / 2\E(1_B) \\
& = & \half \big(\E_{y \in B}1_A(y)\big)^2 \big(\E_{y \in x + B'}1_A(y)\big)^2 \E(1_B) \E(1_{B'})^2 \\
&\geq & \delta^4 K^{-4} \E(1_B) \E(1_{B'})^2/8.
\end{eqnarray*}
Thus we have
\begin{equation}\label{a0}
\sum_{\xi \in \widehat G} |\widehat 1_A(\xi)|^2 |\widehat 1_{A'}(\xi)|^2  \geq  \textstyle \frac{1}{8} \displaystyle \delta^4 \E(1_B) \E(1_{B'})^2.
\end{equation}
Now let \[ R := \{ \xi \in G: |\widehat 1_{A'}(\xi)| \geq \textstyle \frac{1}{8} \displaystyle \delta^{3/2}\E(1_{B'}) \},\] and let $x \in B(R, \frac{1}{10})$.
Then by taking real parts of both sides of \eqref{aaaa}, we conclude
\begin{eqnarray*}
1_A * 1_{A'} * 1_{-A} * 1_{-A'}(x) &= &  \sum_{\xi \in \widehat G} |\widehat 1_A(\xi)|^2 |\widehat 1_{A'}(\xi)|^2 \cos(2\pi \xi \cdot x) \\
&\geq & \sum_{\xi \in R} |\widehat 1_A(\xi)|^2 |\widehat 1_{A'}(\xi)|^2 \cos(2\pi / 10) -
\sum_{\xi \not \in R} |\widehat 1_A(\xi)|^2 |\widehat 1_{A'}(\xi)|^2 \\
& = & \sum_{\xi \in \widehat G}|\widehat 1_A (\xi)|^2 |\widehat 1_{A'}(\xi)|^2 \cos (2\pi /10) \\ & & \qquad - (\cos (2\pi /10) + 1) \sum_{\xi \notin R} |\widehat 1_A(\xi)|^2 |\widehat 1_{A'}(\xi)|^2 \\ 
&\geq & \half (\sum_{\xi \in \widehat G} |\widehat 1_A(\xi)|^2 |\widehat 1_{A'}(\xi)|^2)
- 2 \sum_{\xi \not \in R} |\widehat 1_A(\xi)|^2 |\widehat 1_{A'}(\xi)|^2 \\
&\geq & \textstyle \frac{1}{16} \displaystyle \delta^4 \E(1_B) \E(1_{B'})^2 - \textstyle \frac{1}{32} \displaystyle \delta^3 \E(1_{B'})^2\sum_{\xi \in G} |\widehat 1_A(\xi)|^2
\end{eqnarray*}
using \eqref{a0} and the definition of $R$.  On the other hand, from Plancherel's identity we have
$$ \sum_{\xi \in G} |\widehat 1_A(\xi)|^2  = \E( 1_A) \leq \delta\E(1_B),$$
and hence
$$ 1_A * 1_{A'} * 1_{-A} * 1_{-A'}(x) \geq \big(\frac{\delta^4}{16} - \frac{\delta^4}{32}\big) \E(1_B) \E(1_{B'})^2 > 0,$$
which implies that $x$ is contained in $A+A'-A-A'$ and hence in $2A-2A$.  Hence we have
$$ B(R, \textstyle\frac{1}{10}\displaystyle) \subseteq  2A-2A.$$
We are not done yet, because we do not have good bounds for $|R|$.  Let $\theta > 0$ be a small parameter to be chosen later.
Invoking Corollary \ref{lbi-dual}, we conclude the existence of frequencies $\xi_1,\ldots,\xi_k \in \widehat G$ with $k \leq 128 \delta^3$ such that
$$ R \subseteq  \{ \xi \in \widehat G: \| \xi - \xi_j \|_{B(S,\theta \rho)} \leq 2^{28}\delta^{-6} \theta d  \hbox{ for some } 1 \leq j \leq k \}.$$
Let $S' := \{\xi_1,\ldots,\xi_k\}$.  If $x \in B(S \cup S', \theta \rho)$ then in particular $x \in B(S,\theta \rho)$.  Thus
if $\xi \in R$, then by the preceding inclusion we have
$$ \| \xi \cdot x - \xi_j \cdot x \|_{\R/\Z} \leq 2^{28} \delta^{-6} \theta d \hbox{ for some } j.$$
Also, since $x \in B(S',\theta \rho)$, we get
$$ \| \xi_j \cdot x \|_{\R/\Z} \leq \theta \rho \leq \theta;$$
by the triangle inequality we then obtain
$$ \| \xi \cdot x  \|_{\R/\Z} \leq 2^{29}\delta^{-6} \theta d.$$
We thus conclude that 
$$ B(S \cup S', \theta \rho) \subseteq  B(R, 2^{29} \delta^{-6} \theta d ).$$
Thus if we choose $\theta := 2^{-33}\delta^6/d$, we have
$B(S \cup S', \theta \rho) \subseteq  B(R,\frac{1}{10})$, and since $B(R,\frac{1}{10}) \subseteq 2A - 2A$ the claim follows.
\end{proof}

\section{The general group case}\label{ggc}

\ni We now prove Theorem \ref{mainthm2}, which generalizes Theorem \ref{mainthm1} to the case of arbitrary
finite additive groups $G$.  We begin by disposing of part (ii) of the theorem, which is rather easier to establish than (i).\vs

\ni Recall that in order to establish Theorem \ref{mainthm2} (ii) we are to prove that if $S \subseteq  \widehat{G}$ is a set of $d$ characters, if $B = B(S,\rho)$ is a regular Bohr set, if $f : G \rightarrow \D$ is a bounded function and if $\| f \|_{u^3(y + B)} \geq \eta$ then we have 
\[ \| f \|_{U^3(G)} \geq (\eta^3 \rho^2/Cd^3)^d\]
for some absolute constant $C$. 
By translation invariance we may take $y=0$.  Let $\phi: B \to \R/\Z$ be a locally quadratic phase function
on $B$, and suppose that 
$$ |\E_{x \in B}( f(x) e(-\phi(x)) )| = \eta. $$
It turns out to be convenient to have $\phi$ defined, and to be a quadratic form,  on a slightly larger Bohr set than $B$. This is not in general possible, but the same effect can be achieved by first passing to a smaller Bohr set. In fact, in the argument which follows we will have two smaller Bohr sets $B' = B(S,\rho')$ and $B'' = B(S,\rho'')$. Set $\eps = c\eta/d$, where $c$ is a small constant to be specified later. We will take $\rho' \in [\eps\rho/2,\eps \rho]$ so that $B'$ is regular (this is possible by Lemma \ref{make-regular}) and $\rho'' = \eps \rho'$ (we will not require $B''$ to be regular). It will be convenient to write $\beta' := \E 1_{B'}$ and $\beta'' := \E 1_{B''}$.   
By Lemma \ref{avg} we have
$$ \eta = \E_{z \in B} \E_{x \in z + B'}( f(x) e(-\phi(x))) + O(\eps d ).$$
Observe that the contribution from $z \in B \backslash B(S, (1-10\eps) \rho)$ is at most $O(\eps d)$, thanks to the regularity
of $B$,  Thus we in fact have
$$ \E_{x \in B}(f(x) e(-\phi(x)) ) = \E_{z \in B} 1_{B(S,(1-10\eps)\rho)}(z) \E_{x \in z + B'}( f(x) e(-\phi(x))) + O(\eps d ),$$
and hence by the pigeonhole principle there exists $z \in B(S,(1-10\eps)\rho)$ such that
\begin{equation}\label{efbb}
 |\E_{x \in z + B'}( f(x) e(-\phi(x)) )|  \geq 2\eta/3
 \end{equation}
provided that $c$ is chosen sufficiently small.\vs

\ni We are going to compare $1_{z + B'}$, which is relevant to \eqref{efbb}, with the function
\[ F(x) := \frac{1}{\beta''} \E_h 1_{z + B''}(x + h)1_{z + B' }(x + 2h).\]
Write $B'_{-} = B(S,(1 - 2\epsilon)\rho')$ and $B'_{+} = B(S,(1 + 2\epsilon)\rho')$.
Note that if $x \in z + B'_{-}$ and $x + h \in z + B''$ then $x + 2h = 2(x + h) - x$ is contained in $z + B' $. For such $x$, then, we have $F(x) = 1$. Also, if $F(x) \neq 0$ then there is some $h$ such that $x + h \in z + B''$ and $x + 2h \in z + B'$, which means that $x = 2(x + h) - (x + 2h)$ lies in $z + B'_+$. We have, then, 
\begin{equation}\label{eq9.01} |F(x) - 1_{z + B'}(x)| \leq 2 1_{B'_+ \setminus B'_{-}}(x).\end{equation}
Now note further that if $x + h \in z + B''$ and $x + 2h \in z + B'$ then \[ x \in z + B' + 2B'' \subseteq  B,\] and also 
\[ x + 3h = 2(x + 2h) - (x + h) \in z + 2B' + B'' \subseteq  B.\]
Both of these are consequences of the fact that $z \in B(S,(1 - 10\eps)\rho)$. In such an eventuality, then, all four of the elements $x,x+h,x+2h,x+3h$ lie in $B$ and, since $\phi$ is quadratic, we have 
 $(h \cdot \nabla_x)^3 \phi(x) = 0$, or in other words 
\begin{equation}\label{constraint}
 \phi(x) - 3\phi(x+h) + 3\phi(x+2h) - \phi(x+3h) = 0.
 \end{equation}
Therefore
\begin{eqnarray*} F(x) e(-\phi(x)) & = &  \frac{1}{\beta''}\E_h 1_{z + B''}(x + h) e(3\phi(x + h)) 1_{z + B'}(x + 2h) \times \\ & & \qquad\qquad\qquad\qquad\qquad\qquad \times \; e(-3\phi(x + 2h)) e(\phi(x + 3h))\\ & =: & \frac{1}{\beta''}\E_h g_1(x+h)g_2(x+2h)g_3(x+3h),\end{eqnarray*} say, where the functions $g_1,g_2,g_3$ are all bounded by 1.
It is immediate from \eqref{eq9.01} that 
\[ |F(x) e(-\phi(x)) - 1_{z + B'}(x)e(-\phi(x))| \leq 2 1_{B'_+ \setminus B'_{-}}(x).\]
From \eqref{efbb} we infer, then, that 
\begin{eqnarray*} \frac{1}{\beta''}\E_{x,h} f(x)g_1(x+h) g_2(x+2h)g_3(x + 3h) & = &   \E_x f(x) F(x) e(-\phi(x)) \\ & = & \E_x f(x)1_{z + B'}(x) e(-\phi(x)) + O(\E1_{B'_+ \setminus B'_{-}}) \\ & =  & \E_x f(x)1_{z + B'}(x) e(-\phi(x)) + O(\epsilon d \beta'),\end{eqnarray*}
the penultimate step being a consequence of the regularity of $B'$. If $c$ is chosen small enough, this means in view of \eqref{efbb} that 

\begin{equation}\label{eq9.02} \E_{x,h} f(x)g_1(x+h) g_2(x+2h)g_3(x + 3h) \geq \beta' \beta'' \eta/3.\end{equation}
However Proposition \ref{gvn} implies that we have
\[ \big|\E_{x,h} f(x) g_1(x + h) g_2(x + 2h) g_3(x + 3h)\big| \leq \|f\|_{U^3(G)},\]
and hence from \eqref{eq9.02} we have
\[ \| f \|_{U^3(G)} \geq \beta' \beta'' \eta/3 \geq (\eps^3 \rho^2/4)^d \eta/3 \geq (\eta^3 \rho^2/Cd^3)^d\] for some absolute constant $C$.\endproof\vs

\ni Now we turn to the proof of Theorem \ref{mainthm2} (i).  As in \S \ref{finite-sec} our starting point is Proposition \ref{add-quad-2}, which the reader may care to recall now. The argument is closely analogous to that in \S \ref{finite-sec}, and hence in turn to that in \S \ref{model-sec}.\vs

\ni\textit{Step 1: Linearization of phase derivative.} We begin by carrying out the first major step, which is to show that the function $h \mapsto \xi_h$, which
roughly speaking captures the derivative of the phase of $f$, matches up with a locally linear function.

\begin{proposition}\label{prop9.1a}
Let $H' \subseteq G$, and suppose that $\xi : H' \rightarrow \widehat G$ is a function whose graph 
\[ \Gamma' := \{(h,\xi_h) : h \in H'\} \subseteq G \times \widehat G\]
obeys the estimates
\[ K^{-1} N \leq |\Gamma'| \leq |9\Gamma' - 8\Gamma| \leq KN\]
for some $K \geq 1$. Then there is a set $S \subseteq \widehat G$, 
\[ d_1 := |S| \leq 2^{13}K^{26},\] a regular Bohr set $B_1 := B(S,\rho)$, where $\rho \in [\frac{1}{16},\frac{1}{8}]$, elements $x_0 \in G, \xi \in \widehat G$ and a function $M : B(S,\frac{1}{4}) \rightarrow \widehat G$ satisfying the local linearity condition\footnote{recall that $\|h\|_S := \sup_{\xi \in S} \| \xi_h \|_{\R/\Z}$} 
\begin{equation}\label{madd}
M(h \pm h') = Mh \pm Mh'  \hbox{ whenever } \|h\|_S, \|h'\|_S \leq \textstyle\frac{1}{8}\displaystyle,
\end{equation}
and such that 
\[ \E \big( 1_{H'}(x_0 + h) 1_{\xi_{x_0 + h} = 2M h + \xi_0} | h \in B_1\big) \geq 2^{-6} K^{-13}.\]
\end{proposition}

\ni\textit{Proof.} As in the finite field case, the first step is to refine the graph $\Gamma'$ so that certain of the iterated sum-difference sets
$k \Gamma' - l\Gamma'$ are also graphs.  To do this we use the following generalization of Lemma \ref{random-slice}:

\begin{lemma}\label{more-random-slice}  There 
exists a subset $\Gamma'' = \{ (h,\xi_h) : h \in H'' \}$ of $\Gamma'$ with
$$ |\Gamma''| \geq 2^{-6} K^{-13} N$$
such that $4\Gamma'' - 4 \Gamma''$ is a graph.  
\end{lemma}

\begin{proof}  Let $A \subseteq \widehat G$ be the set of all $\xi$ such that $(0,\xi) \in 8\Gamma' - 8\Gamma'$.  Arguing as in the proof of 
Lemma \ref{random-slice} we conclude that $|A| \leq K^2$.  Applying Lemma \ref{separation}, we can find a set $S \in \widehat{\widehat{G}\,}$
with $|S| \leq 1 + 2\log_2 K$ such that $A \cap B(S, \frac{1}{4}) = \{0\}$.\vs

\ni Let $\Psi: \widehat G \to (\R/\Z)^S$ be the homomorphism $\Psi(\xi) := (s(\xi))_{s \in S}$.
Now let us cover the torus $(\R/\Z)^S$ by $2^{6|S|} \leq 2^6 K^{12}$ cubes of side-length $\frac{1}{64}$.  Since $|\Gamma'| \geq N/K$, the pigeonhole principle implies that there exists one of these cubes $Q$ for which the set
$$ \Gamma'' := \{ (h,\xi_h) \in \Gamma': \Psi(\xi_h) \in Q \}$$
has cardinality at least $2^{-6} K^{-13} N$.  Now observe from the linearity of $s$ that if $(0,\xi) \in 8\Gamma''-8\Gamma''$ then
$\| s(\xi)\|_{\R/\Z} \leq \frac{16}{64}$ for all $s \in S$.  In other words, $\xi \in B(S, \frac{1}{4})$.  But $\xi$ also lies in $A$,
and hence $\xi=0$ by construction.  Since $8\Gamma'' - 8\Gamma''$ is the difference set of $4\Gamma'' - 4\Gamma''$, we conclude
that $4\Gamma'' - 4\Gamma''$ is a graph as desired.
\end{proof}\vs

\ni Define $H''$ so that $\Gamma'' = \{(h,\xi_h) : h \in H''\}$. Applying Lemma \ref{bog-lemma} with $A := H''$, we obtain a set
$S \subseteq \widehat G$ with $|S| \leq 2^{13}K^{26}$ such that the Bohr set $B_0 := B(S, \frac{1}{4})$ is completely contained
inside $2H''-2H''$.  We will now work inside this Bohr set $B_0$ and pass to progressively narrower Bohr sets $B_1, B_2, \ldots$
when necessary. We will eventually end up at $B_5$; the dimension of $B_j$ will be denoted $d_j$, and so in particular $d_0 = |S|$. It will turn out that $d_0 = d_1 = d_2 < d_3 = d_4 = d_5$, that is to say it is only in passing from $B_2$ to $B_3$ that we shall increment the dimension of $B_j$. This is because that passage will involve Lemma \ref{bog-lemma-local}.\vs

\ni Since $2\Gamma''-2\Gamma''$ is a graph, we can find a (unique) function $M: B_0 \to \widehat G$ such that
$$ \{ (h, 2M(h)): h \in B_0 \} \subseteq  2H'' - 2H''.$$ 
Since $2\Gamma''-2\Gamma''$ contains 0, we conclude that $\phi(0) = 0$.  Also, since $8\Gamma''-8\Gamma''$ is a graph we see that
\begin{equation}\label{mfreiman}
M(h_1) + M(h_2) = M(h'_1) + M(h'_2)
\end{equation}
whenever $h_1,h_2, h'_1,h'_2 \in B_0$ is such that $h_1+h_2 = h'_1 + h'_2$;
in other words, $M$ is a \emph{Freiman homomorphism} of order 2.  In particular, since $M0 = 0$, we have the local linearity relationship \eqref{madd}.\vs

\ni By Lemma \ref{make-regular} there is $\rho \in [\frac{1}{16},\frac{1}{8}]$ such that
the Bohr set $B_1 := B(S, \rho_1)$ is regular.  By Lemma \ref{avg-1} there exists $x_0 \in G$ such that
$$ \E_{h \in B_1} 1_{H''}(x_0+h) \geq 2^{-6} K^{-13}.$$
Let us fix this $x_0$, and set $A := \{ h \in B_1: x_0 + h \in H'' \}$, so that we have
$$ |A| \geq 2^{-6} K^{-13}|B_1|.$$
Observe that if $h,h' \in A$, then
$(h-h', \xi_{x_0+h} - \xi_{x_0+h'})$ lies in $\Gamma'' - \Gamma''$, which is a subgraph of $2\Gamma''-2\Gamma''$.  Thus we have
$\xi_{x_0+h} - \xi_{x_0+h'} = 2M(h-h')$.  Combining this with \eqref{madd}, we conclude that there exists $\xi_0 \in \widehat G$ such that
$$ \xi_{x_0+h} = \xi_0 + 2Mh \hbox{ for all } h \in A.$$ This concludes the proof of Proposition \ref{prop9.1a}.\endproof\vs

\ni Combining Proposition \ref{prop9.1a} with Proposition \ref{add-quad-2} leads immediately to the following, which generalizes Proposition \ref{add-quad-3} to arbitrary $G$. 

\begin{proposition}[Large $U^3(G)$-norm implies locally linear phase derivative]\label{local-linear-prop}
Let $G$ be an arbitrary finite additive group, and let $f: G \to \D$ be a bounded function such that $\|f\|_{U^3(G)} \geq \eta$ for some $\eta > 0$.  Then there exists a set $S \subseteq \widehat G$ with 
\[ d_1 := |S| \leq 2^{C_3}\eta^{-C'_3},\]
a regular Bohr set $B_1 :=B(S,\rho) \subseteq B(S,\frac{1}{4}) = B_0$ with
$\rho \in [\frac{1}{16},\frac{1}{8}]$, elements $x_0 \in G$ and $\xi_0 \in \widehat G$, and a function $M: B_0 \to \widehat G$ obeying the local linearity
property \eqref{madd}, such that
\begin{equation}\label{gabor}
\E_{h \in B_1} |\E_{x \in G} T^{x_0+h} f(x)  \overline{f(x)} e(-(\xi_0 + 2Mh) \cdot x)| \geq 2^{-C_4} \eta^{C'_4}.
\end{equation} We could take $C_i,C'_i = 2^{18}$, $i= 3,4$.
\end{proposition}

\ni\textit{Step 2: The symmetry argument.}
Let $S, B_1, x_0, \xi_0, M$ be as in Proposition \ref{local-linear-prop}.  Using the proof of Theorem \ref{mainthm1} as a model, the next step
would be to establish some symmetry property on $M: B_0 \to \widehat G$, in the sense that
the form $\{x,y\} := M(x) \cdot y - M(y) \cdot x$ is small.  More precisely, we shall establish

\begin{lemma}[Symmetry of derivative]  Let the notation be as in Proposition \ref{local-linear-prop}. 
For any $x, y \in B_0$, let $\{x,y\}$ denote the anti-symmetric form
$$ \{x,y\} := M(x) \cdot y - M(y) \cdot x.$$
Then there exists a set $S_3$ of frequencies with $S_3 \supseteq \frac{1}{2} \cdot S$ and
$d_3 := |S_3| \leq 2^{C_5} \eta^{-C'_5}$, and a Bohr set $B_3 = B(S_3, 2^{-C_6}\eta^{C'_6}) \subseteq  B_1$, such that
\begin{equation}\label{sym}
 \| \{ x, z \} \|_{\R/\Z} \leq 2^{C_7} \eta^{-C'_7} \| x \|_{S_3} \hbox{ for all } x, z \in B_3.
\end{equation}
It is permissible to take all of the $C_i, C'_i$, $i = 5,6,7$, equal to $2^{23}$.
\end{lemma}

\begin{proof} Let $\eps_2 = 2^{-C_3 - 2C_4 - 10} \eta^{C'_3 + 2C'_4}$.  By Lemma \ref{make-regular}, we can find $\rho_2 \in [\eps_2,2\eps_2]$ such that $B_2 := B(S, \rho_2) \subseteq  B_1$ is a
regular Bohr set. Of course we have
\[ d_2 = d_1 \leq 2^{C_3} \eta^{-C_3}.\] We write \eqref{gabor} as
\begin{equation}\label{eq9.71}|\E_{h \in B_1; x \in G} \Bounded(x+h) \Bounded(x) \Bounded(h) e(-2Mh \cdot x)| \geq 2^{-C_4} \eta^{C'_4},\end{equation}
absorbing all the phase terms into the functions $\Bounded$. Applying Lemma \ref{avg-1}
we can find $x_1 \in G$ such that
\begin{equation}\label{eq9.72}|\E_{h \in B_1; x \in B_2} \Bounded(x+x_1+h) \Bounded(x+x_1) \Bounded(h) e(-2Mh \cdot (x+x_1))| \geq 2^{-C_4} \eta^{C'_4}.\end{equation}
Absorbing the $x_1$ terms into the functions $\Bounded$ we conclude that
\begin{equation}\label{eq9.73}|\E_{h \in B_1; x \in B_2} \Bounded(h) \Bounded(x+h) \Bounded(x) e(-2Mh \cdot x)| \geq 2^{-C_4} \eta^{C'_4}.\end{equation}
Applying the Cauchy-Schwarz inequality (Lemma \ref{cz}) to eliminate $\Bounded(h)$, we then deduce
\begin{equation}\label{eq9.74}\E_{h \in B_1; x,y \in B_2} \Bounded(x+h) \Bounded(x) \Bounded(y+h) \Bounded(y) e(-2Mh \cdot (y-x)) \geq 2^{-2C_4} \eta^{2C'_4}.\end{equation}
Making the substitution $z := x+y+h$, this becomes
\begin{equation}\label{eq9.75}\E_{x,y \in B_2} \E_{z \in x+y+B_2} \Bounded(z,x) \Bounded(z,y) e(-2M(z-x-y) \cdot (y-x)) \geq 2^{-2C_4} \eta^{2C'_4}.\end{equation}
Absorbing as many phase terms into the functions $\Bounded(z,x)$ and $\Bounded(z,y)$ as we can, we conclude
$$\E_{x,y \in B_2} \E_{z \in x+y+B_1} \Bounded(z,x) \Bounded(z,y) e(2 \{x,y\} ) \geq 2^{-2C_4} \eta^{2C'_4}.$$
Next, by the regularity of $B_1$ and Lemma \ref{avg} (i), we observe that
\begin{eqnarray*} && \big |\E_{x,y \in B_2} \E_{z \in x+y+B_1} \Bounded(z,x) \Bounded(z,y) e(2 \{x,y\} ) \\ && \qquad \qquad\qquad\qquad\qquad-
\E_{x,y \in B_2} \E_{z \in B_1} \Bounded(z,x) \Bounded(z,y) e(2 \{x,y\} )\big| \leq 2^9 \eps_2 d\end{eqnarray*} which, due to the choice of $\eps_2$ and the bound $d \leq 2^{C_3} \eta^{-C'_3}$, implies that 
\begin{equation}\label{eq9.75a}\E_{x,y \in B_2} \E_{z \in B_1} \Bounded(z,x) \Bounded(z,y) e(2 \{x,y\} ) \geq 2^{-2C_4 - 1} \eta^{2C'_4}.\end{equation}
In particular, by the pigeonhole principle in $z$ we have
\begin{equation}\label{eq9.76}|\E_{x,y \in B_2} \Bounded(x) \Bounded(y) e(2 \{x,y\} )| \geq 2^{-2C_4 - 1} \eta^{2C'_4}\end{equation}
for some bounded functions $\Bounded(x), \Bounded(y)$.
At this point we observe the local bilinearity relationships
\begin{equation}\label{bilinear}
 \{ x+x',y\} = \{x,y\} + \{x',y\}; \quad \{x,y+y'\} = \{x,y\} + \{x,y'\} ,
\end{equation}
which hold whenever all four of $\|x\|_S, \|x'\|_S, \|y\|_S, \|y'\|_S$ are at most $\frac{1}{8}$. 
We can then apply Cauchy-Schwarz (Lemma \ref{cz}) to eliminate $\Bounded(x)$ and conclude that
$$|\E_{x,y,y' \in B_2} \Bounded(y,y') e(2 \{x,y'-y\} )| \geq 2^{-4C_4 - 2} \eta^{4C'_4},$$
and hence by the triangle inequality
$$ \E_{y,y' \in B_2} |\E_{x \in B_2} e(2 \{x,y'-y\})|  \geq 2^{-4C_4 - 2} \eta^{4C'_4}.$$
By the pigeonhole principle, there exists $y' \in B_2$ such that
\begin{equation}\label{eq9.79} \E_{y \in B_2} |\E_{x \in B_2} e(2 \{x,y'-y\})| \geq 2^{-4C_4 - 2} \eta^{4C'_4}.\end{equation} Fix this $y'$.
Since $|\E_{x \in B_2} e(2 \{x,y'-y\})|$ is bounded above by 1, we conclude that there exists a set $A \subseteq B_2$
with $|A| \geq 2^{-4C_4 - 3} \eta^{4C'_4}|B_2|$ such that
$$
|\E_{x \in B_2} e(2 \{x,y'-y\}) | \geq 2^{-4C_4 - 3} \eta^{4C'_4}\hbox{ for all } y \in A.
$$
Applying Lemma \ref{gfd} (and recalling that $d_2 \leq 2^{C_3} \eta^{-C'_3}$ and $\rho_2 \geq \eps_2 = 2^{-C_3 - 2C_4 - 10} \eta^{C'_3 + 2C'_4}$) we conclude that
$$ \| 2\{ x, y'-y \} \|_{\R/\Z} \leq  2^{24 + 2C_3 + 8C'_4}\eta^{-2C'_3 - 10 C'_4}\| x \|_S \hbox{ for all } x \in B_2, y \in A.$$
Applying \eqref{bilinear} (and recalling that $\rho_2 \leq 2\eps \leq \frac{1}{32}$), we conclude that
$$ \| 2\{ x, z \} \|_{\R/\Z} \leq  2^{26 + 2C_3 + 8C'_4}\eta^{-2C'_3 - 10 C'_4}\| x \|_S \hbox{ for all } x \in B_2, z \in 2A-2A.$$
On the other hand, by applying Lemma \ref{bog-lemma-local}, we can find $S' \subseteq  \widehat G$ with $|S \cup S'| \leq 2^{12C_4 + 16} \eta^{-12 C'_4}$ and
a Bohr set $B'_3 := B(S \cup S', 2^{-61 - 2C_3 - 26C_4}\eta^{2C'_3 + 26 C'_4})$ which is completely contained in $2A-2A$ and inside $B_2$.  Thus we have
$$ \| 2\{ x, z \} \|_{\R/\Z} \leq 2^{26 + 2C_3 + 8C'_4}\eta^{-2C'_3 - 10 C'_4}\| x \|_S \hbox{ for all } x, z \in B'_3.$$
Let us now eliminate the factor $2$.  Observe that $B_3 := \{ 2x: x \in B'_3 \}$ is also a Bohr set (with the frequency set $S \cup S'$ replaced
by $S_3 := \frac{1}{2} \cdot (S \cup S')$).  Since $\|x\|_S \leq 2 \|x\|_{S_3}$, we also observe that $B_3 \subseteq  B(S,2\rho_2) \subseteq  B_1$.
By \eqref{bilinear}, which implies that $\{ 2x, z\} = 2\{x,z\}$, we conclude \eqref{sym} as desired.
\end{proof}\vs

\ni There are extremely close analogies between the above argument and that of \S \ref{finite-sec}. Equations \eqref{eq9.71}, \eqref{eq9.72}, \eqref{eq9.73}, \eqref{eq9.74}, \eqref{eq9.75}, \eqref{eq9.75a}, \eqref{eq9.76} and \eqref{eq9.79} are analogous to \eqref{tff-1}, \eqref{tff-2}, \eqref{tff-3}, \eqref{tff-4}, \eqref{tff-5}, \eqref{tff-5a}, \eqref{tff-6} and \eqref{tff-9} respectively.\vs

\ni\textit{Step 3: Eliminating the quadratic phase component.} We now return to the conclusion of Proposition \ref{local-linear-prop} and localize the $x$ and $h$ variables to a small Bohr set.
Let \[ \eps_4 = \min \big(2^{-C_3 - C_4 - 10}\eta^{C'_3 + C'_4}, 2^{-5 - C_4 - C_7} \eta^{C'_4} \big),\] let $S_3$ be the set of characters coming from the previous subsection, and let
$B_4 := B(S_3, \rho_4) \subseteq  B_3$ be a regular Bohr set such that $\rho_4 \in [\eps_4,2\eps_4]$.  By the previous estimate we have
\begin{equation}\label{symm-2}
 \| \{ x, z \} \|_{\R/\Z} \leq 2^{C_7} \eta^{-C'_7}\eps_4 \hbox{ for all } x, z \in B_4.
\end{equation}
Let us write \eqref{gabor} as
$$
|\E_{h \in B_1; x \in G} \Bounded(h) \Bounded(x+h)  \overline{f(x)} e(- 2Mh \cdot x)| \geq 2^{-C_4} \eta^{C'_4},$$
where we have absorbed some phase terms into the functions $\Bounded$ as before.
Since $B_1 = B(S,\rho)$ for some $\rho \geq \frac{1}{16}$, we conclude from Lemma \ref{avg} (ii) that
\begin{eqnarray*}
& & \big|\E_{h \in B_1; x \in G} \Bounded(h) \Bounded(x+h)  \overline{f(x)} e(- 2Mh \cdot x) - \\
 & & \qquad \E_{h' \in B_1; h \in B_4; x \in G} \Bounded(h+h') \Bounded(x+h+h')  \overline{f(x)} e(- 2M(h+h') \cdot x)\big| \leq 2^9 \eps_4 d_1.
\end{eqnarray*} This is at most $2^{-C_4 - 1} \eta^{C'_4}$, and therefore
$$ |\E_{h' \in B_1; h \in B_4; x \in G} \Bounded(h+h') \Bounded(x+h+h')  \overline{f(x)} e(- 2M(h+h') \cdot x)|
\geq 2^{-C_4 - 1} \eta^{C'_4}.$$
Hence by the pigeonhole principle, there exists $h' \in B_r$ such that
$$ |\E_{h \in B_4; x \in G} \Bounded(h+h') \Bounded(x+h+h')  \overline{f(x)} e(- 2M(h+h') \cdot x)|
\geq 2^{-C_4 - 1} \eta^{C'_4}.$$
Since $e(-Mh' \cdot x) = e(-2Mh' \cdot (x + h)) e(2Mh' \cdot h)$, we can absorb all the $h'$ terms into the functions $\Bounded$ to conclude that
$$ |\E_{h \in B_4; x \in G} \Bounded(h) \Bounded(x+h)  \overline{f(x)} e(- 2Mh \cdot x)|
\geq 2^{-C_4 - 1} \eta^{C'_4}$$
By Lemma \ref{avg-1} we then have
$$ |\E_{y \in G; x, h \in B_4} \Bounded(h) \Bounded(x+y+h)  \overline{f(x+y)} e(- 2Mh \cdot (x+y))|
\geq 2^{-C_4 - 1} \eta^{C'_4},$$
and hence by the triangle inequality
$$ \E_{y \in G} |\E_{x, h \in B_4} \Bounded(h,y) \Bounded(x+h,y)  \overline{f(x+y)} e(- 2Mh \cdot x)|
\geq 2^{-C_4 - 1} \eta^{C'_4}.$$
Now we observe from \eqref{madd} that
$$ 2Mh \cdot x = M(x+h) \cdot (x+h) - Mx \cdot x - Mh \cdot h - \{x,h\},$$
and hence by \eqref{symm-2}
$$ \big|e(-2Mh \cdot x) - \Bounded(x+h) \Bounded(h) e(Mx \cdot x)\big| \leq 2\pi \cdot 2^{C_7} \eta^{-C'_7} \eps_4 \leq 2^{-C_4 - 2} \eta^{C'_4}.$$
Thus we have
$$ \E_{y \in G} \big|\E_{x, h \in B_4} \Bounded(h,y) \Bounded(x+h,y)  \overline{f(x+y)} e(Mx \cdot x)\big|
\geq 2^{-C_4 - 2}\eta^{C'_4} .$$
It is convenient to localize $x$ further.
Let $\eps_5 = 2^{-C_4 - C_5 - 13} \eta^{C'_4 + C'_5} \eps_4$.  By Lemma \ref{make-regular} we can find a regular Bohr set $B_5 = B(S_3, \rho_5)$ with $\rho_5 \in [\eps_5,2\eps_5]$.  By Lemma \ref{avg} (ii) we have that
\begin{eqnarray*}
& & \bigg| \E_{y \in G} \big|\E_{x, h \in B_4} \Bounded(h,y) \Bounded(x+h,y)  \overline{f(x+y)} e(Mx \cdot x)\big|\\
& & \quad - \E_{y \in G} |\E_{w, h \in B_4; x \in B_5} \Bounded(h,y) \Bounded(x+w+h,y)  \overline{f(x+w+y)} e(M(x+w) \cdot (x+w))\big|\bigg|   
\end{eqnarray*}
is at most $2^{10}\eps_5 d_3/\eps_4$,
which on account of the choice of $\eps_5$ implies that
$$ \E_{y \in G} |\E_{w, h \in B_4; x \in B_5} \Bounded(h,y) \Bounded(x+w+h,y)  \overline{f(x+w+y)} e(M(x+w) \cdot (x+w))|
\geq 2^{-C_4 - 3}\eta^{C'_4}.$$
By the pigeonhole principle there exists $w \in B_4$ such that
$$ \E_{y \in G} |\E_{h \in B_4; x \in B_5} \Bounded(h,y) \Bounded(x+w+h,y)  \overline{f(x+w+y)} e(M(x+w) \cdot (x+w))|
\geq 2^{-C_4 - 3}\eta^{C'_4}.$$
Let us now apply Lemma \ref{fgh-bohr}. Since $B_4$ is regular and $\eps_5$ is so small, we certainly have $\E (1_{B_4})/(\E 1_{B_4 + B_5}) \geq 1/2$, and so that lemma allows us to conclude that 
$$ \E_{y \in G} \| \overline{f(x+w+y)} e(M(x+w) \cdot (x+w)) \|_{u^2(B_5)}
\geq 2^{-C_4 - 4}\eta^{C'_4}.$$ This, of course, implies that 
\[ \E_{y \in G} \| \overline{f(x+w+y)} e(M(x+w) \cdot (x+w)) \|_{u^3(B_5)}
\geq 2^{-C_4 - 4}\eta^{C'_4}.\]
The $u^3$ norm being invariant under translation, conjugation and quadratic phase modulation, we conclude that 
$$ \E_{y \in G} \| f \|_{u^3(y+w +B_5)}
\geq 2^{-C_4 - 4}\eta^{C'_4}.$$
Making the change of variables $y \mapsto y+w$, and completing a small computation, we obtain \eqref{efphi-general} as desired.
\endproof\vs

\ni We remark that we have proved slightly more than \eqref{efphi-general}, in that the quadratic phase functions used to demonstrate
the largeness of the $\|f\|_{u^3(y+B)}$ norm all agree up to lower order (i.e. linear and constant) terms.  However we were unable to find any way
to exploit this additional fact.

\section{Bohr sets and generalized arithmetic progressions}\label{sec10}

\ni Our focus from this point on is largely on the group $G = \Z/N\Z$, as we are working towards connections with themes in ergodic theory, in particular involving $\Z$-actions. As we have stressed, $\Z/N\Z$ is an appropriate group to consider if one is interested in discrete questions concerning the integers. However some of what we have to say, particularly in the present section, can be generalized without undue pain to arbitrary additive groups. \vs

\ni We have obtained an inverse theorem, Theorem \ref{mainthm2}, for the $U^3$ norm which relates that norm to the quadratic bias norm
$u^3(y+B)$ on Bohr sets $B$ (or on subspaces $W$, in finite field cases such as $G=\Ffiven$).  This is a fairly satisfactory state of affairs, except for
the presence of the Bohr set $B$; in particular, it is not clear at present what exactly the locally quadratic phase functions are on
$B$.  In this section we show how the Bohr set can, if desired, be replaced with a \textit{generalized arithmetic progression,} and how to characterize the
locally quadratic phase functions on such progressions.\vs

\ni We begin by recalling what a generalized arithmetic progression is.

\begin{definition}[Generalized arithmetic progression]  A \emph{generalized arithmetic progression} $P$ in an additive group $G$ is any set 
of the form
$$ P := \{ a + l_1 v_1 + \ldots + l_d v_d: 0 \leq l_j < L_j \hbox{ for all } 1 \leq j \leq d \}$$
where $d \geq 0$, $a,v_1,\ldots,v_d \in G$, and $L_1,\ldots,L_d \geq 1$.  We shall abbreviate the right-hand side as
$P = a + [0,L) \cdot v$, where $L := (L_1,\ldots,L_d)$ and $v := (v_1,\ldots,v_d)$.
We call $a$ the \emph{base point} of the progression, $d$ the \emph{rank},
$v_1,\ldots,v_d$ the \emph{generators}, and $L_1,\ldots,L_d$ the \emph{lengths} of the progression.  If all the sums in $P$ are distinct, so that 
thus $|P| = L_1 \ldots L_d$, we say that $P$ is \emph{proper}.  A \emph{coset progression} is any set of the form $P+H$ where $H$ is a subgroup of $G$.  We say that the coset progression $P+H$ is \emph{proper} if $P$ is proper and $|P+H|=|P| |H|$ (i.e. all the sums in $P+H$ are distinct); we define the rank, base point, etc. of the coset progression $P+H$ to be the same as that of its component $P$.
\end{definition}

\ni The need to generalize from generalized arithmetic progressions to coset progressions in the setting of a general group $G$ was first noted 
in \cite{green-ruzsa-freiman}. \footnote{Of course, the classification of finite abelian groups tells us that every coset progression is also a generalized
arithmetic progression (by expanding $H$ as the direct sum of cyclic groups, which can each be interpreted as an arithmetic progression) but in doing so one can cause the rank of the coset progression to increase enormously (by the number of generators needed to span $H$).}\vs

\ni We now use standard facts from the geometry of numbers to
show that every Bohr set contains a large proper coset progression.  The first lemma follows from a result of Mahler (\cite[Chapter VIII, Corollary to Theorem VII]{cassels}) together with Minkowski's Second Theorem (loc. cit, Chapter VIII, Theorem V). This result (in fact, a rather stronger one) was used in an additive-combinatorial context in Bilu's work on Freiman's theorem \cite[Lemma 2.1]{bilu}.

\begin{lemma}\label{wbasis}
Let $\Gamma$ be a lattice of full rank in $\R^d$. Then there exists linearly independent vectors $w_1, \ldots, w_d \in \Gamma$
which generate $\Gamma$, and such that
\begin{equation}\label{w-magic}
|w_1| \ldots |w_d| \leq 2 \cdot d! \cdot  {\rm mes}(\R^d/\Gamma),
\end{equation}
where ${\rm mes}(\R^d/\Gamma)$ is the volume of a fundamental domain of $\Gamma$.
\end{lemma}

\ni Next, we give a ``discrete John's theorem'' which shows that the intersection of a convex symmetric body and a lattice of full rank is essentially equivalent to a progression.

\begin{lemma}[Discrete John's theorem]\label{djt}  Let $B$ be a convex symmetric body in $\R^d$, and let $\Gamma$ be a lattice in $\R^d$ of full rank.
Then there exists a $d$-tuple $$w = (w_1,\ldots,w_d) \in
\Gamma^d$$ of linearly independent vectors in $\Gamma$ and and a
$d$-tuple $L = (L_1,\ldots,L_d)$ of positive integers such that
$$ (d^{-2d} \cdot B) \cap  \Gamma \subseteq  (-L,L) \cdot w \subseteq  B \cap \Gamma \subseteq  ( - d^{2d} L, d^{2d} L ) \cdot w.$$  Here of course
$$ (-L,L) \cdot w := \{ l_1 w_1 + \ldots + l_d w_d: -L_j < l_j < L_j \hbox{ for all } 1 \leq j \leq d \}.$$
\end{lemma}

\begin{proof}
We first observe using John's theorem \cite{john} (see also \cite{bollobas,pisier}) and an invertible linear transformation that we may assume
without loss of generality that $B_d \subseteq  B \subseteq  d \cdot B_d$, where $B_d$ is the unit ball in $\R^d$.
We may also assume $d \geq 2$, since the claim is easy otherwise.\vs

\ni Now let $w = (w_1, \ldots, w_d)$ be as in Lemma \ref{wbasis}.  For each $j$, let $L_j$ be the least integer greater than $1/d|w_j|$.
Then from the triangle inequality we see that $|l_1 w_1 + \ldots + l_d w_d| < 1$ whenever $|l_j| < L_j$, and hence $(-L,L) \cdot w$ is contained
in $B_d$ and hence in $B$.\vs

\ni Now let $x \in B \cap \Gamma$.  Since $w$ generates $\Gamma$, we have $x = l_1 w_1 + \ldots + l_d w_d$ for some integers $l_1,\ldots,l_d$;
since $B \subseteq  d \cdot B_d$, we have $|x| \leq d$.  Applying Cramer's rule to solve for $l_1,\ldots,l_d$ and \eqref{w-magic}, we have
\[ |l_j| = \frac{|x \wedge w_1 \ldots w_{j-1} \wedge w_{j+1} \wedge w_d|}{|w_1 \wedge \ldots \wedge w_d|}
 \leq \frac{|x| |w_1| \ldots |w_d|}{|w_j| |w_1 \wedge \ldots \wedge w_d|} = \frac{|x|{\rm mes}(\R^d/\Gamma)}{|w_j|}
 \leq \frac{2d \cdot d!}{|w_j|}, 
\]
which is certainly at most $d^{2d} L_j$. It follows that $x \in ( - d^{2d} L, d^{2d} L ) \cdot w$, which is what we wanted to prove.  A more-or-less identical argument gives the inclusion
$(d^{-2d} \cdot B) \cap  \Gamma \subseteq  (-L,L) \cdot w$.
\end{proof}\vs

\ni Let $x \mapsto \{x\}$ denote the fractional part map from $\R/\Z$ to the fundamental domain $(-1/2,1/2]$. 

\begin{lemma}[Bohr sets contain large coset progressions]\label{bohr-cp}  Let $S \subseteq \widehat{G}$ be a set of $d$ characters, let $\rho < 1/4$ be a real number, and let $B(S,\rho) \subseteq  G$ be a Bohr set.  Then there exists a proper coset progression $P+H$ of rank $d'$, $0 \leq d' \leq d$, where 
$P = (-L,L) \cdot v$ for some $L_1,\ldots,L_{d'} \geq 1$ and $v_1,\ldots,v_{d'} \in G$, and we have the inclusions
\begin{equation}\label{eq10.1} B(S, d^{\prime -2d'} \rho) \subseteq  P+H \subseteq  B(S,\rho).\end{equation}
In particular, from Lemma \ref{bohr-bound} we have
\begin{equation}\label{eq10.2} |P+H| \geq \rho^d d^{-2d^2}N.\end{equation}
Furthermore, the vectors $(\{ \xi \cdot v_j \})_{\xi \in S} \in \R^S$, $1 \leq j \leq d'$,
can be chosen to be linearly independent, and $H$ can be taken to be the orthogonal complement of $S$, that is to say the group
\begin{equation}\label{H-def}
H := \{ x \in G: \xi \cdot x = 0 \hbox{ for all } S \}.
\end{equation}
\end{lemma}
\noindent\textit{Remark.} The lemma is at the same time a refinement and a weakening of a lemma from \cite{green-ruzsa-freiman}. The refinement, corresponding to the fact that we use Lemma \ref{wbasis} rather than Minkowski's second theorem, is that we obtain the left-hand inclusion in \eqref{eq10.1} and not just the right-hand one. The weakening is that using just Minkowski's second theorem (and thus sacrificing the left-hand inclusion in \eqref{eq10.1}) gives a stronger bound than \eqref{eq10.2}.\vs

\begin{proof}  Let $\phi: G \to (\R/\Z)^S$ be the group homomorphism $\phi(x) := (\xi \cdot x)_{\xi \in S}$.  Observe that $\phi(G)$ is a finite subgroup of the torus $(\R/\Z)^S$, and that $B(S,\rho)$ is the inverse image of the cube $Q := \{ (y_\xi)_{\xi \in S}: |y_\xi| \leq \rho \}$
under $\phi$.  \vs

\ni Let $\Gamma \subseteq \R^S$ be the lattice $\phi(G) + \Z^S$. Though it is a slight abuse of notation, we consider $\phi(G) \cap Q$ to be the same as
$\Gamma \cap Q$. Applying Lemma \ref{djt}, we can find a progression $\tilde P := (-L,L) \cdot w$ for some linearly independent
$w_1,\ldots,w_{d'} \subseteq  \Gamma$ with $0 \leq d' \leq d$
such that 
$$\Gamma \cap d^{\prime -2d'} \cdot Q \subseteq  \tilde P \subseteq  \Gamma \cap Q.$$
Since the $w_j$ are independent, $\tilde P$ is necessarily proper.  The claim now follows by setting $v_j$ to be an arbitrary element of
$\phi^{-1}(w_j)$ for each $1 \leq j \leq d'$, and setting $H$ equal to the kernel of $\phi$, which is of course just \eqref{H-def}.
\end{proof}\vs 

\ni When $G = \Z/N\Z$ is a cyclic group of prime order, the subgroup $H$ has no r\^{o}le to play\footnote{At the other extreme, in the finite field geometry setting
$G = \Ffiven$ it is the progression component $P$ which is irrelevant, because the properness of $P$ forces all the lengths to be less
than 5.  Since the rank of $P$ is also under control, we thus see that $H$ is a substantial portion of $P+H$.  Indeed,
the fact that Bohr sets in finite field geometries contain large subspaces was already exploited in the proof of Theorem \ref{mainthm1}.} and we conclude the following corollary.

\begin{corollary}\label{prog-large}  Let $G = \Z/N\Z$ be a cyclic group of prime order, let $S \subseteq \widehat{G}$ be a set of $d$ characters, and let $\rho < 1/4$ be a parameter. Then there is a proper generalized arithmetic progression $P = (-L,L) \cdot v$ of rank at most $d$ and size at least $\rho^d d^{-2d^2} N$ such that $B(S, d^{-2d} \rho) \subseteq  P \subseteq  B(S,\rho)$.  Furthermore, the vectors
$(\{ \xi \cdot v_j \})_{\xi \in S}$ are linearly independent in $\R^S$.
\end{corollary}

\ni In this paper, it is in general more convenient technically to work with Bohr sets than progressions or coset progressions.
However, there is one task which is much easier to achieve on progressions than on Bohr sets, and that is to classify quadratic phase functions:

\begin{lemma}[Inverse theorem for locally quadratic functions]\label{inv-quad}  Let $G$ be a finite additive group of odd order.
Let $P+H$ be a coset progression in $G$, let $a$ be the base point of $P$, and let $v_1,\ldots,v_d$ be the generators.  Let 
$\phi: P+H \to \R/\Z$ be a locally quadratic phase function on $P+H$.  Then there exists a self-adjoint homomorphism $M: H \mapsto \widehat H$,
elements $\xi_0, \xi_1, \ldots, \xi_d \in \widehat H$, elements $\eta_i, \lambda_{ij} \in \R/\Z$ for $1 \leq i,j \leq d$, and $c \in \R/\Z$
such that $\lambda_{ij} = \lambda_{ji}$ and
\begin{equation}\label{phinv}
\begin{split}
 \phi(a + l_1 v_1 + \ldots + l_d v_d + h) &= Mh \cdot h + 
 2\sum_{i=1}^d l_i \xi_i \cdot h +  \sum_{1 \leq i,j \leq d} l_i l_j \lambda_{ij} \\ 
&+ \xi_0 \cdot h + \sum_{i=1}^d l_i \eta_i + c
\end{split}
\end{equation}
for all $l_1,\ldots,l_d,h$ with $0 \leq l_j < L_j$ for all $1 \leq j \leq d$ and $h \in H$.
\end{lemma}

\ni\textit{Remark.} In the converse direction, it is easy to show that \eqref{phinv} is indeed well-defined and gives a locally quadratic function if
$3P + H$ is a proper coset progression, but we will not need that fact here.\vs

\begin{proof}  We may assume that $L_j \geq 2$ for all $j$, and we may translate so that $a=0$.
Let $\phi|_H$ be the restriction of $\phi$ to $H$.  By Theorem \ref{QHB}, the quadratic extension theorem, we can extend $\phi|_H$ to a globally quadratic phase function $\psi$ on $G$.  Using Lemma \ref{quadratic-classify}, it is easy to see that $\psi$, when restricted to $P+H$, has
the form \eqref{phinv}.  Thus we may subtract off $\psi$ from $\phi$, which means that $\phi$ now vanishes on $H$.  We now claim that under this reduction, $\phi$ takes the simpler form
$$
 \phi(l_1 v_1 + \ldots + l_d v_d + h) = 
 \sum_{i=1}^d l_i \xi_i \cdot h +  \sum_{1 \leq i,j \leq d} l_i l_j \lambda_{ij} + \sum_{i=1}^d l_i \eta_i.$$
 Observe that for any $h \in H$, the function $(h \cdot \nabla) \phi$ is locally linear on $P+H$ and vanishes on $H$, and hence takes the form
$$
 (h \cdot \nabla) \phi(l_1 v_1 + \ldots + l_d v_d + h') =  \sum_{i=1}^d l_i (h \cdot \nabla) \phi(v_i) $$
for all $l_1,\ldots,l_d,k$ with $0 \leq l_j < L_j$ and $h' \in H$.  
It is then easy to see that $h \mapsto (h \cdot \nabla) \phi(v_i)$ is a group homomorphism from $H$ to $\R/\Z$
and thus there exists $\xi_i \in \widehat H$ such that $(h \cdot \nabla) \phi(v_i) = 2\xi_i \cdot h$ for all $h \in H$.  Using this, we thus reduce
to showing that
$$
 \phi(l_1 v_1 + \ldots + l_d v_d) = \sum_{1 \leq i,j \leq d} l_i l_j \lambda_{ij} + \sum_{i=1}^d l_i \eta_i.$$
or equivalently that
$$
 \phi(l_1 v_1 + \ldots + l_d v_d) = \sum_{1 \leq i < j \leq d} 2 l_i l_j \lambda_{ij} + \sum_{i=1}^d l_i \eta_i + l_i^2 \lambda_i.$$
 We induct on $d$.  When $d=0$ there is nothing to prove.  Now suppose that the claim is already proven for $d-1$.  We observe that the derivative $(v_d \cdot \nabla) \phi$ is linear, and hence
$$
(v_d \cdot \nabla) \phi(l_1 v_1 + \ldots + l_d v_d) = \sum_{1 \leq i \leq d} 2 l_i \lambda_{id} + \eta_d$$
for some $\lambda_{1d},\ldots,\lambda_{dd},\eta_d \in \R/\Z$, with the caveat that $l_d$ now must be less than $L_d-1$ rather than $L_d$.
(Note that it is always possible to divide by two in $\R/\Z$, though the value obtained need not be unique).
The claim then follows from the induction hypothesis and a simple ``integration'' argument which we omit.
\end{proof}\vs

\ni Now, let us specialize to the setting of cyclic groups $\Z/N\Z$ of prime order.  We begin by defining some special functions on this set.

\begin{definition}[Bracket polynomials]  Let $\Z/N\Z$ be a cyclic group of prime order.  If $k \geq 0$, we define a \emph{bracket monomial of degree $k$} on $\Z/N\Z$ to be any function $\phi: \Z/N\Z \to \R/\Z$ of the form
$$ \phi(x) = a \{ \xi_1 \cdot x \} \ldots \{ \xi_k \cdot x\} \mod 1$$
where $\xi_1,\ldots, \xi_k \in \widehat{\Z/N\Z}$ and $a \in \R$; we refer to $\xi_1,\ldots,\xi_k$ as the \emph{frequencies} of the monomial.  If $L \geq 0$ and $S \subseteq  \widehat{\Z/N\Z}$, we define a \emph{bracket polynomial of degree at most $k$, length at most $L$ and frequency set $S$} to be any function $\phi: \Z/N\Z \to \R/\Z$ which can be expressed as the sum of $L$ or fewer bracket monomials of degree at most $k$ and frequencies inside $S$. We write $\Freq(\phi) \subseteq S$. Note that if $|S| = d$ then we may always assume that $L \leq kd^k$; for this reason there will be little subsequent discussion of length.
\end{definition}

\ni These bracket polynomials are special cases of the \textit{generalized polynomials} considered in various papers of H{\aa}land, H{\aa}land-Knuth, Bergelson and Leibman. \cite{berg-hal,bergelson-leibman2,hal1,hal2,hal3}, though with the (minor) caveat that our fractional parts take values from $-\half$ to $\half$, while the ones in those papers take values from $0$ to $1$.  \vs

\ni Define a \emph{bracket quadratic} to be a bracket polynomial of degree at most 2.
We can now link locally quadratic phase functions with bracket quadratics.

\begin{proposition}\label{bohr-agree}  Let $G = \Z/N\Z$ be a cyclic group of prime order, let $S \subseteq \widehat G$ be a set of $d$ characters, and suppose that $\rho \in (0,\quarter]$. Let $P$ be the proper progression contained in $B(S,\rho)$ which was constructed in Corollary \ref{prog-large}, and let $\phi: B(S,\rho) \to \R/\Z$ be a locally quadratic phase function on $B(S,\rho)$.  Then there exists a bracket quadratic $\tilde \phi: \Z/N\Z \to \R/\Z$ with $\mbox{\emph{Freq}}(\widetilde{\phi}) \subseteq S$ such that $\widetilde{\phi} = \phi$ on $P$.  
\end{proposition}

\begin{proof}  Let $v_1,\ldots,v_d$ and $L_1,\ldots,L_d$ be the generators and lengths of $P$.  By Lemma \ref{inv-quad} we have
\begin{equation}\label{phin-expand}
\phi(l_1 v_1 + \ldots + l_d v_d ) = \sum_{1 \leq i \leq j \leq d} l_i l_j \lambda_{ij} \\ 
+ \sum_{i=1}^d l_i \eta_i + c \mod 1
\end{equation}
for some real numbers $\lambda_{ij}$, $\eta_i$. \vs 

\ni Next, let $\Phi: P \to \R^S$ be the map
$$ \Phi(x) := \{ ( \{ \xi \cdot x \} )_{\xi \in S} ).$$   Since $P \subseteq  B(S, \rho)$, we see that
$\Phi(P)$ lies inside the cube $Q := \{ (y_\xi)_{\xi \in S}: |y_\xi| \leq \rho \hbox{ for all } \xi \in S \}$.  Since $\rho < \frac{1}{4}$, it is also
easy to verify that 
$$ \Phi(l_1 v_1 + \ldots + l_d v_d) = l_1 \Phi(v_1) + \ldots + l_d \Phi(v_d).$$
From Corollary \ref{prog-large} we know that the $\Phi(v_j)$ are linearly independent.  Thus there exists a vector $u_i \in R^S$ such that
$$ \Phi(l_1 v_1 + \ldots + l_d v_d) \cdot u_i = l_i.$$
Writing $u_i = (u_{i\xi})_{\xi \in S}$ and $x = l_1 v_1 + \ldots + l_d v_d$ we conclude that
$$ l_i = \sum_{\xi \in S} u_{i\xi} \{ \xi \cdot x \}.$$
Inserting this formula into \eqref{phin-expand} we obtain the claim.
\end{proof}\vs

\ni We can now give a version of Theorem \ref{mainthm2}, the $U^3(G)$ inverse theorem, in the case $G = \Z/N\Z$, which involves bracket quadratic functions. In this theorem $C_0,C_1,\dots$ and $c$ denote absolute constants which \textit{do not} vary from line to line.

\begin{theorem}[Inverse theorem for $U^3(\Z/N\Z)$, bracket quadratic functions]\label{u3-quad} 
Let $\eta \in$

\ni $[0,c)$, and let $\Z/N\Z$ be a cyclic group of prime order.
If $f: \Z/N\Z \to \D$ is a bounded function such that $\|f\|_{U^3(\Z/N\Z)} \geq \eta$, then there exists a set $S \subseteq \widehat{\Z/N\Z}$ of size $d \leq \eta^{-C}$ and a proper progression $P$ of rank at most $d$ and size $|P| \geq \exp(-\eta^{-C})N$
with the inclusions
\begin{equation}\label{bohr-p}
B(S, \exp(-\eta^{-C_1})) \subseteq  P \subseteq  B(S, \textstyle\frac{1}{8}\displaystyle),
\end{equation}
and there exists a generalized quadratic $\phi: \Z/N\Z \to \R/\Z$ with $\mbox{\emph{Freq}}(\phi) \subseteq S$ such that we have the local quadratic bias estimate
\begin{equation}\label{thip}
 |\E_{x \in P}(T^h f(x) e(-\phi(x)) )| \geq \eta^{C_1}
\end{equation}
for some $h \in \Z/N\Z$.  More generally, for any non-empty set $A \subseteq  P$ there exists an $h_A \in \Z/N\Z$ such that
\begin{equation}\label{thip-better}
 |\E_{x \in A}(T^{h_A} f(x) e(-\phi(x)))| \geq \eta^{C_1}.
\end{equation}
Furthermore, there exists another generalized quadratic $\tilde \phi: \Z/N\Z \to \R/\Z$ with  $\mbox{\emph{Freq}}(\phi) \subseteq S$ such that we have the global quadratic bias estimate
\begin{equation}\label{thip-tilde}
 |\E_{x \in \Z/N\Z}(f(x) e(-\tilde \phi(x)))| \geq \exp(-\eta^{-C_2}).
\end{equation}
 
\ni Conversely, suppose that $S$ is a set of $d$ frequencies and that $\widetilde \phi : \Z/N\Z \rightarrow \R/\Z$ is a bracket quadratic with $\mbox{\emph{Freq}}(\phi) \subseteq S$. Suppose that $f : \Z/N\Z \rightarrow \D$ is a function such that $|\E (f(x) e(-\widetilde{\phi}(x))| \geq \eta$. Then we have
\begin{equation}\label{eq10.99} \Vert f \Vert_{U^3} \geq (\eta^3 \rho^2/C_3d^3)^d.\end{equation}
\end{theorem}

\begin{proof}
Applying Theorem \ref{mainthm2}, we can find a regular Bohr set $B := B(S,\rho)$ in $\Z/N\Z$ 
with $d = |S| \leq \eta^{-C_0}$ and $\rho \geq \eta^{C_0}$, a $y \in \Z/N\Z$,
and a locally quadratic phase function $\phi_0: B \to \Z/N\Z$ on $B$ such that
$$ |\E_{x \in B}(T^y f(x) e(-\phi_0(x)))| \geq \eta^{C_0}.$$ We could take $C_0 = 2^{25}$.
Next, let $\eps = \eta^{2C_0 + 11}$ and apply Corollary \ref{prog-large} to find a proper progression $P$ of rank at most $d$ such that
$$ B(S, \eps \exp(-\eta^{-C_0 - 1})\rho) \subseteq  P \subseteq  B(S, \eps \rho).$$
By Lemma \ref{avg} (iii), we can find $w \in B(S,(1 - \eps)\rho)$ such that 
\[ |\E_{x \in w + P} (T^yf(x) e(-\phi_0(x))) | \geq \eta^{C_0 + 1},\]
and \eqref{thip} follows after translating by $w$ and applying Proposition \ref{bohr-agree}.  A very similar argument 
gives \eqref{thip-better}.  Note that the inclusions on $P$
follow by choice of $\eta$, and the lower bound on $|P|$ follows from Lemma \ref{bohr-bound}.\vs

\ni Now we prove \eqref{thip-tilde}. From \eqref{bohr-p} and Lemma \ref{make-regular} we can find a regular Bohr set $B' := B(S,\rho')$ with
$\rho' \geq \exp(-\eta^{-C_1 - 3})$ which is contained in $P$.  Applying \eqref{thip-better} we can find a shift $h' \in \Z/N\Z$ such that
$$  |\E(T^{h'} f(x) e(-\phi(x)) 1_{B'}(x))| \geq \eta^{C_1 + 1}\E 1_{B'}.$$
Let $\eps' = \eta^{2C_1 + 9}$, and let $\chi: \R \to [0,1]$ be a smooth cutoff 
such that $\chi(s) = 1$ when $|s| \leq \rho'(1 - \eps')$,
$\chi(s) = 0$ when $|s| \geq \rho'(1 + \eps')$ and such that the derivative estimate $\|\chi''\|_{\infty} \leq 100/\eps^{\prime 2}\rho^{\prime 2}$ holds true.  By the regularity of $B'$ we see that
$$  |\E_x T^{h'} f(x) e(-\phi(x)) \prod_{\xi \in S} \chi( \{ \xi \cdot x \} )|
\geq \big( \eta^{C_1 + 1} - 200\eps' d \big) \E 1_{B'} \geq \eta^{C_1 + 2} \E 1_{B'},$$ and this is at most $\exp(-\eta^{-2C_1  - 4})$ by Lemma \ref{bohr-bound}.
Next, we use Fourier expansion on $\R$ to write
$$ \chi(s) = \int_\R \widehat \chi(t) e( t s )\ dt \qquad \mbox{where} \qquad \widehat \chi(t) := \int_\R \chi(s) e(-ts)\ ds.$$
This allows us to conclude that
\begin{equation}\label{eq10.10} \big|\int_\R \ldots \int_\R [\E_x T^{h'} f(x) e(-\phi(x) + \sum_{\xi \in S} t_\xi \{ \xi \cdot x \} )] \prod_{\xi \in S} \widehat \chi(t_\xi) dt_\xi \big| \geq \exp(-\eta^{-2C_1 - 4}).\end{equation}
Now set $\lambda := \eps^{\prime -1/2} \rho^{\prime - 1}$. Then by integration by parts, applied twice, we have
\[ \int^{\infty}_{\lambda} |\widehat \chi(t)| \, dt = \frac{1}{4\pi^2} \int^{\infty}_{\lambda} \frac{dt}{t^2} \big| \int_{\R} \chi''(s) e(-ts) \, ds \big| \leq \frac{100}{\eps' \rho'} \int^{\infty}_{\lambda} \frac{dt}{t^2} \leq \frac{100}{\eps^{1/2}}\]
and whence
\[ \int_{\R} |\widehat \chi (t)|\, dt \leq 2\lambda \| \widehat \chi \|_{\infty} + \frac{200}{\eps^{1/2}} \leq 8\lambda \rho' + \frac{200}{\eps^{1/2}} \leq 2^8 \eps^{-1/2}.\]
This implies that 
\[ \int_{\R} \dots \int_{\R} \prod_{\xi \in S} |\widehat \chi(t_{\xi}) |\, dt_{\xi} \leq \exp(-\eta^{-C_1 - 1}).\]
Comparing this with \eqref{eq10.10}, we conclude the existence of real numbers $t_\xi$ for $\xi \in S$ such that
$$ |\E_x T^{h'} f(x) e(-\phi(x) + \sum_{\xi \in S} t_\xi \{ \xi \cdot x \} )| \geq \exp(-\eta^{-3C_1 - 5})$$
and \eqref{thip-tilde} follows.\vs

\ni Finally, we prove \eqref{eq10.99}. Assume then that $\widetilde{\phi}$ is a bracket quadratic with $\Freq(\phi) \subseteq S$, $|S| = d$, and that $f : \Z/N\Z \rightarrow \D$ is a function with $|\E (f(x) e(-\widetilde{\phi}(x))| \geq \eta$. Set $\eps := \eta/160 d$, and select a $\rho \in [\eps,2\eps]$ such that $B := B(S,\rho)$ is regular. Introducing an averaging over translates of $B$, we see that
\begin{equation}\label{eq10.100} |\E_y \E_{x \in y + B} (f(x) e(-\widetilde \phi(x)))| \geq \eta
.\end{equation}
Now write $U := B(S,\half - 10 \rho)$. We wish to exclude from \eqref{eq10.100} those $y$ which lie in the complement $U^{c}$ of $U$, since the bracket functions $\{\xi \cdot x\}$, $\xi \in S$, fail to be linear when $\{\xi \cdot x\} \approx \pm\half$. To this end, we estimate
\begin{eqnarray*}
\big|\E_y \E_{x \in y + B} (f(x) e(-\widetilde{\phi}(x)) 1_{U^c} (y)\big| & \leq & \E_y 1_{U^c}(y) \\ & \leq & \E_y \sum_{\xi \in S} 1_{\| y \cdot \xi \| \geq \half - 10 \rho}(y) \\ & \leq & 40 d\rho \leq 80 d \eps,
\end{eqnarray*}
the penultimate estimate following from the fact that $N$ is prime, so that $\xi \cdot x$ takes on the values $r/N$, $r \in \Z/N\Z$, precisely once each as $x$ varies. Combining this with \eqref{eq10.100} we see that 
\[ |\E_y \E_{x \in y + B} (f(x) e(-\widetilde{\phi}(x)) 1_U(y)| \geq \eta/2,\]
and hence there is $y \in U$ such that 
\begin{equation}\label{eq10.101} |\E_{x \in y + B} (f(x) e(-\widetilde{\phi}(x)))| \geq \eta/2.\end{equation}
We claim that $\widetilde{\phi}$ is a quadratic phase function on $y + B$. To check this, we must show that if the cube $(x + \omega_1 h_1 + \omega_2 h_2 + \omega_3 h_3)_{\omega \in \{0,1\}^3}$ is contained in $y + B$ then $(h_1 \cdot \nabla_x)(h_2 \cdot \nabla_x)(h_3 \cdot \nabla_x) \phi(x) = 0$. This is easy to prove once one appreciates that (for example) if $x, x + h_1 \in y + B$ and $\xi \in S$ then $\{\xi \cdot (x + h_1)\} = \{ \xi \cdot x\} + \{ \xi \cdot h_1\}$. Indeed, this identity is patently true $\md{1}$, and furthermore one has the bounds $|\{ \xi \cdot x\}| \leq \half - 9\rho$ and $|\{\xi \cdot h_1\}| \leq 2\rho$, whence $|\{\xi \cdot (x + h_1)\}| \leq \half - 7\rho$.\vs

\ni Equation \eqref{eq10.101}, then, implies that 
\[ \| f \|_{u^3 (y + B)} \geq \eta/2.\] The result is now an immediate consequence of Theorem \ref{mainthm2} (ii).
\end{proof}

\section{Application: a bound for $r_4(G)$.}\label{sz-g}

\ni As an application of Theorem \ref{mainthm2}, we obtain a bound for $r_4(G)$, the size of the largest set $A \subseteq G$ with no 4-term arithmetic progressions.

\begin{theorem}[Szemer\'edi's theorem for $G$]\label{sz4-g} 
Let $G$ be a finite additive group of order $N$, where $(N,6) = 1$. Then we have the bound 
\[ r_4(G) \ll N(\log \log N)^{-c}\] for some absolute constant $c > 0$.
\end{theorem} 
\ni 
The reader may find it helpful to recall Gowers' argument \cite{gowers-4-aps} in the case $G = \Z/N\Z$ as explained, for example, in \cite{gowers-icm}. Our argument here will be similar, though we must handle torsion in $G$. Gowers did not use the so-called ``symmetry argument'' of \S \ref{ggc}, since he was able to apply the weak inverse theorem, Theorem \ref{weak-inverse}, where we shall apply Theorem \ref{mainthm2}. It is likely that if our only interest was in proving Theorem \ref{sz4-g} then we could do likewise. However, as remarked in the introduction, our work is ultimately directed towards a study of 4-tuples $p_1 < p_2 < p_3 < p_4 \leq N$ of primes in arithmetic progression, and potentially towards the bound $r_4(\Z/N\Z) \ll N(\log N)^{-c}$. For these applications one does need the full strength of Theorem \ref{mainthm2} (when $G = \Z/N\Z$).\vs

\ni The key to the proof of Theorem \ref{sz4-g} is the following density increment result, which is in the spirit of Proposition \ref{sec7-dens-inc} but rather more complicated.

\begin{proposition}\label{sec11-dens-inc} Let $G$ be an abelian group of size $N$, and
suppose that all elements of $G$ have order at most $e^{\sqrt{\log N}}$. Suppose that $(6,N) = 1$, that $N$ is sufficiently large, that $\delta \geq 1/\log \log N$ is smaller than some absolute constant, and that $A \subseteq G$ has size at least $\delta N$. Suppose that $A$ contains no 4-term arithmetic progression. Then there is some subgroup $G' \leq G$, $|G'| \geq N^{1/8}$, together with a coset $t + G'$ such that 
\[ \E_{x \in t + G'}1_A(x) \geq \E_{x \in G} 1_A(x) + \delta^C,\] where $C$ is some absolute constant.
\end{proposition}
\ni\textit{Proof.} 
\ni Write $\alpha := \E 1_A$ and $f := 1_A - \alpha$. Applying Corollary \ref{prog-unif} we have $\|f\|_{U^3(G)} \geq \delta^3/8$; applying Theorem \ref{mainthm2}, 
we can then  find a regular Bohr set $B(S,\rho)$ in $G$ with
$|S| \leq \delta^{-C}$ and $\rho \geq \delta^C$ such that
$$ \E_{y \in G} \| f \|_{u^3(y+B)} \geq \delta^C.$$
Applying Lemma \ref{bohr-cp}, we see that $B$ contains a proper coset progression $P+H$ of rank at most $ \delta^{-C}$
such that
$$ |P+H| \geq \exp(-\delta^{-C} ) N.$$
Since every element in $G$ has order at most $e^{\sqrt{\log N}}$, we see that the lengths of the proper progression $P$ are also at most $e^{\sqrt{\log N}}$.
Thus $|P| \leq N^{1/2}$, which implies that
\begin{equation}\label{eq11.101} |H| \geq N^{1/3}.\end{equation}
By the definition \eqref{eq10.1} of $H$, we see that $B$ can be partitioned into cosets of $H$, and whence
$$ \E_{y \in G} \| f \|_{u^3(y+H)} \geq \delta^C.$$
Set $\alpha'(y) := \E_{x \in y + H} 1_A(x)$.
From the triangle inequality we have
\begin{eqnarray*} \|f\|_{u^3(y+H)} & \leq & \| 1_A - \alpha'(y) \|_{u^3(y + H)} + |\alpha'(y) - \alpha| \\ & = & \| 1_A - \alpha'(y) \|_{u^3(y + H)} - 8 |\alpha'(y) - \alpha| + 9 |\alpha'(y) - \alpha|,\end{eqnarray*}
and so either 
\begin{equation}\label{ayah}
\E_{y \in G} \| 1_A- \alpha'(y) \|_{u^3(y+H)} - 8|\alpha'(y) - \alpha| \geq \delta^C/2
\end{equation}
or
$$ \E_{y \in G} |\alpha'(y) - \alpha| \geq \delta^C/18.$$
Suppose the latter inequality holds. From Lemma \ref{avg-1} we have
$$ \E_{y \in G} (\alpha'(y) - \alpha) = 0;$$
adding this to the preceding estimate and applying the pigeonhole principle, we conclude that there exists $y$ such that 
\begin{equation}\label{density-inc}
 \E_{x \in y + H}1_A(x) = \alpha'(y) \geq \alpha + \delta^C/36 = \E_{x \in G}1_A(x) + \delta^C/36,
 \end{equation}
which implies the proposition (with a change to the absolute constant $C$). Suppose, then, that \eqref{ayah} holds.  By the pigeonhole principle, we
can find $y \in G$ such that
$$\| 1_A(x) - \alpha'(y) \|_{u^3(y+H)}  \geq 8|\alpha'(y) - \alpha| + \delta^C/2.$$
By translating $A$ we may take $y=0$. Writing $\alpha' := \E_{x \in H}1_A(x)$ we conclude the existence of a quadratic phase function $\phi: H \to \R/\Z$ such that 
$$|\E_{x \in H}( f_H(x) e(-\phi(x)) )| \geq 8|\alpha' - \alpha| + \delta^C/2,$$
where $f_H(x) := 1_A(x) - \alpha'$.  
Applying Lemma \ref{quadratic-classify} (or Lemma \ref{inv-quad}), we may thus
find a self-adjoint homomorphism $M: H \to \widehat H$ and $\xi \in \widehat H$ such that
\begin{equation}\label{xdensity}
 |\E_{x \in H} f_H(x) e(-Mx \cdot x) e(-\xi \cdot x)| \geq 8|\alpha' - \alpha| + \delta^C/2.
 \end{equation}
As in the proof of Theorem \ref{sz4-ff}, the next step is to locate a large subgroup of $H$ on which $M$ vanishes.  To achieve this we need some preliminary algebraic (and Fourier-analytic) lemmas, of similar flavour to Lemma \ref{gauss}.

\begin{lemma}[Orthogonal complements]\label{orthocomp}  Let $K$ be any subgroup of $H$, and let $K^\perp \subseteq  H$ be the subgroup
$$ K^\perp := \{ y \in H: Mx \cdot y = 0 \hbox{ for all } x \in K \}.$$
Then $|K^\perp| \geq |H|/|K|$.
\end{lemma}

\ni\textit{Proof.}  Let $\phi: H \to \widehat K$ be the homomorphism $\phi(y)(x) := Mx \cdot y$.  Then $K^\perp$ is precisely the kernel of $\phi$.  But
since $\phi$ is a homomorphism with a domain of size $|H|$ and a range of size at most $|\widehat K| = |K|$, the claim follows.\endproof

\begin{lemma}[Gauss sum lemma]\label{gauss-2}  
Let $K$ be any finite group with $(6,K) = 1$, and suppose that every non-zero element of $K$ has order at most $t$ for some
$t > 2$.  Let $M: K \to \widehat K$ be a self-adjoint homomorphism.  Then, if $|K| \geq 100 t^4$,
there exists a non-zero element $x \in K \backslash \{0\}$ such that $Mx \cdot x = 0$.
\end{lemma}

\ni\textit{Proof.}  We can assume that $M$ is injective (and hence bijective), since otherwise we can just set $x$ to equal a non-zero element in the kernel of $M$.  Using the classification of finite abelian groups, we can write $K$ as the direct sum of cyclic groups of odd prime power order.  Note that
if $p$ is a prime such that at least three cyclic groups of order equal to a power of $p$ appear in this direct sum, then $K$ contains a subgroup isomorphic to $\mathbb{F}_p^3$.  Restricting $M$ to $\mathbb{F}_p^3$ (note that $M$ will still be self-adjoint) and applying Lemma \ref{gauss} we can then conclude the existence of a non-zero $x \in K$ such that $Mx \cdot x = 0$.  Thus we may assume that for each $p \geq 5$ there are at most two cyclic groups of order equal to a power of $p$ in the direct sum decomposition of $K$, in which case we may write
\[ K = \otimes_{j=1}^k (\Z/p_j^{u_j}\Z) \times (\Z/p_j^{u'_j}\Z)\] for some distinct primes $p_1,\dots,p_k$ and exponents $u_j,u'_j$. Let $n \neq 0$ be any integer, and define
\[ n^{\perp} := \{x \in K : nx = 0\}.\]
Writing
\[ n = (-1)^v p_1^{v_1} \dots p_k^{v_k},\]
one confirms the estimate
\begin{equation}\label{11-star} |n^{\perp}| = \prod_{j=1}^k \min(p_j^{u_j}, p_j^{v_j}) \min(p_j^{u'_j},p_j^{v_j}) \leq |n|^2.
\end{equation}
Let $\chi: \R/\Z \to \R$ be a smooth bump function such that $\chi(s) = 1$ when $\|s\|_{\R/\Z} < 1/2t$, $\chi(s) = 0$ when $\| s \|_{\R/\Z} \geq 1/t$, and for which the derivative estimate $\Vert \chi''' \Vert_{\infty} \leq 100/t^3$ holds true.   Observe that if $Mx \cdot x$ is non-zero, then $\|Mx \cdot x\|_{\R/\Z} \geq 1/t$ since $x$ (and hence $Mx \cdot x$) has order at most $t$.
Thus 
$$ \E_{x \in K} 1_{Mx \cdot x = 0} = \E_{x \in K} \chi( Mx \cdot x ).$$ Expanding in a Fourier series gives
$$ \E_{x \in K} 1_{Mx \cdot x = 0} = \sum_{n \in \Z} \widehat \chi(n) \E_x e(n Mx \cdot x) \quad \mbox{where} \quad \widehat\chi(n) = \int_{\R/\Z} \chi(s) e(-ns)\ ds.$$ Isolating the term $n = 0$ we obtain the inequality
$$ \big|\E_{x \in K} 1_{Mx \cdot x = 0} - \widehat \chi(0)\big| \leq \sum_{n \in \Z \backslash \{0\}} |\widehat \chi(n)| |\E_x e(nMx \cdot x)| .$$
Now by \eqref{vdc}, Fourier inversion, the injectivity and self-adjointness of $M$, and the hypothesis that $|K|$ is odd, we have
\begin{eqnarray*}
|\E_x e(nMx \cdot x)| & \leq & |\E_h \E_x e(nM(x+h) \cdot (x+h) - nMx\cdot x)|^{1/2} \\ & \leq &  (\E_h |\E_x e(2nMx \cdot h)|)^{1/2} \\ 
&= & (\E_h 1_{M 2nh = 0})^{1/2} = (\E_h 1_{nh = 0})^{1/2} = \bigg(\frac{|n^{\perp}|}{|K|}\bigg)^{1/2} \leq \frac{n}{\sqrt{K}},
\end{eqnarray*}
the last estimate following from \eqref{11-star}. It follows that 
\begin{equation}\label{eq11.102} \big|\E_{x \in K} 1_{Mx \cdot x = 0} - \widehat \chi(0)\big| \leq \sum_{n \in \Z \backslash \{0\}} |\widehat \chi(n)| |n| |K|^{-1/2} .\end{equation}
Now by integrating by parts three times and using the bound on $\| \chi''' \|_{\infty}$ one sees that $|\widehat{\chi}(n)| \leq t^2/|n|^3$. In combination with the trivial bound $\| \widehat \chi \|_{\infty} \leq 2/t$, we obtain
\[ \sum_{n \neq 0} |\widehat{\chi} (n)| |n| \leq \sum_{|n| < t} \frac{|n|}{t} + \sum_{|n| \geq t} \frac{t^2}{n^2} \leq 6t.\]
Thus, since $|K| \geq 100t^4$, \eqref{eq11.102} implies that
\[\E_{x \in K} 1_{Mx \cdot x = 0} \geq \frac{1}{3t} > \frac{1}{|K|},\] which immediately implies the result.\endproof
 
\begin{corollary}  Let $H$ be a finite additive group, $(|H|,6) = 1$, such that every element has order at most $t$, $t \geq 2$, and let $M: H \to \widehat H$ be a self-adjoint homomorphism.  Then there exists a subgroup $K$ of $H$ such that 
\begin{equation}\label{eq11.100}
|K| \geq |H|^{1/2}/10t^2
\end{equation} and $Mx \cdot y = 0$
for all $x,y \in K$.
\end{corollary}

\begin{proof}  Let $K$ be a subgroup of $H$ on which the quadratic form $Mx \cdot y$ vanishes (i.e. $Mx \cdot y = 0$ for all $x,y \in K$), and which
is maximal with respect to set inclusion.  Observe that the orthogonal complement $K^\perp$ of $K$ contains $K$, and hence by Lemma \ref{orthocomp}
the quotient group $K^\perp/K$ has cardinality at least $|H| / |K|^2$.  Also, every element in this group has order at most $t$.  Since $Mx \cdot y = My \cdot x = 0$ whenever $x \in K$ and $y \in K^\perp$ we see that the bilinear form $Mx \cdot y$ descends to a bilinear form
on $K^\perp/K$.  If the associated quadratic form vanished for at least one non-zero element of $K^\perp/K$, then by adjoining this element to $K$ we could contradict the maximality of $K$.  Thus we may assume that there is no such form. But then by Lemma \ref{gauss-2} we have
$|K^\perp/K| \leq 100 t^4$.  Combining this with our lower bound for $|K^\perp/K|$ we obtain the result.
\end{proof}\vs

\ni Let us return now to the situation \eqref{xdensity}, and let $K$ be the subgroup obtained by the above Corollary.
By Lemma \ref{avg-1} we have
$$
 |\E_{y \in H} \E_{x \in y+K} f_H(x) e(-Mx \cdot x) e(-\xi \cdot x)| \geq 8|\alpha' - \alpha| + \delta^C/2.
$$ Setting $\alpha''(y) := \E_{x \in y + K}1_A(x)$, the triangle inequality implies that either
\begin{equation}\label{eq11.3a} \E_{y \in H} |\alpha''(y) - \alpha'| \geq 2|\alpha' - \alpha| + \delta^C/12\end{equation}
or
\begin{equation}\label{ayeka}
 \E_{y \in H} \big(\;|\E_{x \in y+K} f_H(x) e(-Mx \cdot x) e(-\xi \cdot x)| - 3|\alpha''(y)-\alpha'| \; \big) \geq 2|\alpha' - \alpha| + \delta^C/4.
\end{equation}
Suppose that \eqref{eq11.3a} holds.  From Lemma \ref{avg-1} we have
$$ \E_{y \in H} (\alpha''(y) - \alpha) = 0$$
and hence
\[ 2 \E_{y \in H} \max(\alpha''(y) - \alpha',0) = \E_y |\alpha''(y) - \alpha'| \geq 2|\alpha' - \alpha| + \delta^{C}/12.\]
By the pigeonhole principle we thus conclude that there exists $y \in H$ such that
$$ 2 (\alpha''(y) - \alpha') \geq 2|\alpha' - \alpha| + \delta^C/12,$$
and hence by the triangle inequality
\begin{equation}\label{density-inc-2}\E_{x \in y + K}1_A(x) = \alpha''(y) \geq \alpha' + |\alpha' - \alpha| + \delta^C/24 \geq \alpha + \delta^C/24 = \E_{x \in G} 1_A(x) + \delta^C/24.
 \end{equation}
This, together with the lower bound
\[ |K| \geq |H|^{1/2}/10 e^{2\sqrt{\log N}} > N^{1/8}\]
(cf. \eqref{eq11.100}) implies the proposition under the assumption that \eqref{eq11.3a} holds.\vs

\ni Suppose, then, that 
\eqref{ayeka} holds instead.  By the pigeonhole principle, we can find $y \in H$
such that
$$ |\E_{x \in y+K} f_H(x) e(-Mx \cdot x) e(-\xi \cdot x)| \geq 3 |\alpha''(y) - \alpha'| + 2|\alpha' - \alpha| + \delta^C/4.$$
Splitting $f_H$ as $(1_A - \alpha''(y)) + (\alpha''(y) - \alpha')$ and using the triangle inequality, we conclude
$$ |\E_{x \in y+K} ((1_A(x) - \alpha''(y)) e(-Mx \cdot x) e(-\xi \cdot x)| \geq 2 |\alpha''(y) - \alpha| + \delta^C/4.$$
We write $x = y+z$ and use the fact that the bilinear form $My \cdot z$ is symmetric and vanishes on $K$ to conclude that
$$ |\E_{z \in K} ((1_A(y+z) - \alpha''(y)) e(-(2My+\xi) \cdot z)| \geq 2 |\alpha''(y) - \alpha| + \delta^C/4.$$
Consider now the homomorphism $\phi: K \to \R/\Z$ defined by $\phi(x) := (2My + \xi) \cdot x$.  A simple pigeonhole argument shows that
there exist at least $|K|e^{-\sqrt{\log N}}$ elements $x$ of $K$ for which $\| \phi(x) \|_{\R/\Z} < e^{\sqrt{\log N}}$.  But since every element of $K$ has order at most $e^{\sqrt{\log N}}$, we conclude that $\phi(x) = 0$.  Thus if we set $K'$ to be the kernel of $\phi$, then $K'$ is a subgroup of $K$ with
\begin{equation}\label{eq11.102A}
|K'| \geq  |K|e^{-\sqrt{\log N}} > N^{1/8}.\end{equation}  We then apply Lemma \ref{avg-1} again to conclude that
$$ |\E_{w \in K} \E_{z \in w+K'} ((1_A(y+z) - \alpha''(y)) e(-(2My+\xi) \cdot z)| \geq 2 |\alpha''(y) - \alpha| + \delta^C/4.$$
Since the phase $e(-(2My+\xi) \cdot z)$ is constant for $z \in w + K$, for fixed $w$, we conclude
$$ \E_{w \in K} |\E_{z \in w+K'} \big(1_A(y+z) - \alpha''(y)\big)| \geq 2 |\alpha''(y) - \alpha| + \delta^C/4$$
and thus, writing $\alpha'''(t) := \E_{x \in t + K'}1_A(x)$, 
\begin{equation}\label{eq11.5a} \E_{w \in K} |\alpha'''(w+y) - \alpha''(y)| \geq 2 |\alpha''(y) - \alpha| +  \delta^C/4.\end{equation}
Now from Lemma \ref{avg-1} we have
$$ \E_{w \in K} \big(\alpha'''(w + y) - \alpha''(y)\big) = 0.$$
Together with \eqref{eq11.5a} this implies that 
\[ 2 \E_{w \in K} \max(\alpha'''(w + y) - \alpha''(y),0) \geq 2|\alpha''(y) - \alpha| + \delta^C/4,\]
and so there exists $w \in K$ such that
$$ \alpha'''(w+y) - \alpha''(y) \geq |\alpha''(y) - \alpha| + \delta^C/8.$$
One final application of the triangle inequality gives at last that \[
\E_{x \in w + y + K'}1_A(x) = \alpha'''(w + y) \geq  \alpha + \delta^C/8 = \E_{x \in G}1_A(x) + \delta^C/8.
\]
Together with the lower bound \eqref{eq11.102A}, this concludes the proof of Proposition \ref{sec11-dens-inc}.\endproof\vs

\ni\textit{Proof of Theorem \ref{sz4-g}.} Let $G$ be an abelian group, and let $A \subseteq G$ be a set with cardinality at least $\delta N$ which contains no 4 distinct elements in arithmetic progression. We wish to show that $\delta \ll (\log \log N)^{-c}$, for some absolute constant $c$; thus we may certainly suppose that $\delta \geq 2/\log \log N$. \vs

\ni Suppose that $G$ has an element $g$ of order greater than $\exp((\log N)^{1/3})$. Writing $H = \langle g \rangle$, we see that there is some coset $y + H$ such that $\E_{x \in y + H} 1_A(x) \geq \delta$. The result of Gowers \cite{gowers-4-aps} then immediately implies that 
\[ \delta \ll (\log \log |H|)^{-c} \ll (\log \log N)^{-c}.\]
Suppose, then, that all elements of $G$ have order at most $\exp((\log N)^{1/3})$. We will define a sequence $G = G_0 \geq G_1 \geq \dots $ of subgroups of $G$ with cardinalities $N = N_0 \geq N_1 \geq \dots $.
The sequence will be defined in such a way that 
\begin{equation}\label{eq11.103} N_j \geq \exp ((\log N)^{2/3}),\end{equation} which means that no element in $G_j$ has order greater than $\exp((\log N_j)^{1/2})$, and also that $\delta \geq 1/\log\log N_j$. This is to enable us to apply Proposition \ref{sec11-dens-inc}.\vs

\ni Suppose that we have defined $G_j$. For any $x_j$, the set $(x_j + A) \cap G_j$ does not contain a 4-term arithmetic progression. Suppose that $x_j$ is such that $\E_{x \in x_j + G_j} 1_A(x) \geq \delta$. Then, applying Proposition \ref{sec11-dens-inc}, we see that there is some $G_{j+1}$, $N_{j+1} := |G_{j+1}| \geq N_j^{1/8}$, together with some $x_{j+1}$ so that 
\[ \E_{x \in x_{j+1} + G_{j+1}}1_A(x) \geq \E_{x \in x_j + G_j}1_A(x) + \delta^C.\]
Iterating this construction leads to a contradiction for some $j \leq \delta^{-C}$ unless \eqref{eq11.103} is violated. Since $N_j \geq N^{(1/8)^j}$, we must therefore have
\[ N^{(1/8)^{\delta^{-C}}} < \exp ((\log N)^{2/3}),\]
which implies the required bound $\delta \ll (\log \log N)^{-c}$. \endproof\vs

\ni\textit{Remarks.} We hope to prove a bound of the form $r_4(G) \ll N(\log N)^{-c}$ in a future paper by combining the ideas of \cite{green-tao-ffszem} with nested Bohr set technology in the spirit of that used in \S \ref{sec8} and \S \ref{ggc} of the present paper. These methods were first introduced by Bourgain \cite{Bou}, who obtained the bound $r_3(G) \ll N(\log N)^{-1/2 + \epsilon}$, which is still the best currently known when $G = \Z/N\Z$. A feature of this approach is that, unlike in the present section, it is no easier to deal with $\Z/N\Z$ than it is with an arbitrary abelian $G$. We note that the celebrated Erd\H{o}s-Tur\'an conjecture \cite{erdos} is roughly equivalent to a bound of the form $r_k(\Z/N\Z) \ll_k N(\log N)^{-1 + \epsilon}$, and so even in the case $k = 3$ there is an awful lot left to be done.

\section{An ergodic theory interpretation}\label{ergodic-sec}

\ni We now connect the $U^3$ inverse theorems discussed earlier to ergodic theory, and in particular to the recent work of Host-Kra
\cite{host-kra2} and Ziegler \cite{ziegler}.\vs

\ni Define a \emph{measure-preserving system} $(X, {\B}, T, \P)$ to be a probability space
$(X, {\B}, \P)$ with an invertible measure-preserving (i.e. probability-preserving) shift operator $T: X \to X$.  This 
induces a shift operator $T$ on random variables $f: X \to \R$ by the formula
$Tf(x) := f(T^{-1} x)$, and more generally $T^n f(x) := f(T^{-n} x)$ for any $n \in \Z$.  
We use $\E_X(f)$ to denote the expectation of $f$.\vs

\ni The \emph{Furstenberg correspondence principle} (see e.g. \cite{furst-book}) equates combinatorial theorems such as Szemer\'edi's theorem to 
recurrence results in ergodic theory.  In particular, Theorem \ref{szemeredi} is logically equivalent (using the axiom of choice) to the following theorem.

\begin{theorem}[Furstenberg recurrence theorem \cite{furst,FKO}]  Let $(X, {\B}, T, \P)$ be a measure-preserving system, and let
$f \in L^\infty(X, {\B})$ be any non-negative random variable with $\E_X f > 0$.  Then for every $k \geq 1$ we have
$$ \liminf_{N \to \infty} \E_{-N \leq n \leq N} \E_X f T^n f \ldots T^{(k-1)n f} > 0.$$
In particular, if $A \in {\B}$ is any event with positive probability $\P(A) > 0$, then
$$ \liminf_{N \to \infty} \E_{-N \leq n \leq N} \P(A \cap T^n A \cap \ldots \cap T^{(k-1)n} A) > 0.$$
\end{theorem}

\ni Recently, it was shown in \cite{host-kra2} and \cite{ziegler} that this limit inferior can in fact be replaced by a limit; earlier work related to the $k=4$ case can be found in  \cite{conze,furstenberg-weiss,host-kra1,hk-cubes}.  The two approaches are slightly different; the argument in \cite{host-kra2} proceeds by establishing the ergodic theory analogue of an inverse theorem for
the Gowers uniformity norm $U^d(\Z/N\Z)$.  Indeed, if $f \in L^\infty(X, {\B})$ is any complex-valued random variable, define the quantity $\|f\|_{U^d(T)}$ for $d \geq 0$ by the formula
$$ \|f\|_{U^d(T)} := \lim_{N \to \infty}
\big(\E_{-N \leq h_1,\ldots,h_d \leq N} \E_X \prod_{\omega \in \{0,1\}^d} {\Conj}^{|\omega|} T^{\omega \cdot h} f(x)\big)^{1/2^d}.$$
It can be shown \cite{hk-cubes,host-kra2} that this limit actually exists, and it can also be shown that $\|f\|_{U^d(T)}$ is in fact
a semi-norm on bounded random variables for any $d \geq 1$; see \cite{host-kra2}.  This semi-norm is clearly related to the $U^d(G)$ norms defined in 
Definition \ref{gundef}.  For instance, if $T$ is periodic of order $N$ then it is easy to see that
the $U^d(T)$ norm of $f(x)$ is the $L^{2^d}$ average of the $U^d(\Z/N\Z)$ norm of $f$ restricted to the orbits of $T$:
$$ \|f\|_{U^d(T)}^{2^d} = \E_x \| (T^h f(x))_{h \in \Z/N\Z} \|_{U^d(\Z/N\Z)}^{2^d}.$$
Here of course we take advantage of the periodicity of $T$ to define $T^h$ for $h \in \Z/N\Z$ in the obvious manner. \vs

\ni In \cite{host-kra2} it was observed that this semi-norm controls expressions such as those appearing in the Furstenberg recurrence theorem.  Indeed,
there is an analogue of Proposition \ref{gvn} which asserts that if $f_0,\ldots,f_{k-1}$ are bounded random variables and at least one of them has vanishing $U^{k-2}(T)$ norm, then
$$  \lim_{N \to \infty} \E_{-N \leq n \leq N} \E_X f_0 T^n f_1 \ldots T^{(k-1)n} f_n = 0;$$
in fact one can make the slightly stronger claim that $\E_{-N \leq n \leq N} T^n f_1 \ldots T^{(k-1)n} f_n$ converges to zero in (for instance)
the $L^2(X)$ sense.  Informally, this fact shows that functions with vanishing $U^{k-2}(T)$ norm are irrelevant for understanding $k$-fold recurrence.\vs

\ni It is thus of interest to determine when the $U^{k-2}(T)$ norm is positive.  This question is answered in \cite{host-kra2} using the language of nilsystems.  We first recall some notation.  If $G$ is a (not necessarily abelian) group written multiplicatively and if $g, h \in G$, we let $[g,h] := g^{-1} h^{-1} gh$ 
denote the commutator of $g$ and $h$.  If $G'$ and $G''$ are subgroups of $G$, we let $[G',G''] = [G'',G']$ be the subgroup generated by 
the commutators $\{ [g',g'']: g' \in G', g'' \in G'' \}$.  We then define the \emph{lower central series}
$$ G = G_1 \supseteq G_2 \supseteq G_3 \supseteq \ldots$$
of subgroups of $G$ by the recursive definition $G_1 := G$; $G_{k+1} := [G,G_k]$.  We say that $G$ is \emph{$(k-2)$-step nilpotent} for some $k \geq 2$
if $G_{k-1}$ is trivial.  Thus for instance a group is $1$-step nilpotent if and only if it is abelian.\vs

\ni A \emph{$(k-2)$-step nilmanifold} is defined to be a manifold of the form $G/\Gamma := \{ x\Gamma: x \in G \}$, where $G$ is a finite-dimensional nilpotent Lie group, and
$\Gamma$ is a discrete subgroup of $G$ which is co-compact (i.e. the nilmanifold $G/\Gamma$ is compact).  Note that we do not assume $\Gamma$ to be normal, and hence a nilmanifold need not have a group structure.  It is however a compact symmetric space, with a left-action of the group $G$.  Thus there is a unique invariant Haar measure $\P$ on a nilmanifold, which we normalize to be a probability measure.  Thus every nilmanifold is a probability space, taking the $\sigma$-algebra to be the Borel $\sigma$-algebra.\vs

\ni If $g \in G$ then we write $T_g$ for the \textit{shift operator} from $G$ to itself defined by $T_g(x) = gx$, and also (by abuse of notation) for the map from $G/\Gamma$ to itself defined by $T_g(x\Gamma) := gx \Gamma$. This latter map is measure-preserving and invertible.  Let us call a $(k-2)$-step nilmanifold with one of these shift operators a \emph{$(k-2)$-step nilflow}.  A \emph{$(k-2)$-step nilfunction} is defined to be any continuous function $F: G/\Gamma \to \C$ on a $(k-2)$-step nilmanifold; given such a nilfunction, a point $x_0 \in G/\Gamma$ and a group element $g \in G$, we define the associated \emph{basic $(k-2)$-step nilsequence} $F_{g,x_0}: \Z \to \C$ by the formula\footnote{A general $(k-2)$-step nilsequence is defined as the uniform limit of basic $(k-2)$-step sequences; see \cite{bhk} for further analysis of these nilsequences.}
$$ F_{g,x_0}(n) := F(T_g^n x_0) \hbox{ for all } n \in \Z.$$
We can truncate this to $\Z/N\Z$ and define the truncated nilsequence $F_{N,g,x_0}: \Z/N\Z \to \C$ by the formula
$$ F_{N, g,x_0}(n) := F_{g,x_0}(n) = F(T_g^n x_0) \hbox{ for all } -N/2 < n \leq N/2 $$
where we identify the integers from $-N/2$ to $N/2$ with $\Z/N\Z$ in the usual manner.\vs

\ni We now give three key examples of nilflows and nilsequences. 

\begin{example}[The circle nilflow]\label{nil1}  Let $G$ be the one-dimensional matrix group
$$ G := \left( \begin{array}{ll}
1 & \R \\
0 & 1 
\end{array} \right) := \{ \left( \begin{array}{ll}
1 & x \\
0 & 1 
\end{array} \right) : x \in \R \}$$
and let $\Gamma$ be the discrete subgroup
$$ \Gamma := \left( \begin{array}{ll}
1 & \Z \\
0 & 1 
\end{array} \right) := \{ \left( \begin{array}{ll}
1 & n \\
0 & 1 
\end{array} \right) : n \in \Z \}.$$
Then $G/\Gamma$ is a 1-step nilmanifold (and hence also a 2-step nilmanifold), 
indeed we can easily identify it with the unit circle $\R/\Z$.  A shift $T_g$ on this 
nilmanifold then corresponds to
a simple translation $x \mapsto x+\alpha$, where $\alpha \in \R$ is the upper right matrix entry of $g$.
In particular, we observe that if $F: \R/\Z \to \C$ is any function, we see that the sequence
$$ n \mapsto F(T^n_g x) = F(x + n\alpha)$$
is a basic 1-step nilsequence. Thus, for instance, the linear phase function $n \mapsto e(n\alpha)$ is a basic $1$-step nilsequence (and hence also a $2$-step nilsequence).  More generally, any quasiperiodic sequence is a basic $1$-step nilsequence, and any almost periodic sequence can be expressed as the uniform limit of basic $1$-step nilsequences.
\end{example} 

\begin{example}[The skew shift nilflow]\label{nil2} Now we consider the example
$$ G := \left( \begin{array}{lll}
1 & \Z & \R\\
0 & 1  & \R\\
0 & 0  & 1
\end{array} \right); \quad
\Gamma := \left( \begin{array}{lll}
1 & \Z & \Z\\
0 & 1  & \Z\\
0 & 0  & 1
\end{array} \right). $$
Then $G/\Gamma$ is a 2-step nilmanifold, and one can identify it topologically with the 2-torus $(\R/\Z)^2$ by the identification
$$ (x,y) \equiv \left( \begin{array}{lll}
1 & 0  & y\\
0 & 1  & x\\
0 & 0  & 1
\end{array} \right) \Gamma.$$
If we let
$$ g := \left( \begin{array}{lll}
1 & m  & \beta \\
0 & 1  & \alpha \\
0 & 0  & 1
\end{array} \right) $$
be a typical element of $G$ (thus $m \in \Z$ and $\alpha,\beta \in \R$) then the shift $T_g$ is then given by
$(x,y) \mapsto (x + \alpha, y + \beta + mx)$, and thus if $F: (\R/\Z)^2 \to \C$ is any function then the sequence
$$ n \mapsto F(T^n_g (x,y)) = F( x + n\alpha, y + n\beta + \half m n(n+1) \alpha )$$
is a basic $2$-step nilsequence.  Thus, for instance, the quadratic phase function $n \mapsto e( \half \alpha n(n+1))$ is a basic $2$-step nilsequence.  More generally, any quadratic phase $n \mapsto e( \alpha n^2 + \beta n + \gamma)$, or finite linear combination of such phases,
is a basic $2$-step nilsequence.
\end{example}

\begin{example}[The Heisenberg nilflow]\label{nil3}  Now we consider the example
$$ G := \left( \begin{array}{lll}
1 & \R & \R\\
0 & 1  & \R\\
0 & 0  & 1
\end{array} \right); \quad
\Gamma := \left( \begin{array}{lll}
1 & \Z & \Z\\
0 & 1  & \Z\\
0 & 0  & 1
\end{array} \right). $$
Then $G/\Gamma$ is a 2-step nilmanifold.  By using the identification
$$ (x,y,z) \equiv \left( \begin{array}{lll}
1 & z  & y\\
0 & 1  & x\\
0 & 0  & 1
\end{array} \right) \Gamma,$$
we can identify $G/\Gamma$ (as a set) with $\R^3$, quotiented out by the equivalence relations
$$ (x,y,z) \sim (x+a, y+b+az, z+c) \hbox{ for all } a,b,c \in \Z.$$
This can in turn be coordinatized by the cylinder $[-1/2,1/2] \times (\R/\Z)^2$ with the identification
$(-1/2,y,z) \sim (1/2,y+z,z)$.\vs

\ni Let $F : G/\Gamma \rightarrow \C$ be a function. We may lift this to a function $\widetilde{F} : G \rightarrow \C$, defined by $\widetilde{F}(g) := F(g\Gamma)$. In coordinates, this lift takes the form
\[ \widetilde{F}(x,y,z) = F(\{x\}, y - [x] z \md{1}, z \md{1})\]
where $[x] = x - \{x\}$ is the nearest integer to $x$ (we round half-integers up).
If we let
$$ g := \left( \begin{array}{lll}
1 & \gamma  & \beta \\
0 & 1  & \alpha \\
0 & 0  & 1
\end{array} \right) $$
be an element of $G$, then the shift $T_g : G \rightarrow G$ is given by
\[ T_g(x,y,z) = (x + \alpha, y + \beta + \gamma x, z + \gamma),\]
from which a short induction confirms that 
\[ T_g^n(x,y,z) = (x + n\alpha, y + n\beta + \half n(n+1)\alpha , z + n\gamma).\]
 Therefore if $F : G/\Gamma \rightarrow G/\Gamma$ is any function, written as a function $F: [-\half,\half] \times (\R/\Z)^2 \to \C$ 
with $F(-1/2,y,z) = F(1/2,y+z,z)$, then we have
$$F(T^n_g (x,y,z)) = F( \{x + n\alpha\}, y + n\beta + \half n(n+1) \alpha\gamma - [x+n\alpha] (z+n\gamma), z + n \gamma).$$ This, of course, is a basic $2$-step nilsequence.  We see, for instance, that the generalized quadratic phase 
function $n \mapsto e( \half n(n+1) \alpha\gamma - [n\alpha] n\gamma)$ is a basic $2$-step nilsequence.
\end{example}

\ni We call the three basic examples just discussed the \textit{fundamental} $2$-step nilsequences. They may be used in a straightforward product construction to construct further nilsequences, as we now describe.\vs

\ni If $(G/\Gamma,T_g)$ and $(G'/\Gamma,T_{g'})$ are 2-step nilflows, then so is the
direct sum $( (G \oplus G')/(\Gamma \oplus \Gamma'), T_{(g,g')})$.  Also, if $F: G/\Gamma \to \C$ and $F': G'/\Gamma' \to \C$ are functions,
and we define the tensor product $F\otimes F': (G \oplus G')/(\Gamma \oplus \Gamma')$ in the usual manner as
$$ F \otimes F'(x,x') := F(x) F'(x')$$
then we see that the function $F \otimes F'(T_{(g,g')}^n (x,x') )$ factors as
$$ F \otimes F'(T_{(g,g')}^n (x,x') ) = F(T_g^n x) F'(T_{g'}^n x').$$

\ni Now suppose we take $n_1$ nilsequences coming from circle nilflows, $n_2$ nilsequences coming from skew shift nilflows, and $n_3$ nilsequences coming from Heisenberg nilflows, and tensor them all together. What results is a nilsequence on a $2$-step nilmanifold $G/\Gamma$, which is topologically the cube $[-\half,\half]^{n_1 + 2n_2 + 3n_3}$ with faces identified. $2$-step nilsequences of this type, that is to say tensor products of fundamental nilsequences, are in a sense the only important ones if one is interested in the Gowers $U^3$ norm. We call them the \textit{elementary $2$-step nilsequences} (we also refer to \textit{elementary $2$-step nilmanifolds} and \textit{elementary $2$-step nilflows}).\vs

\ni Given an elementary $2$-step nilsequence, it is natural to refer to $n_1 + 2n_2 + 3n_3$ as its \textit{dimension}. It is also of interest to have a notion of how continuous the underlying function $F : G/\Gamma \rightarrow G/\Gamma$ is. We adopt a rather low-brow approach to this concept which is sufficient for our purposes. Let $\delta = 1/m$ be the reciprocal of an integer. Then we may divide $(-\half,\half)^d$ into $m^d$ cubes of sidelength $\delta$, for any $d$, and in fact these subdivisions respect the quotienting of Examples \ref{nil1}, \ref{nil2} and \ref{nil3} (for $d = 1,2,3$ respectively), giving what we refer to as $\delta$-\textit{nets} on the three fundamental $2$-step nilmanifolds $G/\Gamma$. We refer to the building blocks of these nets as the $\delta$-atoms. Taking products, we may obtain $\delta$-atoms and a $\delta$-net on any elementary $2$-step nilmanifold $G/\Gamma$. Finally, if $F : G/\Gamma \rightarrow \C$ is a function and if $K > 0$ is a constant, we say that $F$ is $K$-Lipschitz if for all $\delta = 1/m$ and for all $\delta$-atoms $A$ we have
\[ |F(x) - F(x')| \leq K\delta\]
whenever $x,x' \in A$. The following lemma is straightforward.

\begin{lemma}\label{prod-lipschitz} Suppose that $F_i : G_i/\Gamma_i \rightarrow \D$, $i = 1,\dots,k$, are $K$-Lipschitz functions on elementary $2$-step nilmanifolds $G_i/\Gamma_i$. Then the tensor product $F_1 \otimes \dots \otimes F_k$ is $Kk$-Lipschitz.\endproof
\end{lemma}

\ni Now it has been known since the work of Furstenberg and Weiss \cite{furstenberg-weiss}
that random variables on a $(k-2)$-step nilflow $(G/\Gamma, T_g)$ have a non-trivial behavior with respect to $k$-term recurrence.  Indeed, given any bounded random variable $f_0: G/\Gamma \to \C$ which is not identically zero, one can find $f_1,\ldots,f_{k-1}$ such that
the averages \[ \E_{-N \leq n \leq N} \E_X f_0 T^n f_1 \ldots T^{(k-1)n} f_n\] do not converge to zero; this is basically due to non-trivial algebraic relations between the $k$ points $x\Gamma$, $g^n x\Gamma, \ldots, g^{(k-1)n} x\Gamma$.  In particular, the $U^{k-1}(T)$ norm is non-degenerate on
this nilmanifold (and is thus a genuine norm); see \cite{host-kra2,ziegler} for some further discussion of this fact.  We will prove a variant of this statement.

\begin{proposition}[Nilsequences obstruct uniformity]\label{tgnx}
Let $k \geq 3$, and let $(G/\Gamma,T)$ be a $(k-2)$-step nilsystem. Let $F: G/\Gamma \to \C$ be a continuous function on $G/\Gamma$ which is not identically zero. Suppose that $N > k - 1$ is a prime, and that $f : \Z/N\Z \rightarrow \D$ is a function such that 
\[ |\E_{-N/2 \leq n \leq N/2} f(n) \overline{F(T_g^n x)}| \geq \eta.\]
Then we have 
\[ \| f \|_{U^3} \geq c_{F, G/\Gamma}(\eta) > 0\]
uniformly in $x \in G/\Gamma, g \in G$ and $N$.
\end{proposition}
\ni\textit{Proof.} In proving this proposition we will use the following lemma to the effect that 
the point $T_g^{n+(k-1)r} x$ is completely constrained
by the cosets $T_g^n x \Gamma, T_g^{n+r} x \Gamma, \ldots, T_g^{n+(k-2)r} x \Gamma$. 
\begin{lemma}\label{constraint-lem}
Let $(G/\Gamma, T)$ be a $(k-2)$-step nilsystem. Then there is a compact set $\Sigma \subseteq (G/\Gamma)^{k-1}$ and a continuous function $P: \Sigma \to G$ such that for all $n,r \in \Z$, $g \in G$ and $x \in G/\Gamma$ we have
\[ (T_g^n x, T_g^{n+r} x, \ldots, T_g^{n+(k-2)r} x) \in \Sigma\]
and 
\[
P(T_g^n x, T_g^{n+r} x, \ldots, T_g^{n+(k-2)r} x) = T_g^{n+(k-1)r} x \Gamma .\]
\end{lemma}
\ni The existence of a constraint of the type is discussed in several places in the ergodic theory literature \cite{bhk,furstsurvey1,furstsurvey2,ziegler02}. For the convenience of the reader we supply a self-contained proof in Appendix \ref{appendix-sec}.\vs

\ni The proof of Proposition \ref{tgnx} is not dissimilar to that of Theorem \ref{mainthm2} (ii), which was given at the start of \S \ref{ggc}, but here we use Lemma \ref{constraint-lem} in place of \eqref{constraint}, and the technical details are rather different.\vs

\ni Assume without loss of generality that $\| F \|_{\infty} \leq 1$, and let $\eps \leq 1/100k$ be a small constant (it will be chosen to be a small multiple of $\eta$). Observe that the function $F \circ P: \Sigma \to \C$ is continuous, hence uniformly continuous, on the compact set $\Sigma$.  In particular
we can find a neighbourhood $V \subseteq G$ of the identity $1$ which depends on $\eps$, $F$, $G/\Gamma$ such that
\begin{equation}\label{f-approx}
F( P(x_0,\ldots,x_{k-2}) ) = F(P(y_0,\ldots,y_{k-2})) + O(\eps) 
\end{equation}
whenever $(x_0,\ldots,x_{k-2}), (y_0,\ldots,y_{k-2}) \in \Sigma$
and $x_j \in V y_j$ for all $0 \leq j \leq k-2$.
Applying Lemma \ref{constraint-lem} and exploiting compactness again, we conclude that there exists another neighbourhood $V' \subseteq  V$ of the identity 1
such that given any $z_1,\ldots,z_{k-2} \in G/\Gamma$, there exists a bounded function $Q_{z_1,\ldots,z_{k-2}}: G/\Gamma \to \D$
such that
$$
F( T_g^{n+(k-1)r} x ) = Q_{z_1,\ldots,z_{k-2}} (T_g^n x) + O(\eps) \hbox{ whenever } T_g^{n+jr}x \in V' z_j \hbox{ for all } 1 \leq j \leq k-2.
$$
It is not hard to ensure that the function $Q_{z_1,\ldots,z_{k-2}}$ depends in a measurable manner on $z_1,\ldots,z_{k-2}$.
In particular we see that
\begin{eqnarray}\nonumber
& & F( T_g^{n+(k-1)r} x ) f(n+(k-1)r) \prod_{j=1}^{k-2} 1_{V' z_j}(T_g^{n+jr} x)  \\ \nonumber & & 
\qquad\qquad = Q_{z_1,\ldots,z_{k-2}}(T_g^n x) f(n+(k-1)r) \prod_{j=1}^{k-2} 1_{V' z_j}(T_g^{n+jr} x) \\
& & \qquad\qquad\qquad\qquad + O(\eps F( T_g^{n+(k-1)r} x ) f(n+(k-1)r) \prod_{j=1}^{k-2} 1_{V'z_j}(T_g^{n+jr} x) .\label{eq10.102}\end{eqnarray}
for all $n,r \in \Z$, $g \in G$ and $x \in G/\Gamma$. We are going to average this over $n$ and $r$, but for technical reasons related to the difference between $\{n : |n| \leq N/2\}$ and $\Z/N\Z$ as additive objects, we shall restrict the range of $r$. To this end, take a function $\chi : \Z/N\Z \rightarrow \R^+$ such that $\E_r \chi(r) = 1$, $\mbox{Supp}(\chi) \subseteq [-\eps N, \eps N]$ and for which we have the Fourier estimate
\begin{equation}\label{eq10.103} \sum_{m \in \Z/N\Z} |\widehat{\chi}(m)| \leq 100/\eps.\end{equation}
Such a function can easily be constructed, for example by convolving an interval with itself. Averaging \eqref{eq10.102} over $|n| \leq N/2 - \eps k N$ and over $r$ weighted by $\chi$, one obtains
\begin{eqnarray}\nonumber
& & \E_{|n| \leq N/2 - \eps k N} \E_{r \in \Z/N\Z} \chi(r) F( T_g^{n+(k-1)r} x ) f(n+(k-1)r) \prod_{j=1}^{k-2} 1_{V' z_j}(T_g^{n+jr} x)  \\ \nonumber & & 
= \qquad \E_{|n| \leq N/2 - \eps k N} \E_{r \in \Z/N\Z}\chi(r) h(n)  f(n+(k-1)r) \prod_{j=1}^{k-2} h_j(n + jr) \\ & & \qquad\qquad\qquad\qquad + O(\eps  \prod_{j=1}^{k-2} 1_{V'z_j}(T_g^{n+jr} x)) ,\label{eq10.104}\end{eqnarray}
where \begin{equation}\label{eq10.105} h(n) := Q_{z_1,\dots,z_{k-2}} (T_g^n x) \qquad \mbox{and} \qquad h_j(n) := 1_{V'z_j}(T_g^n x).\end{equation}
Note that as a consequence of the restrictions we have made on the support of $n$ and of $r$, the expressions $n + jr$, $j = 0,1,\dots, k-1$ are the same whether we regard the addition as taking place in $\Z$ or in $\Z/N\Z$. In particular, \eqref{eq10.104} remains valid if one imagines that these additions are made in $\Z/N\Z$.\vs

\ni Now the second line of \eqref{eq10.104} can be written as 
\[ \E_{|n| \leq N/2} \E_{r \in \Z/N\Z}\chi(r) \widetilde{h}(n)  f(n+(k-1)r) \prod_{j=1}^{k-2} h_j(n + jr),\]
where 
\[ \widetilde{h}(n) := (1 - 2\eps k)^{-1}h(n) 1_{|n| \leq N/2 - \eps k N},\] and in particular $\Vert \widetilde{h} \Vert_{\infty} \leq 2$. Writing $\chi(r)$ in terms of its Fourier transform on $\Z/N\Z$ and using Proposition \ref{gvn}, we can bound this expression above as follows, where $e_N(r) := e(r/N)$:
\begin{eqnarray*} 
&& \big|\E_{|n| \leq N/2} \E_{r \in \Z/N\Z}\chi(r) \widetilde{h}(n)  f(n+(k-1)r) \prod_{j=1}^{k-2} h_j(n + jr)\big| \\
&= & \big|\E_{n \in \Z/N\Z} \E_{r \in \Z/N\Z} \sum_{m \in \Z/N\Z} \widehat{\chi}(m) e_N(-mr) \widetilde{h}(n) f(n + (k-1)r) \prod_{j=1}^{k-2} h_j(n + jr)\big| \\ & = &  \big|\E_{n \in \Z/N\Z} \E_{r \in \Z/N\Z} \sum_{m \in \Z/N\Z} \widehat{\chi}(m) \widetilde{h}(n) e_N(mn) h_1(n + r) e_N(-m(n+r)) \times \\ & & \qquad \qquad \qquad \qquad \qquad \qquad \times \prod_{j=2}^{k-2} h_j(n + jr)\cdot f(n + (k-1)r)\big| \\ & \leq & 2 \sum_{m \in \Z/N\Z} |\widehat{\chi}(m)| \Vert f \Vert_{U^{k-1}(\Z/N\Z)} \leq \frac{200}{\eps} \Vert f \Vert_{U^{k-1}(\Z/N\Z)}.
\end{eqnarray*}
\ni Now we substitute this into \eqref{eq10.105} and average over $z_1,\dots,z_{k-2}$ (picking these $k-2$ elements uniformly according to the Haar measure $\mathbb{P}$ on $G/\Gamma$). This yields
\begin{equation}\label{eq10.106} \big|\E_{|n| \leq N/2 - \eps k N} \E_r \chi(r) F(T_g^{n + (k-1)r} x) f(n + (k-1)r)\big| \leq \frac{200}{\eps \mathbb{P}(\pi(V'))^{k-2}} \Vert f \Vert_{U^{k-1}(\Z/N\Z)} + O(\eps).\end{equation}
Write $G(n) := F(T_g^n x) f(n)$. Then $\Vert G \Vert_{\infty} \leq 1$, and so we have
\begin{eqnarray*}
\E_{|n| \leq N/2 - \eps k N}\E_r \chi(r)  G(n + (k-1)r) & = & \E_r \chi(r) \E_{-N/2 + \eps k N - (k-1)r \leq n' \leq N/2 - \eps k N - (k-1)r} G(n') \\ & = & \E_r \chi(r) \big( \E_{|n| \leq N/2} G(n') + O(\eps) \big) \\ & = & \E_{|n| \leq N/2} G(n) + O(\eps).\end{eqnarray*}
Comparing this with \eqref{eq10.106} gives
\[ \big| \E_{|n| \leq N/2} f(n) F(T_g^n x) \big| \leq \frac{200}{\eps \mathbb{P}(\pi(V'))^{k-2}} \Vert f \Vert_{U^{k-1}(\Z/N\Z)} + O(\eps),\]
which implies the result if $\eps = c \eta$ for $c = c_k$ sufficiently small.\endproof\vs

\ni\textit{Remark.} An alternative way to obtain this lemma is to establish that the $(k-2)$-step 
nilsequence $n \mapsto F(T_g^n x)$ can be approximated to high accuracy by
a function which is uniformly almost periodic of order $k-2$ in the sense of \cite{tao:ergodic}; this approach has the advantage of not requiring
an explicit algebraic constraint such as that given in Lemma \ref{constraint-lem}, but we do not pursue it here.  
This approach corresponds closely to the observation 
that a $k-2$-step nilflow can be constructed as a tower of $k-2$ compact extensions of the trivial measure-preserving system, see \cite{furstenberg-weiss,host-kra2,ziegler02} for further discussion.  \vs

\ni Proposition \ref{tgnx} shows (essentially) that the basic $(k-2)$-step nilsequences 
form ``obstructions to quadratic uniformity'', in the
sense that functions $f: \Z/N\Z \to \D$ which have a large inner product with such functions cannot have small $U^{k-1}(\Z/N\Z)$ norm.  The remarkable result of
Host and Kra \cite{host-kra2} asserts, roughly speaking, that these are in fact the \emph{only} obstructions to having small $U^{k-1}$ norm.  More precisely, they work in the infinitary setting of arbitrary measure-preserving systems (as opposed to the shift on $\Z/N\Z$) and show that this system contains as an invariant factor an inverse limit of $(k-2)$-step nilflows, such that the $U^{k-1}(T)$ norm vanishes on
the orthogonal complement of this inverse limit.  In particular, this inverse limit is a characteristic factor for the $U^{k-1}(T)$ norm, and for all quantities controlled by this norm, including the $k$-term recurrence expressions appearing for instance in the Furstenberg recurrence theorem; this fact is crucial in establishing the convergence of these recurrence expressions.
We remark that the work of Ziegler \cite{ziegler02,ziegler} achieves a very similar result, but avoids use of the $U^{k-1}(T)$ norm and obtains a characteristic factor (and convergence results) for the recurrence expressions directly.  Also, the subsequent work of Bergelson, Host, and Kra \cite{bhk} gives a further discussion of the connection between the $U^{k-1}(T)$ norm and $(k-2)$-step nilsequences.\vs

\ni We now use Theorem \ref{u3-quad} to obtain a finitary (and reasonably quantitiative) version of the Host-Kra theorem in the case $k=4$, with very explicit nilsequences; in fact, they will be none other than the elementary $2$-step nilsequences defined earlier.

\begin{theorem}[Inverse theorem for $U^3(\Z/N\Z)$, elementary nilsequence version]\label{nilsequence-thm}  Let

\ni $N > 2$ be a prime, let $0 < \eta \leq 1$ be sufficiently small, and suppose that $f : \Z/N\Z \rightarrow \D$ is a function with $\| f \|_{U^3(\Z/N\Z)} \geq \eta$. Then there exists an elementary $2$-step nilsystem $G/\Gamma$ of dimension $\leq \eta^{-C}$, an $\eta^{-C}$-Lipschitz function $F : G/\Gamma \rightarrow \D$, and elements $g \in G$, $x_0 \in G/\Gamma$ and $h \in \Z/N\Z$ such that 
\[ |\E_{n \in \Z/N\Z} (T^h f \overline{F_{N,g,x_0}}) | \geq \exp(-\eta^{-C}).\]
Here, we define $F_{N,g,x_0} : \Z/N\Z \rightarrow \C$ by
\[ F_{N,g,x_0}(n) := F(T_g^n x_0) \; \mbox{for all} \; -N/2 < n < N/2.\]
\end{theorem}
\ni\textit{Remarks.} The function $F_{N,g,x_0}$ is just a $2$-step nilsequence, adapted to $\Z/N\Z$. The analogous theorem for $U^2(\Z/N\Z)$ is a trivial consequence of Proposition \ref{u2-inverse}; the only linear nilfunction that
needs to be considered is the function $F(x) := e(x)$ on the unit circle $\R/\Z$ from Example \ref{nil1}, with $h = 0$ and $x_0 = 0$.  
A modification of Example \ref{fw-ex}
can be used to show that in formulating this theorem we must take into account Example \ref{nil3}; the other two fundamental $2$-step nilsystems are in fact embedded inside this one and one could have dispensed with them altogether, but we have kept them for expository purposes.  One could also easily eliminate the r\^ole of $x_0$ (which is harmless anyway, since it ranges over a compact set) and of the shift $h$, but the parameter $g$ ranges over a genuinely non-compact set and cannot be eliminated from this theorem (this can be seen even in the linear case; the frequency $\xi$ in Proposition \ref{u2-inverse} is not restricted to a bounded set of values independently of $N$).  \vs

\ni\textit{Proof.} 
Applying Theorem \ref{u3-quad} (and Lemma \ref{make-regular}), we obtain a set $S \subseteq \Z/N\Z$ with $d := |S| \leq \eta^{-C}$, a regular Bohr set $B = B(S,\rho)$ with $\rho \geq \eta^{C}$, and a bracket quadratic 
\begin{equation}\label{phi-expand}\phi(n) := \sum_{\xi,\xi' \in S} a_{\xi,\xi'} \{ \xi \cdot n\} \{ \xi' \cdot n\} + \sum_{\xi \in S} a_{\xi} \{\xi \cdot n\}\end{equation} with $\Freq(\phi) \subseteq S$
such that 
$$ |\E_{n \in B}(T^h f(n) e(-\phi(n)) )| \geq \eta^C$$
and thus
\begin{equation}\label{eq10.107}|\E_{n \in \Z/N\Z}( T^{h} f(n) e(-\phi(n)) 1_B(n) )| \geq \eta^C \E 1_B.\end{equation}
Now let $\eps := \eta^C/400d$, and let
$\chi: \R/\Z \to [0,1]$ be a continuous function such that $\chi(x)=1$ when $|x| < \rho(1 - \eps)$ and $\chi(x) = 0$ when $1/2 \geq |x| > \rho(1 + \eps)$.
Consider the function\footnote{The large power of $\chi$ here is so that we can distribute the cutoff $\chi$ among various factors later.}
\[ \prod_{\xi \in S} \chi^{1 + 4d}(\xi \cdot n) - 1_B(n).\]
It is supported on $B(S,\rho(1 + \eps)) \setminus B(S,\rho(1 - \eps))$, which by the regularity of $B$ has cardinality no more than $200d\eps |B|$. Thus \eqref{eq10.107} implies that 
\[  \big|\E_n T^hf(n) e(-\phi(n)) \prod_{\xi \in S} \chi^{1 + 4d}(\xi \cdot n)\big| \geq (\eta^C - 200d\eps )\E 1_B \geq \half \eta^C \E 1_B \geq \exp(-\eta^{-C}),\]
the latter inequality being a consequence of Lemma \ref{bohr-bound}.\vs

\ni Expanding out $\phi$ as in \eqref{phi-expand}, we see that our task is to show that the function
\begin{equation}\label{eq10.111} n \mapsto \bigg(\prod_{\xi \in S} \chi(\xi \cdot n) e(a_\xi \{ \xi \cdot n \} ) \bigg) 
\bigg(\prod_{\xi, \xi' \in S} \chi^2(\xi \cdot n) \chi^2(\xi' \cdot n)
 e( a_{\xi,\xi'} \{ \xi \cdot n \} \{\xi' \cdot n\} \bigg)\end{equation}
is an elementary $2$-step nilsequence, for which it is enough to handle each of the functions in the product separately as in the following lemma.

\begin{lemma}\label{construct-lemma}
Each of the individual functions
\begin{equation}\label{linear-xi}
n \mapsto \chi(\xi \cdot n) e(a_\xi \{ \xi \cdot n \} )
\end{equation}
and 
\begin{equation}\label{quadratic-xi}
 n \mapsto \chi^2(\xi \cdot n) \chi^2(\xi' \cdot n) e( a_{\xi,\xi'} \{ \xi \cdot n \} \{\xi' \cdot n\} )
 \end{equation}
can be written as an elementary $2$-step nilsequence with dimension no more than $9$ and Lipschitz constant at most $50$.\end{lemma}

\ni\textit{Proof.} We begin by considering the functions \eqref{linear-xi}, which are easier (corresponding to linear nilcharacters rather than quadratic ones).
Split $a_\xi = q + s$, where $q$ is an integer and $|s| \leq 1/2$, and
observe that $e( q \{ \xi \cdot n \} ) = e(q\xi n)$ if we identify
$\Z/N\Z$ and $\widehat{\Z/N\Z}$ with the integers from $-N/2$ to $N/2$.  Thus the function \eqref{linear-xi} takes the form
$$ n \mapsto \chi(\xi n / N) e( s \{ \xi n / N \} ) e(q\xi n / N).$$
This function may be identified as the elementary nilsequence $F_{N,(\xi/N,q\xi/N),0}$, where the underlying nilmanifold 
$G/\Gamma$ is the direct sum of two copies of the unit circle shift (i.e. it is the torus $(\R/\Z)^2$) and
$F: (\R/\Z)^2 \to \C$ is the function
$$ F(x,y) := \chi(x) e(sx) e(y)$$
where we identify $x \in \R/\Z$ with a real number from $-1/2$ to $1/2$ in the usual manner. It is not hard to check that $F$ is $50$-Lipschitz.\vs

\ni Now consider the functions \eqref{quadratic-xi}.  We split $a_{\xi,\xi'}$ as $q + s$, much as before, so that \eqref{quadratic-xi} becomes
$$ n \mapsto \chi( \alpha n)^2 \chi(\gamma n)^2 e(s \{\alpha n \} \{\gamma n \}) 
 e( q \{ \alpha n \} \{ \gamma n \} ),$$
where $\gamma := \xi/N$ and $\gamma := \xi'/N$.   Observe that
since $\{\alpha n\} = \alpha n - [\alpha n]$ and $\{ \gamma n \} = \gamma n - [\gamma n]$, we have the identity
$$ q \{ \alpha n \} \{\gamma n\} = q\alpha\gamma n^2 - q\alpha n [\gamma n] - q\gamma n [\alpha n] + q[\alpha n] [\gamma n].$$
The last term is an integer, and hence
$$ e(q \{ \alpha n \} \{\gamma n\}) = e(q\alpha\gamma n^2) e(-q\alpha n [\gamma n]) e(-q\gamma n [\alpha n]).$$
Thus it suffices to exhibit the three functions
\begin{align*}
n &\mapsto \chi(\alpha n) \chi(\gamma n) e(s \{\alpha n \} \{\gamma n \})\\
n &\mapsto e(q\alpha\gamma n^2) \\
n &\mapsto \chi( \gamma n) e(-q\alpha n [\gamma n]) \; , \;  
n \mapsto \chi( \alpha n) e(-q\gamma n [\alpha n])
\end{align*}
as elementary $2$-step nilsequences (note that the last two functions are essentially the same).\vs

\ni The first function can be obtained from a direct sum of two copies of the unit circle shift (Example \ref{nil1}) by repeating
the analysis of \eqref{linear-xi}, with $F$ now defined by $F(x,y) := \chi(x) \chi(y) e(sxy)$
when $-1/2 < x,y \leq 1/2$.\vs

\ni The second function can easily be obtained from the skew shift (Example \ref{nil2}), by writing
$$ q\alpha\gamma n^2 = - q\alpha\gamma n + 2q\alpha\gamma \frac{n(n+1)}{2},$$
and then taking $F(x,y) := e(y)$, $x_0 = (0,0)$, and
$$ g := \left( \begin{array}{lll}
1 & 1  & -q\alpha\gamma \\
0 & 1  & 2q\alpha\gamma \\
0 & 0  & 1
\end{array} \right).$$
\ni Finally let us consider the third function.  We write
$$ e(-q\gamma n [\alpha n]) = e( \half n(n+1) \alpha q\gamma - [n\alpha] n q\gamma ) e( - \half n(n+1) \alpha q\gamma ).$$
The second factor can be generated using the skew shift as before\footnote{Indeed we could simply have rewritten our factorization of $e(q \{\alpha n \} \{\gamma n\})$ to incorporate these factors (and a linear phase correction), so as to then dispense with the second factor and the skew shift altogether.  However, we have left this example in here to emphasize that purely quadratic phase functions such as $e(q\alpha\gamma n^2)$ are indeed examples of nilcharacters.}.  We are thus left with
$$ \chi(\alpha n) e( \frac{n(n+1)}{2} \alpha q\gamma - [n\alpha] n q\gamma ).$$
But this can be generated from the Heisenberg shift (Example \ref{nil3}) with $x_0 = (0,0,0)$, 
$F(x,y,z) := \chi(x) e(y)$ on $[-1/2,1/2] \times (\R/\Z)^2$, and 
$$ g := \left( \begin{array}{lll}
1 & 0  & q\gamma \\
0 & 1  & \alpha \\
0 & 0  & 1
\end{array} \right).$$
It is easy to check that all of the functions $F$ used in these constructions are $2\pi$-Lipschitz, a bound which together with Lemma \ref{prod-lipschitz} completes the proof of the lemma.\endproof\vs

\ni Lemma \ref{construct-lemma}, together with another application of Lemma \ref{prod-lipschitz}, confirms that the function \eqref{eq10.111} is an elementary $2$-step nilfunction with dimension at most $18d^2$ and Lischitz constant no more than $100 d^2$. This completes the proof of Theorem \ref{nilsequence-thm}.\endproof\vs

\ni\textit{Remark.} A pleasant reformulation of Theorem \ref{nilsequence-thm} may be obtained by considering the $\delta$-atoms of $G/\Gamma$. Suppose that $F : G/\Gamma \rightarrow \D$ is $K$-Lipschitz. Let $(A_{\omega})_{\omega \in \Omega}$ be the $\delta$-atoms of $G/\Gamma$, pick arbitrary points $x_{\omega} \in A_{\omega}$ for each $\omega$, and write 
\[ \widetilde F(n) := \sum_{\omega \in \Omega} F(x_{\omega}) 1_{A_{\omega}}(n).\]
Since $F$ is $K$-Lipschitz we clearly have the bound $\| F - \widetilde F\| \leq K\delta$. Taking $\delta \leq \half K^{-1}\exp(-\eta^{-C})$ we may replace the conclusion of Theorem \ref{nilsequence-thm} by
\[ |\E_{n \in \Z/N\Z} (T^h f \overline{\widetilde{F}_{N,g,x_0}}) | \geq \half\exp(-\eta^{-C}).\] Removing the sum over the atoms (of which there are at most $\delta^{-\eta^{-C}}$) by the pigeonhole principle this implies that 
\[ |\E_{n \in \Z/N\Z} (T^h f \overline{(1_{A_{\omega}})_{N,g,x_0}}) | \geq \exp(-2\eta^{-2C}).\]
That is, if $\Vert f \Vert_{U^3}$ is large then $f$ correlates with the set of return times of a $2$-step nilsequence to an atom.\vs

\ni\textit{Remark.} The space of quadratic nilsequences forms an algebra, being closed under multiplication, addition, subtraction, and
conjugation.  This allows one to employ an ``energy incrementation'' argument of the type used in \cite[\S 7]{green-tao-primes} in order to decompose an arbitrary bounded function on $\Z/N\Z$ as the sum of a
bounded function with $U^3$ norm smaller than some specified $\eta$, plus a $2$-step nilsequence with dimension and Lipschitz constant controlled by functions of $\eta$. In the ergodic theory setting an extremely similar decomposition was obtained in \cite{bhk}.  Informally speaking, the quadratic nilsequences form a \emph{characteristic factor} for the $U^3$ norm, and hence for any expression controlled by that norm, and
thus many questions involving such expressions can be reduced to questions concerning $2$-step nilsequences.\vs

\ni It is perhaps of interest to briefly discuss a decomposition of this type for the $U^2$ norm, where the availability of harmonic analysis allows one to proceed more directly. If $f : \Z/N\Z \rightarrow \D$ is a function then we write $R:= \{ r : |\widehat{f}(r)| \geq \eps_1\}$ for some suitable $\eps_1$, and define $\beta(x)$, $\E \beta = 1$, to be a normalized and suitably smoothed version of $1_B(x)$, where $B := B(R,\eps_2)$. We then
decompose
\begin{equation}\label{eq10.115} f(x) = f \ast \beta + (f - f \ast \beta).\end{equation} It is easy to check that $\Vert f - f \ast \beta \Vert_{U^2}$ is small, but to write $f \ast \beta$ as a $1$-step nilsequence it must be modified slightly. To do this, write
\begin{equation}\label{eq10.116} f \ast \beta(x) = \sum_{r \in \overline{R}} \widehat{f}(r)\widehat{\beta}(r) e(-rx/N) + \sum_{r \notin \overline{R}} \widehat{f}(r)\widehat{\beta}(r) e(-rx/N),\end{equation}
where $\overline{R} := \{ m_1 r_1 + \dots + m_d r_d | r_j \in R , |m_j| \leq M\}$ for some $M$. Now if $\beta$ is sufficiently smoothed and if $M$ is sufficiently large then \[ \sum_{r \notin \overline{R}} |\widehat{\beta}(r)| \leq \eps_3,\] and so the second term in \eqref{eq10.116} is bounded by $\eps_3$, and in particular has small $U^2$ norm. The first term can be written as a $1$-step nilfunction, the underlying nilmanifold being $(\R/\Z)^d$ and the rotation $T$ being $(x_1,\dots,x_d) \mapsto (x_1 + r_1/N,\dots,x_d + r_d/N)$. See \cite{green-roth-primes} for an application of such a decomposition (there the language of nilsystems and ergodic theory did not feature, and the simpler decomposition \eqref{eq10.115} was used).

\section{Future prospects and open questions}\label{sec13}

\ni It is natural to ask whether there are inverse theorems for the higher $U^k(G)$-norms, $k \geq 4$, which generalize Theorems \ref{mainthm2}, \ref{u3-quad} and \ref{nilsequence-thm}. We are certain that the answer to this question is ``yes''. It is easy to guess at the correct generalization of Theorem \ref{nilsequence-thm}, which should simply involve replacing $2$-step nilmanifolds by $(k-2)$-step ones. Guessing at the generalization of Theorem \ref{mainthm2} is a bit harder. We suspect that the correct objects to consider for the $U^4(G)$-norm are of the form $1_Q(x) e(-\lambda(x))$, where now $Q := \{ x : \psi_1(x),\dots,\psi_d(x) \approx 1 \}$ is a \textit{quadratic Bohr set}, each $\psi_j$ being of the form $\psi_j(x) = 1_{B_j}(x) e(-\phi_j(x))$ appearing in Theorem \ref{mainthm2}. The phase $\lambda : Q \rightarrow \R/\Z$ is now cubic. It is easy to guess how functions appropriate for the $U^k(G)$ norm may be constructed inductively.\vs

\ni We think it likely that most of the ingredients necessary to prove such inverse theorems may be found in \cite{gowers-long-aps}, and we intend to pursue this direction. The same major difficulty that Gowers encountered in dealing with the $U^k(G)$ norm for $k \geq 4$ is also present here. 
Suppose that $f$ has large $U^4(G)$ norm. This means that $T^h f \overline{f}$ has large $U^3(G)$ norm for many values of $h$.  Applying Theorem \ref{mainthm2} 
we obtain a Bohr set $B_h$ for each of these $h$, such that $T^h \overline{h}$ has large quadratic bias on several shifts of this Bohr set $B_h$.  
The problem is that we do not, \textit{a priori}, have any control on how the Bohr set $B_h$ depends on $h$.\vs

\ni Another interesting issue is that of obtaining better bounds in Theorems \ref{mainthm1} and Theorem \ref{mainthm2}. It is quite possible that the codimension of $W$ in Theorem \ref{mainthm1} can be taken to be $O(\log(1/\eta))$ rather than $O(\eta^{-C})$. 
This would give bounds of the form
\begin{equation}\label{eq10.117} \Vert f \Vert_{U^3(\Ffiven)} \leq \Vert f \Vert_{u^3(\Ffiven)} \leq \Vert f \Vert_{U^3(\Ffiven)}^c\end{equation}
for some absolute constant $c$. Such a bound would be a consequence of the \textit{Polynomial Freiman-Ruzsa Conjecture} (PFR), which is discussed in detail in \cite{green-fin-field}, together with some mild adjustments to the arguments of \S \ref{finite-sec}. We refer to the statement \eqref{eq10.117} as the \textit{Polynomial Gowers Inverse Conjecture} (PGI) for $\Ffiven$.\vs

\ni It would be nice to have a version of \eqref{eq10.117} in a general $G$. What we mean by this is a statement of the form
\[ \Vert f \Vert_{U^3(G)} \geq \eta \quad \Longrightarrow \quad |\E f(x) 1_B(x) e(-\phi(x))| \geq c(\eta),\]
$B = B(S,\rho)$, which may be reversed with only polynomial losses in the constants, that is to say
\[ |\E f(x) 1_B(x) e(-\phi(x))| \geq c(\eta) \quad \Longrightarrow \quad \Vert f \Vert_{U^3(G)} \geq \eta^C.\]
This would seem to require that we can take $|S| = O(\log(1/\eta))$ and $\rho$ greater than some absolute constant\footnote{even then one would need to replace $1_B$ with something smoother, to avoid the losses in the argument at the beginning of \S \ref{ggc}.}. The methods of this paper seem to fall a long way short of proving such a statement. Even if one had an appropriate analogue of PFR (which might take the form of a stronger version of Lemma \ref{bog-lemma} in which the size of $S$ is logarithmic in $\delta$), we would have to find a way to avoid repeatedly passing to smaller Bohr sets as in \S \ref{ggc}. Each such passage causes too much degradation in $\rho$.\vs

\ni Let us conclude this section by remarking that working out how to drop the restriction $(|G|,6) = 1$ in Theorem \ref{mainthm2} would be a diverting exercise at least, though we cannot think of any applications. The case $G = \mathbb{F}_2^n$ probably captures the essence of the problem\footnote{The authors have recently learnt that Samorodnitsky \cite{samorod} has resolved this issue.}.

\section{Appendix: algebraic constraints on nilmanifolds}\label{appendix-sec}

\ni In this section\footnote{The authors are indebted to Sasha Liebman and Tamar Ziegler for conversations which were very helpful in preparing this appendix.} we prove Lemma \ref{constraint-lem}  For some further discussion of issues related to such constraints, see \cite{bhk,furstsurvey1,furstsurvey2,furstsurvey3,ziegler02}.\vs

\ni By replacing $g, x$ by $g^r, T_g^n x$ respectively our task is to demonstrate, given any $k-2$-step nilmanifold $G/\Gamma$, the existence of a 
compact set $\Sigma \subseteq  (G/\Gamma)^{k-1}$
and a continuous map $P: \Sigma \to G/\Gamma$ such that
\begin{equation}\label{px-alt}
(x, T_g x, \ldots, T_g^{k-2} x) \in \Sigma; \quad 
P(x, T_gx, \ldots, T_g^{k-2} x) = T_g^{n+(k-1)r} x \hbox{ for all } g \in G, x \in G/\Gamma.
\end{equation}

\ni\textit{Remark.} One can prove \eqref{px-alt} by direct algebraic computation in the cases $k=3,4$.  Indeed, when $k=3$ the $1$-step nilpotent group $G$ is abelian, 
as is the subgroup $\Gamma$, so $G/\Gamma$ is also a group (indeed it is a torus).  One can then 
take $\Sigma = (G/\Gamma)^2$ and $P(a,b) := a (a^{-1} b)^2$.  When $k=4$, so that $G$ is a $2$-step nilpotent group,
things are a little more complicated.  One needs to take
$$ \Sigma := \{ (x_0 \Gamma, x_1 \Gamma, x_2 \Gamma): x_2 \in x_0 (x_0^{-1} x_1)^2 G_2 \};$$
this reflects the fact that $G/G_2$ is a $1$-step nilpotent group and thus obeys the $k=3$ constraints.  Note that $G_2$ commutes with all elements of $G$ and is thus easy to quotient out.  One can then define $P: \Sigma \to G/\Gamma$
by setting
$$ P( a\Gamma, b\Gamma, c\Gamma) := a (a^{-1} b)^3 ((a^{-1} b)^{-2} a^{-1} c)^3 \Gamma \hbox{ whenever } (a\Gamma, b\Gamma, c\Gamma) \in \Sigma;$$
one can verify with some effort (using of course the fact that all commutators lie in $G_2$, which commute with all elements of $G$)
that this function is well-defined, continuous on $\Sigma$, and obeys \eqref{px-alt}.  Unfortunately in the $k>4$ case it seems that the function $P$
is significantly messier, and in particular requires choosing a partial inverse for projection maps $G_j \mapsto G_j/\Gamma$ for all $j < k-2$,
which in general cannot be done canonically when $j \geq 2$.\vs

\ni To prove \eqref{px-alt} in the general case, we first need some notation.

\begin{definition}[Continuous right invertibility]
Let $M, N$ be compact spaces, let $\pi: M \to N$ be a continuous map, and let $\Sigma \subseteq M$.  We say $\pi$ is \emph{continuously 
right-invertible on} $\Sigma$
if for every $w \in \overline{\pi(\Sigma)}$ there exists there exists a neighbourhood $V_{w} \subseteq  N$ of $w$ and a continuous map $\pi^{-1}_{w}: V_{w} \to M$ such that
$\pi^{-1}_{w} \circ \pi$ is the identity on $\Sigma \cap \pi^{-1}(V_{w})$. 
\end{definition}

\begin{lemma}  Let $G/\Gamma$ be a $(k-2)$-step nilmanifold, and let $\pi: (x_0,\ldots,x_k) \to (x_0,\ldots,x_{k-1})$ be
the canonical projection from $(G/\Gamma)^k$ to $(G/\Gamma)^{k-1}$. Then $\pi$ is continuously right-invertible on the set
$$ A := \{ (x, T_g x, \ldots, T_g^{k-1} x, T_g^k x): x \in G/\Gamma, g \in G \}.$$
\end{lemma}

\ni The existence of a constraint \eqref{px-alt} then follows by taking $\Sigma$ to be the closure of $\pi(A)$ and using a compactness argument (exploiting the fact that $\pi(A)$ is dense in the compact set $\Sigma$) to glue the various local right-inverses $\pi^{-1}_w$ together.\vs

\begin{proof} 
For any $0 \leq i \leq k-1$, we define the \emph{Hall-Petresco groups} $HP_{k,i} \subseteq  G^k$ to be the sets
$$ HP_{k,i} := \{ ( x_i^{\binom{n}{i}} \ldots x_{k-2}^{\binom{n}{k-2}})_{0 \leq n \leq k-1} : x_i \in G_i, \ldots, x_{k-2} \in G_{k-2} \},$$
where we have the conventions that $G_0 := G_1 = G$, and that $\binom{n}{j} = 0$ if $j > n$.
Thus for instance when $k=4$ we have
\begin{align*}
HP_{4,3} &= \{ (1,1,1,1) \}\\
HP_{4,2} &= \{ (1,1, x_2, x_2^3): x_2 \in G_2 \} \\
HP_{4,1} &= \{ (1, x_1, x_1^2 x_2, x_1^3 x_2^3): x_1 \in G_1, x_2 \in G_2 \} \\
HP_{4,0} &= \{ (x_0, x_0 x_1, x_0 x_1^2 x_2, x_0 x_1^3 x_2^3): x_0 \in G_0, x_1 \in G_1, x_2 \in G_2 \}.
\end{align*}
It is well known (see \cite{leibman}) that the $HP_{k,i}$ are all subgroups of $G^k$, and we also have the nesting
$HP_{k,i+1} \subseteq  HP_{k,i}$ for all $0 \leq i < k-1$.  \vs

\ni Observe that $\Gamma^k$ is a subgroup of $G^k$, so we may form the quotient space $G^k/\Gamma^k$, which we identify with the compact
manifold $(G/\Gamma)^k$.  Inside this space we have the submanifolds $HP_{k,i} / \Gamma^k$ for all $0 \leq i \leq k-1$.  Observe that if $x \in G/\Gamma$ and $g \in G$, and $y \in x\Gamma \subseteq G$ is any representative of $x$ in $G$, then
\begin{align*}
(x,T_g x, \ldots, T_g^k x) \Gamma^k &= (y, g y, \ldots, g^k y) \Gamma^k \\
&= (y,\ldots,y) (1,g,\ldots,g^k) \Gamma^k
\subseteq  HP_{k,0} HP_{k,0} \Gamma^k \\
\subseteq  HP_{k,0} \Gamma^k
\end{align*}
and hence
$$ (x,T_g x, \ldots, T_g^k x) \in HP_{k,0} / \Gamma^k.$$
It thus suffices to show that $\pi$ is right-invertible on $HP_{k,0} / \Gamma^k$.\vs

\ni We shall show inductively, by backwards induction on $i$, that $\pi$ is continuously right-invertible on $HP_{k,i} / \Gamma^k$
for all $0 \leq i \leq k-1$.  The case $i=k-1$ is trivial since $HP_{k,i} / \Gamma^k$ is just a point.  Now suppose
inductively that $0 \leq i < k-1$ and that $\pi$ was already shown to be continuously right-invertible over $HP_{k,i+1} / \Gamma^k$.\vs

\ni Let $(z_0,\ldots,z_{k-1}) \in \overline{\pi(HP_{k,i}/\Gamma^k)}$.  Observe that the first $i$ coefficients of
$z$ must be the origin $O \in G/\Gamma$, defined as the image of the identity $1 \in G$, and the coefficient $z_i$ lies
in the closed manifold $G_i/\Gamma$.\vs

\ni The projection map $\pi_i: G_i \mapsto G_i / \Gamma$ is continuous and surjective from the manifold $G_i$ to the manifold $G_i/\Gamma$, which is
a sub-manifold of $G/\Gamma$.  Thus we may find a continuous function $f: V_{z_i} \to G_i$ defined on a neighbourhood $V_{z_i} \subseteq G/\Gamma$ of $z_i$ such that $\pi_i \circ f$ is the identity on $V_{z_i} \cap (G_i/\Gamma)$.  \vs

\ni We need to right-invert $\pi$ on $HP_{k,i}/\Gamma^k$ in a neighbourhood of $\pi(z)$.  To this end, let
$x := (x_1,\ldots,x_k) \in HP_{k,i}/\Gamma^k$ be such that $\pi(x)$ be close to $\pi(z)$; in particular we may take $x_i \in V_{z_i}$.  
As before we have $x_n = O$ for $n < i$,
and $x_n \in G_i/\Gamma$ for all $n \geq i$.  Thus $x_i \in V_{z_i} \cap (G_i/\Gamma)$ and hence $\pi_i \circ f(x_i) = x_i$, or
in other words $x_i = f(x_i) \Gamma$.  Now let $F(x_i) \in HP_{k,i}$ be the group element
$$ F(x_i) :=  ( ( f(x_i)^{\binom{n}{i}})_{0 \leq n \leq k-1} ),$$
and observe that this depends continuously on $x_i$, and hence on $\pi(x)$, if $\pi(x)$ lies in a neighborhood of $\pi(z)$.\vs

\ni On the other hand, since $x \in HP_{k,i}/\Gamma^k$, there exists
$g = (1,\ldots,1,g_i,\ldots,g_k) \in HP_{k,i}$ such that $g \Gamma^k = x$; in particular, $g_i \in G_i$ and $g_i \in x_i \Gamma$.
Since $f(x_i) \in G_i$ and $f(x_i) \in x_i \Gamma$, we conclude that $f(x_i)^{-1} g_i$ lies in both $G_i$ and in $\Gamma$.  Thus if
we let $\tilde g \in HP_{k,i}$ be the group element
$$\tilde g = ( ( (f(x_i)^{-1} g_i)^{\binom{n}{i}})_{0 \leq n \leq k-1} )$$
then $\tilde g$ lies in both $HP_{k,i}$ and $\Gamma^k$.  Thus we can factorize
$$ g = F(x_i) h \tilde g$$
where $h$ lies in $HP_{k,i}$, and also has $i^{th}$ component equal to the identity.  Thus $h$ in fact lies in $HP_{k,i+1}$.  Multiplying
on the right by $\Gamma^k$, we conclude that
$$ x = g \Gamma^k = F(x_i) h \tilde g \Gamma^k = F(x_i) h \Gamma^k$$
and hence
$$ F(x_i) x = h\Gamma^k \in HP_{k,i+1}/\Gamma^k.$$
Since $\pi(x)$ is close to $\pi(z)$, and $F(x_i)$ depends continuously on $\pi(x)$, we see that $\pi(F(x_i) x)$ is close to $\pi(F(z_i) z)$.
In particular, by the induction hypothesis we can find a continuous map $\pi^{-1}_{\pi(F(z_i) z)}$ mapping a neighborhood $V_{F(z_i) z} \subseteq (G/\Gamma)^k$ of
$\pi(F(z_i) z)$ to $HP_{k,i+1}/\Gamma^k$ which is a local right-inverse of $\pi$ on $HP_{k,i+1}/\Gamma^k \cap \pi^{-1}(V_{F(z_i) z})$.
Thus we have
$$ F(x_i) x = h \Gamma^k = \pi^{-1}_{F(z_i) z}( \pi(F(x_i) x) )$$
and hence
$$ x = F(x_i)^{-1} \pi^{-1}_{F(z_i) z}( \pi(F(x_i) x) ).$$
Observe that $\pi(F(x_i) x) = \tilde \pi(F(x_i)) \pi(x)$, where $\tilde \pi: G^k \to G^{k-1}$ is the canonical projection.  Since $x_i$ of
course depends continuously on $\pi(x)$, the
right-hand side then depends continuously on $\pi(x)$ when $\pi(x)$ lies in a sufficiently small neighbourhood of $\pi(z)$. We have
achieved a right-inverse for $\pi$ on $HP_{k,i}/\Gamma^k$ in a neighborhood of $\pi(z)$, thus closing the induction.
\end{proof}

\providecommand{\bysame}{\leavevmode\hbox to3em{\hrulefill}\thinspace}

\end{document}